\documentclass[11pt]{amsart}

\date{\today}

\usepackage{amsmath, amsthm, amscd, amsfonts, amssymb, color}
\usepackage[bookmarksnumbered, colorlinks]{hyperref}

\theoremstyle{theorem}
\newtheorem{theorem}{Theorem}[section]
\newtheorem{cor}[theorem]{Corollary}
\newtheorem{prop}[theorem]{Proposition}
\newtheorem{lemma}[theorem]{Lemma}
\theoremstyle{remark}
\newtheorem{remark}[theorem]{Remark}
\newtheorem*{cremark*}{Concluding remarks}
\newtheorem{example}[theorem]{Example}
\newtheorem{definition}[theorem]{Definition}

\DeclareMathOperator*{\esssup}{ess\,sup}

\date{\today}
\title{Weak frames in Hilbert $C^*$-modules with application in Gabor analysis}
\author[D. Baki\' c]
{Damir Baki\' c$^{*}$}
\thanks{$^{*}$ This work has been fully supported by the Croatian Science Foundation under the project IP-2016-06-1046.}

\address{Department of Mathematics, University of Zagreb,
Bijeni\v cka cesta 30, 10000 Zagreb, Croatia.}
\email{bakic@math.hr}

\begin{document}

\begin{abstract}
In the first part of the paper we describe the dual $\ell^2(\textsf{A})^{\prime}$ of the standard Hilbert $C^*$-module $\ell^2(\textsf{A})$ over an arbitrary (not necessarily unital) $C^*$-algebra \textsf{A}.  When $\textsf{A}$ is a von Neumann algebra, this enables us to construct explicitly a self-dual Hilbert $\textsf{A}$-module $\ell^2_{\text{strong}}(\textsf{A})$ that is isometrically isomorphic to $\ell^2(\textsf{A})^{\prime}$, which contains $\ell^2(\textsf{A})$, and whose $\textsf{A}$-valued inner product extends the original inner product on $\ell^2(\textsf{A})$. This serves as a concrete realization of a general construction for Hilbert $C^*$-modules over von Neumann  algebras introduced by W. Paschke.

Then we introduce a concept of a weak Bessel sequence and a weak frame in Hilbert $C^*$-modules over von Neumann algebras. The dual $\ell^2(\textsf{A})^{\prime}$ is recognized as a suitable target space for the analysis operator. We describe fundamental properties of weak frames such as the correspondence with surjective adjointable operators, the canonical dual, the reconstruction formula, etc; first for self-dual modules and then, working in the dual, for general modules.

In the last part of the paper we describe a class of Hilbert $C^*$-modules over $L^{\infty}(I)$, where $I$ is a bounded interval on the real line, that appear naturally in connection with Gabor (i.e.~Weyl-Heisenberg) systems. We then demonstrate that Gabor Bessel systems and Gabor frames in $L^2(\Bbb R)$ are in a bijective correspondence with weak Bessel systems and weak frames of translates by $a$ in these modules over $L^{\infty}[0,\frac{1}{b}]$, where $a,b>0$ are the lattice parameters. In this setting some well known results on Gabor systems are discussed and some new are obtained.
\end{abstract}

\subjclass[2010]{42C15, 46L08, 42C40,}

\keywords{Hilbert $C^*$-module, von Neumann algebra, frame, Bessel sequence, Gabor frame}
\maketitle

\vspace{0.2in}

\section{Introduction}

\vspace{.1in}

Frame theory for Hilbert $C^*$-modules is now about two decades old. It has been introduced by M.~Frank and D.~Larson in  the late 1990's and since then it serves as a useful tool and, at the same time, as a subject of research interest in its own. It turned out that frames in Hilbert $C^*$-modules share many properties with classical frames for Hilbert spaces. However, this is limited only to the class of, in the Frank-Larson terminology, standard frames. Those are the frames for which the corresponding analysis operator takes values in the standard Hilbert module $\ell^2(\textsf{A})$. Frames in some weaker sense were not studied by now. There are some generalizations such as modular g-frames (there is a number of recent articles on g-frames for Hilbert $C^*$-modules) and outer frames (see \cite{AB}). The later class is well suited for Hilbert $C^*$-modules over non-unital $C^*$-algebras, but in case the underlying $C^*$-algebra possesses a unit, the outer frames simply become the usual standard frames.

So, the question of an appropriate concept of a weak frame (in any sense weaker than with respect to the norm) is still open. On the other hand, researchers have encountered situations in which such kind of modular frames could play a role (see \cite{CasL}, p.~97 and \cite{CocoL}, pp.~11,12).

The reason why up to date no theory of weak modular frames has been developed, even for some special classes of $C^*$-algebras, lies in the fact that in order to study such frames one has to introduce a suitable target space (a Hilbert $C^*$-module) for the analysis operator. And since no such target module was available, the whole concept remained unfounded.

Here, in the first part of the paper, we describe the dual $\ell^2(\textsf{A})^{\prime}$ of the standard Hilbert $C^*$-module $\ell^2(\textsf{A})$ over an arbitrary $C^*$-algebra \textsf{A}. It is to some extent surprising that there is no complete description of this dual available from the existing literature, having in mind importance of the standard module $\ell^2(\textsf{A})$ (which goes back to G.G. Kasparov, \cite{K}). It turns out that this description is particularly nice when the underlying algebra is a von Neumann algebra. When this is the case, the dual $\ell^2(\textsf{A})^{\prime}$ is realized in a concrete way precisely as we know (from the theoretical viewpoint; see \cite{P}) it should
be: as a self-dual Hilbert $\textsf{A}$-module $\ell^2_{\text{strong}}(\textsf{A})$ that is isometrically isomorphic to $\ell^2(\textsf{A})^{\prime}$, which contains $\ell^2(\textsf{A})$, and whose $\textsf{A}$-valued inner product extends the original inner product on $\ell^2(\textsf{A})$); see Theorem \ref{dual} and Corollary \ref{Paschke vN} below.

Once the dual module was identified and described, it was natural to try to introduce a concept of a weak frame for Hilbert $C^*$-modules over von Neumann algebras; see Definition \ref{very first fd} below. In fact, the convergence of the series in the frame definition is dictated by the inner product in $\ell^2(\textsf{A})^{\prime}$; it was of substantial importance to ensure that the analysis operators for such weak frames take values in the dual module $\ell^2(\textsf{A})^{\prime}$. Another necessary property of the analysis (as well as the synthesis) operator is adjointability. This is the reason why we restricted ourselves in the first step to the class of self-dual modules over von Neumann algebras. In the second step we extend the theory to general Hilbert modules over von Neumann algebras by working in the dual. It turns out that weak Bessel sequences and weak modular frames have properties, with respect to some topology weaker than the norm topology, similar to those of standard Bessel sequences and frames.

In the last part of the paper this new concept is applied. It turns out that our weak frames are well suited for application in Gabor analysis. This was already indicated in \cite{CasL} and \cite{CocoL}. Weak frames (resp.~weak Bessel sequences) of translates in the Lebesgue-Bochner module $L_{\frac{1}{b}}^{\infty}(\ell^2)$ of the form $(T_{na}\frac{1}{\sqrt{b}}g)_{n\in \Bbb Z}$ are in a bijective correspondence with Gabor frames (resp.~Bessel sequences) $G(g,a,b)$ in $L^2(\Bbb R)$, where $a$ and $b$ are the translation and the modulation parameter. Thus, weak frames are non-standard modular frames for a Hilbert $L^{\infty}[0,\frac{1}{b}]$-module $L_{\frac{1}{b}}^{\infty}(\ell^2)$ in the sense that perfectly fits into Gabor analysis. Note that in \cite{CocoL}, pp.~11,12, M.~Coco and M.C.~Lammers have already observed the correspondence of Gabor frames and certain non-standard modular frames in $L_{\frac{1}{b}}^{\infty}(\ell^2)$, while noticing that a concept of non-standard modular frames has not been developed at the time, so the idea of using Hilbert $C^*$-modules in that context could not be exploited.
Here, after developing the theory of weak modular frames, we explore this connection by demonstrating in this new (modular) light simple proofs of some well known results concerning Gabor Bessel sequences and Gabor frames; see Theorem \ref{CC revisited}, Remark \ref{karakterizacija Gabor Parseval}, and Remark \ref{WR} below.  In addition, a modular form of the Walnut representation is  proved (Proposition \ref{bolji walnut}) and, as a consequence, a result of Walnut type for classical Gabor Bessel sequences (Theorem \ref{mozda novi walnut}) is obtained.

\vspace{.2in}

The paper is organized as follows. In the rest of this introduction we shall fix our terminology and notation regarding the modular part of the paper. All necessary preliminaries and notation concerned with Gabor systems will be introduced at the beginning of Section 5.

Section 2 is completely devoted to the identification of the dual $\ell^2(\textsf{A})^{\prime}$ of the standard Hilbert $C^*$-module $\ell^2(\textsf{A})$. We have included some of the arguments already known from the literature trying to make the exposition self-contained as much as possible. The main results are Theorem \ref{dual} and its Corollary \ref{Paschke vN}.

In Section 3 weak Bessel sequences and weak frames are introduced and their fundamental properties are derived. As mentioned above, it was necessary to restrict our discussion to the class of Hilbert $C^*$-modules over von Neumann algebras in order to ensure that the target space for the analysis operators is again a Hilbert $C^*$-module over the same algebra. In this section we have founded the theory of weak frames (resp.~weak Bessel sequences) by obtaining all fundamental results such as unconditional convergence of the series representing synthesis operator (with respect to certain weak topology), invertibility of the frame operator, the reconstruction formula, etc.

In Section 4 we discuss the correspondence of adjointable operators on the standard Hilbert $C^*$-module $\ell^2(\textsf{A})$ over a von Neumann algebra \textsf{A} with infinite matrices with coefficients from \textsf{A}. Among other results, a generalization of the Schur test is proved. It should be pointed out that the results in this section are obtained under the additional assumption that the underlying von Neumann algebra is commutative.

Section 5 introduces Hilbert $C^*$-modules and spaces relevant for Gabor analysis and discusses some of their properties. The central space in our study is
$L^{\infty}[0,\frac{1}{b}]$-module $L_{\frac{1}{b}}^{\infty}(\ell^2)$ which is as a Hilbert $C^*$-module unitarily equivalent to the dual $\ell^2(L^{\infty}[0,\frac{1}{b}])^{\prime}$ of $\ell^2(L^{\infty}[0,\frac{1}{b}])$. Here is $b$ an arbitrary positive number that will in fact play the role of the modulation parameter for a given Gabor system. The main result of the section is Theorem \ref{CasCocLem} in which we establish a bijective correspondence of Gabor frames/Bessel sequences and weak frames/weak Bessel sequences of translates in $L_{\frac{1}{b}}^{\infty}(\ell^2)$.

Finally, in Section 6 we discuss some consequences. Various results on Gabor systems are (re)obtained. In particular a discussion on Wallnut representation is included.

The readers who are primarily interested in the Gabor part of the paper may prefer to start reading beginning with Section 5 and to turn back to the theoretical background concerning weak modular frames (i.e.~Sections 3 and 4), when needed.

\vspace{.2in}

Throughout the paper we work with left Hilbert $C^*$-modules.
Recall that a Hilbert $C^*$-module over a $C^*$-algebra $\textsf{A}$ is a complex vector space $X$ that is also a left $\textsf{A}$-module equipped with an inner product $\langle \cdot,\cdot \rangle: X\times X \rightarrow \textsf{A}$ that is linear in the first, anti-linear in the second variable and satisfies $\langle x,x\rangle \geq 0$ for all $x\in X$,
$\langle x,x\rangle = 0$ if and only if $x=0$, $\langle ax,y\rangle=a\langle x,y\rangle $ and $\langle y,x\rangle= \langle x,y\rangle^*$ for all $x,y\in X$, $a\in \textsf{A}$, and such that $X$ is complete with respect to the norm $\|x\|=\|\langle x,x\rangle\|^{\frac{1}{2}}$.

We opted (between two equivalent choices) to work with left Hilbert $C^*$-modules despite some technical difficulties that arise in establishing correspondence of operators with infinite matrices. This was motivated by possible applications in Gabor analysis where passing back and forth between Gabor systems in $L^2(\Bbb R)$ and modular Bessel sequences is needed, and hence it is more convenient to work with inner products that are linear in the same (i.e.~the first) argument in both structures.

If $X$ and $Y$ are Hilbert $C^*$-modules over the same $C^*$-algebra we denote by $\Bbb B(X,Y)$ the Banach space of all adjointable operators from $X$ into $Y$.

The key role in our considerations in Section 2 will play the standard Hilbert $C^*$-module $\ell^2(\textsf{A})$ that is defined as
$$
\ell^2(\textsf{A})=\left\{(a_n)_n : a_n \in
\textsf{A},\,\,
\sum_{n=1}^{\infty}a_na_n^* \mbox{ converges in norm in } \textsf{A} \right\}$$
and equipped with the inner product
$$
\langle (a_n)_n, (b_n)_n \rangle = \sum_{n=1}^{\infty}a_nb_n^*.
$$
It turns out that this series converges in $\textsf{A}$. In fact, it converges unconditionally since by polarization it can be written as a linear combination of four series of the form $\sum_{n=1}^{\infty}c_nc_n^*$, with $(c_n)_n \in \ell^2(\textsf{A})$ and it is well known that when a series of positive elements converges in a $C^*$-algebra, it necessarily converges unconditionally.

This implies that there is no loss of generality in working with countable systems indexed by natural numbers (i.e.~sequences), although in the second part of the paper we will naturally use indexation over the integers when working with Gabor systems.

We shall often assume (in particular, in Section 2) that our $C^*$-algebra $\textsf{A}$ acts non-degenerately on a Hilbert space $H$. This is not a restriction since any $C^*$-algebra can be faithfully and non-degenerately represented on a Hilbert space. We denote by $LM(\textsf{A})$, $RM(\textsf{A})$ and $M(\textsf{A})$ the sets of left multipliers, right multipliers, and (two-sided) multipliers of $\textsf{A}$, respectively. Recall that $LM(\textsf{A})$ and $RM(\textsf{A})$ are Banach algebras, while $M(\textsf{A})$ is a $C^*$-algebra. Note also that all these algebras are contained in the closure $\overline{\textsf{A}}^s$ of $\textsf{A}$ in the strong operator topology which in this situation coincides with the bicommutant $\textsf{A}^{\prime \prime}$ of $\textsf{A}$.

We refer the reader to \cite{Ped} and \cite{WO} for general facts on $C^*$-algebras and Hilbert $C^*$-modules.

\vspace{.1in}

Further notations will we explained in the course of exposition.

\vspace{.3in}

\section{$\ell^2(\textsf{A})^{\prime}$}\label{druga sekcija}

\vspace{.1in}

Given a Hilbert $C^*$-module $X$ over a $C^*$-algebra $\textsf{A}$, we denote by $X^*$ its adjointable dual; i.e.~$X^*=\Bbb B(X,\textsf{A})$, and by $X^{\prime}$ its dual, that is, the Banach space of  all bounded module maps from $X$ into $\textsf{A}$. By a module map we understand a linear operator $L:X \rightarrow \textsf{A}$ which satisfies $L(ax)=aL(x)$ for all $x$ in $X$ and $a$ in $\textsf{A}$.

The dual $X^{\prime}$ becomes a right Banach module over $\textsf{A}$ with the action of $\textsf{A}$ on $X^{\prime}$ defined by $(L,b) \mapsto Lb$, $(Lb)(x)=L(x)b$, $x$ in $X$.

Clearly, $X^*$ is a closed subspace (in fact, a submodule) of $X^{\prime}$.  We also know that $X$ is embedded in $X^{*}\subseteq X^{\prime}$ in a standard way: for $x$ in $X$ we define the map $L_x : X \rightarrow \textsf{A}$ by $L_x(y)=\langle y,x \rangle$. It is evident that $L_x$ belongs to $X^*=\Bbb B(X,\textsf{A})$. Indeed, $L_x^*$ is given by $L_x^*(a)=ax$, $a \in \textsf{A}$. So, we have the map
\begin{equation}\label{e1}
\varphi : X \rightarrow X^* \subseteq X^{\prime},\,\,\,\varphi(x)=L_x,\,\,\,L_x(y)=\langle x,y\rangle.
\end{equation}
It is easy to show that $\varphi$ is an anti-linear isometry.  What is more, it turns out with the help of the Cohen-Hewitt factorization (each $x$ in $X$ is of the form $x=b^*v$ for some $v$ in $X$ and $b^*$ in $\textsf{A}$) that every $L_x$ is in fact a "rank-one" operator $\theta_{v,b}$, where $\theta_{v,b}(y)=\langle y,v\rangle b$, $y\in X$. It is also known that the image of $\varphi$ coincides with the space of all "compact" operators $\Bbb K(X,\textsf{A})$ (see Lemma 2.32 in \cite{RW}).

A Hilbert $C^*$-module $X$ is said to be self-dual if each bounded module map (i.e.~an element of $X^{\prime}$) arises by taking the inner product with some fixed element of $X$. Thus, $X$ is self-dual if for each $L\in X^{\prime}$ there exists $x\in X$ such that $L=L_x$.

It is always desirable to determine precisely both the dual and the adjointable dual of a Hilbert $C^*$-module under consideration. Our goal in this section is to identify $\ell^2(\textsf{A})^{*}$ and $\ell^2(\textsf{A})^{\prime}$ where $\ell^2(\textsf{A})$ is the standard Hilbert $C^*$-module over $\textsf{A}$ (sometimes called the Hilbert space over $\textsf{A}$). We do not impose any restrictions on the underlying $C^*$-algebra $\textsf{A}$; in particular, we do not assume that $\textsf{A}$ is unital.

\vspace{.1in}

Let us begin by recalling the definition of the multiplier module $M(\ell^2(\textsf{A}))$ of $\ell^2(\textsf{A})$ (cf.~\cite{BG}):
\begin{equation}\label{glavni modul 0}
M(\ell^2(\textsf{A}))=\left\{(c_n)_n : c_n \in
M(\textsf{A}),\,\,
\sum_{n=1}^{\infty}c_nc_n^* \mbox{ converges strictly}\right\}.
\end{equation}
It is known (\cite{BG}, Theorem 2.1) that $M(\ell^2(\textsf{A}))$ is a Hilbert $C^*$-module over $M(\textsf{A})$ with the $M(\textsf{A})$-valued inner product on $M(\ell^2(\textsf{A}))$ given by
$$\langle (c_n)_n,(d_n)_n\rangle=\mbox{(strict)}\,\sum_{n=1}^{\infty}c_nd_n^*.$$
Here the strict topology on $M(\textsf{A})$ is the locally convex topology generated by the seminorms $x\mapsto \|ax\|$ and $x\mapsto \|xa\|$, $x\in M(\textsf{A})$, for all $a \in \textsf{A}$. Since we assumed that \textsf{A} acts non-degenerately on $H$, each strictly convergent net converges also in the strong operator topology to the same limit.

Note that $\ell^2(\textsf{A})$ is contained in $M(\ell^2(\textsf{A}))$ and the norm on $\ell^2(\textsf{A})$ inherited from $M(\ell^2(\textsf{A}))$ coincides with the original Hilbert $C^*$-norm defined on $\ell^2(\textsf{A})$. Clearly, when \textsf{A} is unital, we have $M(\textsf{A})=\textsf{A}$ and consequently $M(\ell^2(\textsf{A}))=\ell^2(\textsf{A})$.

\vspace{.1in}

We will work in another, even larger Hilbert $C^*$-module over $\overline{\textsf{A}}^s=\textsf{A}^{\prime \prime}$. Let
\begin{equation}\label{glavni modul}
\ell^2_{\text{strong}}(\overline{\textsf{A}}^s)=
\left\{ (c_n)_n: c_n\in \overline{\textsf{A}}^s,\,\,\sup_{N}\left\|\sum_{n=1}^Nc_nc_n^*\right\|<\infty\right\}.
\end{equation}
Observe that the condition $\sup_{N}\left\|\sum_{n=1}^Nc_nc_n^*\right\|<\infty$ implies that the series $\sum_{n=1}^{\infty}c_nc_n^*$ converges strongly. On the other hand, applying the uniform boundedness principle we conclude that the converse is also true. Thus, \eqref{glavni modul} can be rewritten as
\begin{equation}\label{glavni modul 1}
\ell^2_{\text{strong}}(\overline{\textsf{A}}^s)=
\left\{ (c_n)_n: c_n\in \overline{\textsf{A}}^s,\,\,\sum_{n=1}^{\infty}c_nc_n^* \mbox{ converges strongly}\right\}.
\end{equation}
Since $M(\textsf{A})\subseteq \overline{\textsf{A}}^s$ and the strong operator topology is weaker than the strict topology on $M(\textsf{A})$, by comparing \eqref{glavni modul 0} and \eqref{glavni modul 1}, we see that $M(\ell^2(\textsf{A}))$ is contained in $\ell^2_{\text{strong}}(\overline{\textsf{A}}^s)$. There is another space of our interest in between:
\begin{equation}\label{glavni modul 2}
\ell^2_{\text{strong}}(LM(\textsf{A}))=\left\{ (c_n)_n: c_n\in LM(\textsf{A}),\,\,\sup_{N}\left\|\sum_{n=1}^Nc_nc_n^*\right\|<\infty\right\}.
\end{equation}
Thus, we have the following chain of inclusions:

\begin{equation}\label{inkluzije}
\ell^2(\textsf{A}) \subseteq M(\ell^2(\textsf{A})) \subseteq \ell^2_{\text{strong}}(LM(\textsf{A})) \subseteq \ell^2_{\text{strong}}(\overline{\textsf{A}}^s).
\end{equation}

\vspace{.2in}

\begin{prop}\label{Paschke}
If $\textsf{A}\subseteq \Bbb B(H)$ is a $C^*$-algebra that acts non-degenerately on $H$, $\ell^2_{\text{strong}}(\overline{\textsf{A}}^s)$ is a Hilbert $C^*$-module over $\overline{\textsf{A}}^s$ with the inner product defined by
$$\langle (x_n)_n,(y_n)_n\rangle=\text{(strong)}\,\sum_{n=1}^{\infty}x_ny_n^*,$$
where $\text{(strong)}$ refers to the strong operator topology on $B(H)$.
The norm on
$M(\ell^2(\textsf{A}))$ inherited from $\ell^2_{\text{strong}}(\overline{\textsf{A}}^s)$ coincides with the original Hilbert $C^*$-norm on $M(\ell^2(\textsf{A}))$.
In addition, $\ell^2_{\text{strong}}(\overline{\textsf{A}}^s)$ contains $\ell^2_{\text{strong}}(LM(\textsf{A}))$  and $M(\ell^2(\textsf{A}))$ as closed subspaces.
\end{prop}
\proof
Let $(x_n)_n$ and $(y_n)_n$ be any two sequences in $\ell^2_{\text{strong}}(\overline{\textsf{A}}^s)$.

We first claim that the series $\sum_{n=1}^{\infty}x_ny_n^*$ converges strongly.
Let $C=\sup\left\{\left\| \sum_{n=1}^Nx_nx_n^*\right\|: N\in \Bbb N \right\}$. Since the sequence $\left(\left\| \sum_{n=1}^Ny_ny_n^*\right\|\right)_N$ is bounded, the series
$\sum_{n=1}^{\infty}y_ny_n^*$ converges strongly. Take any $\xi$ in $H$ and fix $\epsilon>0$. Then there exists $N_0$ such that
$$
N_2>N_1\geq N_0 \Rightarrow \left\| \left(\sum_{n=N_1+1}^{N_2}y_ny_n^*\right)\xi\right\|<\epsilon.
$$
Recall the strong version of the Cauchy-Schwarz inequality that holds true in any Hilbert $C^*$- module $X$: $\langle x,y\rangle^*\langle x,y\rangle \leq \|\langle x,x\rangle\| \langle y,y\rangle$ (\cite{L}, Proposition 1.1). Applying this inequality in the Hilbert $C^*$-module $(\overline{\textsf{A}}^s)^{N_2-N_1}$ we get
$$
\left( \sum_{n=N_1+1}^{N_2}x_ny_n^*\right)^*\left( \sum_{n=N_1+1}^{N_2}x_ny_n^*\right)\leq \left\|\sum_{n=N_1+1}^{N_2}x_nx_n^* \right\|\left(\sum_{n=N_1+1}^{N_2}y_ny_n^*\right).
$$
We now apply the operators from both sides of this inequality to $\xi$ and take the inner product in $H$ by $\xi$. In this way we obtain
$$
\left\|\left(\sum_{n=N_1+1}^{N_2}x_ny_n^*\right)\xi \right\|^2\leq C \left\langle \left(\sum_{n=N_1+1}^{N_2}y_ny_n^*\right)\xi, \xi\right\rangle \leq C\epsilon \|\xi\|.
$$
This proves our claim: the series $\sum_{n=1}^{\infty}x_ny_n^*$  converges strongly. Observe that this implies two things. First, we conclude that the series $\sum_{n=1}^{\infty}(x_n+y_n)^*(x_n+y_n)$ converges strongly and hence by the uniform boundedness principle the sequence $\left(\sum_{n=1}^{N}(x_n+y_n)(x_n+y_n)^*\right)_N$ is bounded; thus,
$\ell^2_{\text{strong}}(\overline{\textsf{A}}^s)$ is closed under addition.
Secondly, we now have a well defined $\overline{\textsf{A}}^s$-valued inner product on
$\ell^2_{\text{strong}}(\overline{\textsf{A}}^s)$ given by
$$
\langle (x_n)_n, (y_n)_n\rangle = (\text{strong})\, \sum_{n=1}^{\infty}x_ny_n^*.
$$
Notice that the strong convergence of the series $\sum_{n=1}^{\infty}x_ny_n^*$
does not  imply a priori the strong convergence of the series $\sum_{n=1}^{\infty}y_nx_n^*$. However, by the preceding discussion both series do converge strongly and hence weakly, which implies $\langle (x_n)_n, (y_n)_n\rangle=\langle (y_n)_n, (x_n)_n\rangle^*$ for all sequences $(x_n)_n$, $(y_n)_n$ from $\ell^2_{\text{strong}}(\overline{\textsf{A}}^s)$.

So, $\ell^2_{\text{strong}}(\overline{\textsf{A}}^s)$ has the structure of an inner-product $\overline{\textsf{A}}^s$-module.

Let us now show that this module is complete. Take a Cauchy sequence $(c_n)_n$ in $\ell^2_{\text{strong}}(\overline{\textsf{A}}^s)$ and put $c_n=(x_k^n)_k$, $n\in \Bbb N$. Fix $\epsilon >0$. Then there exists $n_0$ with the property
\begin{equation}\label{e20}
m,n\geq n_0 \Rightarrow \|c_n-c_m\|<\epsilon
\end{equation}
which means
\begin{equation}\label{e21}
m,n\geq n_0 \Rightarrow \left\|(\text{strong})\,\sum_{k=1}^{\infty}(x_k^n-x_k^m) (x_k^n-x_k^m)^* \right\|<\epsilon^2.
\end{equation}
Since the sequence $(\sum_{k=1}^{N}(x_k^n-x_k^m) (x_k^n-x_k^m)^*)_N$ is an increasing strongly convergent sequence of positive operators, we conclude that
\begin{equation}\label{e22}
\left\|(x_k^n-x_k^m) (x_k^n-x_k^m)^* \right\|<\epsilon^2,\,\,\,\forall k \in \Bbb N
\end{equation}
and hence
$$
\left\| x_k^n-x_k^m\right\|<\epsilon,\,\,\,\forall k \in \Bbb N.
$$
This shows that $(x_k^n)_n$ is a \emph{norm}-Cauchy sequence in $\overline{\textsf{A}}^s$.
Put
\begin{equation}\label{sic}
x_k^0=\lim_{n\rightarrow \infty}x_k^n \in \overline{\textsf{A}}^s,\,\,\,k \in \Bbb N
\end{equation}
and
$$
c_0=(x_1^0,x_2^0,x_3^0,\ldots).
$$

Consider again $n_0$ for which we have \eqref{e20} and \eqref{e21} and fix $n\geq n_0$. Then we have for each $K$
\begin{eqnarray}\label{e31}
\left\| \sum_{k=1}^K(x_k^n-x_k^0)(x_k^n-x_k^0)^*\right\|&=&\lim_{m\rightarrow \infty}\left\| \sum_{k=1}^K(x_k^n-x_k^m)(x_k^n-x_k^m)^*\right\|\\\nonumber
 &\leq&\limsup_{m\rightarrow \infty}\left\|\text{(strong)}\, \sum_{k=1}^{\infty}(x_k^n-x_k^m)(x_k^n-x_k^m)^*\right\|\\\nonumber
 &=&\limsup_{m\rightarrow \infty}\|c_n-c_m\|^2\stackrel{\eqref{e20}}{\leq}\epsilon^2.
\end{eqnarray}

This tells us that the sequence $\left(  \sum_{k=1}^K(x_k^n-x_k^0)(x_k^n-x_k^0)^*\right)_K$ is bounded, and hence $c_n-c_0$ belongs to $\ell^2_{\text{strong}}(\overline{\textsf{A}}^s)$. In particular, $c_0$ is in $\ell^2_{\text{strong}}(\overline{\textsf{A}}^s)$.
Moreover,
$$
\|c_n-c_0\|^2=\left\| (\text{strong}) \sum_{k=1}^{\infty}(x_k^n-x_k^0)(x_k^n-x_k^0)^*\right\|
\leq\sup_K\left\| \sum_{k=1}^{\infty}(x_k^n-x_k^0)(x_k^n-x_k^0)^*\right\|\stackrel{\eqref{e31}}{\leq}\epsilon^2.
$$
Thus, $\ell^2_{\text{strong}}(\overline{\textsf{A}}^s)$ is complete, so it is a Hilbert $C^*$-module over $\overline{\textsf{A}}^s$.

\vspace{.1in}

To prove the second assertion, take any $(t_n)_n\in M(\ell^2(\textsf{A}))$. Its original norm arising from the Hilbert $M(\textsf{A})$-module structure on $M(\ell^2(\textsf{A}))$ is equal to $\left\|(\text{strict})\,\sum_{n=1}^{\infty}t_nt_n^* \right\|^{\frac{1}{2}}$. Since $\textsf{A}$ acts non-degenerately on $H$, strict convergence in $M(\ell^2(\textsf{A}))$ induced by $\textsf{A}$ implies convergence in the strong operator topology. So, $(\text{strict})\,\sum_{n=1}^{\infty}t_nt_n^*=(\text{strong})\,\sum_{n=1}^{\infty}t_nt_n^*$ and hence $\left\|(\text{strict})\,\sum_{n=1}^{\infty}t_nt_n^*\right\|^{\frac{1}{2}}=\left\|(\text{strong})\,\sum_{n=1}^{\infty}t_nt_n^*\right\|^{\frac{1}{2}}$.

Finally, recall our observation preceding \eqref{sic}. (See also \eqref{sic}.)
If we have a sequence $(c_n)_n$, $c_n=(x_k^{n})_k$ such that each $x_k^{n}$ belongs to some norm-closed subalgebra $\textsf{B}$ of $\overline{\textsf{A}}^s$ and if $(c_n)_{n}$ converges
in our Hilbert $C^*$-module $\ell^2_{\text{strong}}(\overline{\textsf{A}}^s)$ to $c_0=(c_k^0)_k$, then all $c_k^0$ (i.e.~the component-wise limits) also belong to $\textsf{B}$.
This proves that $\ell^2_{\text{strong}}(LM(\textsf{A}))$ is closed in
$\ell^2_{\text{strong}}(\overline{\textsf{A}}^s)$. As for $M(\ell^2(\textsf{A}))$, if $c_0=\lim_{n\rightarrow \infty}c_n$ in $\ell^2_{\text{strong}}(\overline{\textsf{A}}^s)$ with $c_n\in M(\ell^2(\textsf{A}))$ for every $n$, then $(c_n)_n$ is a Cauchy sequence in $M(\ell^2(\textsf{A}))$ and since it is a Hilbert $C^*$-module (hence, complete) and its original norm coincides with the norm inherited from $\ell^2_{\text{strong}}(\overline{\textsf{A}}^s)$, we must have $c_0\in M(\ell^2(\textsf{A}))$.
\qed

\vspace{.2in}

\begin{remark}\label{nije modul}
Observe that, although $\ell^2_{\text{strong}}(LM(\textsf{A}))$ has the structure of a Banach $\textsf{A}$-module, it is only a closed subspace, and not a submodule of $\ell^2_{\text{strong}}(\overline{\textsf{A}}^s)$.
\end{remark}

\vspace{.1in}

We can now state our main theorem in this section.

\vspace{.2in}

\begin{theorem}\label{dual}
Let $\textsf{A}\subseteq \Bbb B(H)$ be a $C^*$-algebra that acts non-degenerately on $H$. Then the map
$$
\varphi : \ell^2_{\text{strong}}(LM(\textsf{A})) \rightarrow \ell^2(\textsf{A})^{\prime},\,\,\,\varphi((t_n)_n)=L_{(t_n)_n},\,\,\,L_{(t_n)_n}((a_n)_n)=\sum_{n=1}^{\infty}a_nt_n^*,
$$
where this series converges in norm in $\textsf{A}$, is an anti-linear isometric isomorphism of Banach $\textsf{A}$-modules. Moreover, its restriction to $M(\ell^2(\textsf{A}))$
$$
\varphi |_{M(\ell^2(\textsf{A}))} : M(\ell^2(\textsf{A}))\rightarrow \ell^2(\textsf{A})^*
$$
is an anti-linear isometric isomorphism of Hilbert $M(\textsf{A})$-modules. In particular,  $\varphi |_{\ell^2(\textsf{A})}$ extends the embedding introduced in \eqref{e1}: $\varphi(x)=L_x$, $x\in \ell^2(\textsf{A})$.

\end{theorem}

\vspace{.2in}

Observe that the definition of $\varphi$ makes sense since $t\mapsto t^*$ is a bijection from $LM(\textsf{A})$ to $RM(\textsf{A})$ and hence $at^*\in \textsf{A}$ for each $a\in \textsf{A}$ and $t\in LM(\textsf{A})$.

\vspace{.1in}

Note that in the non-unital case the adjointable dual $\ell^2(\textsf{A})^*$ is much larger than $\ell^2(\textsf{A})$ - a remarkable, but somewhat surprising fact. We will explain later in Remark \ref{objasnjenje} the reason why the idea of passing from $\ell^2(\textsf{A})$ to $\ell^2(\tilde{\textsf{A}})$, where $\tilde{\textsf{A}}$ is the minimal unitization of \textsf{A}, and trying to describe the dual of $\ell^2(\textsf{A})$ in terms of  $\tilde{\textsf{A}}$ turns out to be rather na\" ive and of no help.

\vspace{.2in}

To prove this theorem we shall need a couple of auxiliary results. But first we proceed with some comments and consequences.

When $A$ is unital we have $LM(\textsf{A})=M(\textsf{A})=\textsf{A}$ and the sequence of inclusions \eqref{inkluzije} becomes
\begin{equation}\label{inkluzije e}
\ell^2(\textsf{A}) \subseteq \ell^2_{\text{strong}}(\textsf{A}) \subseteq \ell^2_{\text{strong}}(\overline{\textsf{A}}^s).
\end{equation}
Also note that the assumed non-degenerate action of \textsf{A} on $H$ implies that the unit in \textsf{A} is the identity operator on $H$. Thus, we have

\begin{cor}\label{dual e}
Let $\textsf{A}\subseteq \Bbb B(H)$ be a $C^*$-algebra that contains the identity operator on $H$. Then the map
$$
\varphi : \ell^2_{\text{strong}}(\textsf{A}) \rightarrow \ell^2(\textsf{A})^{\prime},\,\,\,\varphi((t_n)_n)=L_{(t_n)_n},\,\,\,L_{(t_n)_n}((a_n)_n)=\sum_{n=1}^{\infty}a_nt_n^*,
$$
where this series converges in norm in $\textsf{A}$, is an anti-linear isometric isomorphism of Banach $\textsf{A}$-modules. Moreover, its restriction to $\ell^2(\textsf{A})$
$$
\varphi |_{\ell^2(\textsf{A})} : \ell^2(\textsf{A})\rightarrow \ell^2(\textsf{A})^*
$$
is an anti-linear isometric isomorphism of Hilbert $\textsf{A}$-modules.
\end{cor}

\vspace{.2in}

We note that in the unital case the isometric $1-1$ correspondence of sequences in $\ell^2_{\text{strong}}(\textsf{A})$ and bounded module maps from $\ell^2(\textsf{A})^{\prime}$ is proved in Proposition 2.5.5 in \cite{MT}; see also Lemma 4.1 and Corollary 4.2 in \cite{F}. Observe also that one easily concludes from Theorem \ref{dual} that $\ell^2(\textsf{A})$ is self-dual if and only if
$\textsf{A}$ is finite-dimensional - a fact that is first proved in \cite{F}.
\vspace{.2in}

Finally, we can make further specialization by assuming that $\textsf{A}$ is a von Neumann algebra. When this is the case, we have $\overline{\textsf{A}}^s=\textsf{A}$ and the above sequence of inclusions \eqref{inkluzije e} reduces to just two Hilbert $\textsf{A}$-modules:
\begin{equation}\label{inkluzije von}
\ell^2(\textsf{A}) \subseteq \ell^2_{\text{strong}}(\textsf{A}).
\end{equation}

The preceding corollary applies. Here we see that $\ell^2(\textsf{A})^{\prime}$ is isometrically isomorphic to the Hilbert $\textsf{A}$-module $\ell^2_{\text{strong}}(\textsf{A})$. (Note that without the assumption that $\textsf{A}$ is strongly closed $\ell^2_{\text{strong}}(\textsf{A})$ is not a Hilbert $\textsf{A}$-module.)
In fact, even more is true.

\begin{cor}\label{Paschke vN}
Let $A$ be a von Neumann algebra. Then $\ell^2(\textsf{A})^{\prime}$ is anti-linearly isometrically isomorphic to the Hilbert $\textsf{A}$-module $\ell^2_{\text{strong}}(\textsf{A})$. Moreover, $\ell^2_{\text{strong}}(\textsf{A})$ is self-dual.
\end{cor}
\proof
We only need to prove self-duality.
Let $L\in \ell^2_{\text{strong}}(\textsf{A})^{\prime}$. Then by Corollary \ref{dual e} the restriction of $L$ to $\ell^2(\textsf{A})$ is of the form
$(a_n)_n \mapsto \sum_{n=1}^{\infty}a_nt_n^*$ for some sequence $(t_n)_n\in \ell^2_{\text{strong}}(\textsf{A})$. Consider $L_{(t_n)_n}:\ell^2_{\text{strong}}(\textsf{A}) \rightarrow \textsf{A}$, $L_{(t_n)_n}((x_n)_n)=\langle  (x_n)_n,(t_n)_n\rangle = \mbox{(strong)}\,\sum_{n=1}^{\infty}x_nt_n^*$. Clearly, $L$ and $L_{(t_n)_n}$ coincide on $\ell^2(\textsf{A})$. In other words, $L-L_{(t_n)_n}\in \ell^2_{\text{strong}}(\textsf{A})^{\prime}$ satisfies $(L-L_{(t_n)_n})|_{\ell^2(\textsf{A})}=0$. We claim that $L-L_{(t_n)_n}=0$. This can be seen in the following way.

By Theorem 3.2 in \cite{P} there is an $\textsf{A}$-valued inner product $[\cdot,\cdot]$ on  $\ell^2_{\text{strong}}(\textsf{A})$ which extends the original inner product on $\ell^2(\textsf{A})$ and such that $\ell^2_{\text{strong}}(\textsf{A})$ is a self-dual Hilbert $C^*$-module. Moreover, the norm on $\ell^2_{\text{strong}}(\textsf{A})$ arising from $[\cdot,\cdot]$ coincides with the operator norm $\ell^2(\textsf{A})^{\prime}$ and hence, by our Theorem \ref{dual}, with the norm arising from our inner product introduced in Proposition \ref{Paschke}. This is enough to conclude (see the last part of the proof of Theorem 3.2 in \cite{P}) that $L-L_{(t_n)_n}=0$.
\qed

%


\vspace{.1in}
Notice that we now have a concrete realization of Paschke's construction of a self-dual Hilbert $C^*$-module structure on $\ell^2(\textsf{A})^{\prime}$ extending that which is defined on $\ell^2(\textsf{A})$.

\vspace{.1in}

Let us now turn to the proof of  Theorem \ref{dual}.
We first need a lemma that is known (for example, see Lemma 1.4 in \cite{AB}). The proof is included here for completeness.

\vspace{.1in}

\begin{lemma}\label{Heuser}
Let $\textsf{A}\subseteq \Bbb B(H)$ be a $C^*$-algebra that acts non-degenerately on $H$. Let $(t_n)_n$ be a sequence of operators in $LM(\textsf{A})$. Then the following two conditions are equivalent:
\begin{itemize}
\item[(a)] The series $\sum_{n=1}^{\infty}a_nt_n^*$ is norm-convergent for all $(a_n)_n$ from $\ell^2(\textsf{A})$.
\item[(b)] The sequence $\left(\sum_{n=1}^Nt_nt_n^*\right)_N$ is bounded.
\end{itemize}
If (a) and (b) are satisfied then $L:\ell^2(\textsf{A}) \rightarrow \textsf{A}$, $L((a_n)_n)=\sum_{n=1}^{\infty}a_nt_n^*$ is a bounded module map for which we have $\|L\|=\lim_{N\rightarrow \infty}\|L^N\|=\lim_{N\rightarrow \infty}\left\|\sum_{n=1}^Nt_nt_n^*\right\|^{\frac{1}{2}}$ where, for each natural number $N$, the operator $L^N:\ell^2(\textsf{A}) \rightarrow \textsf{A}$ is defined as the $N$-th partial sum, $L^N((a_n)_n)=\sum_{n=1}^{N}a_nt_n^*$.
\end{lemma}
\proof
Assume (a). First note that each $L^N$ is bounded and $\|L^N\|\leq \left\| \sum_{n=1}^Nt_nt_n^* \right\|^{\frac{1}{2}}$. This is a direct consequence of the Cauchy-Schwarz inequality in the Hilbert $\Bbb B(H)$-module $\Bbb B(H)^N$ (the direct sum of $N$ copies of $\Bbb B(H)$). Since by the assumption the sequence $(L^N)_N$ strongly converges to $L$, the uniform boundedness principle tells us that $L$ is bounded.
It should also be observed that $L$ and all $L^N$ do take values in $\textsf{A}$ since each $t_n$ is a left centralizer of $\textsf{A}$ and hence $t_n^*$ is a right centralizer, so $at_n^*\in \textsf{A}$ for all $a$ in $\textsf{A}$.

Since each $L^N$ is a restriction of $L^{N+1}$ and, at the same time a restriction of $L$, we have
\begin{equation}\label{e2}
\|L^1\|\leq\|L^2\|\leq\|L^3\|\leq \ldots \leq \|L\|.
\end{equation}
On the other hand, we have for each $x\in \ell^2(\textsf{A})$
$$
\|Lx\|=\lim_{N\rightarrow \infty}\|L^Nx\|\leq\liminf_N\|L^N\|\,\|x\|
$$
which gives us
\begin{equation}\label{e3}
\|L\|\leq \liminf_N\|L^N\|.
\end{equation}
From \eqref{e2} and \eqref{e3} we conclude that $\|L\|=\lim_{N\rightarrow \infty}\|L^N\|$.

Let us now prove that $\|L^N\|= \left\| \sum_{n=1}^Nt_nt_n^* \right\|^{\frac{1}{2}}$. Denote for simplicity $\left\| \sum_{n=1}^Nt_nt_n^* \right\|$ by $c^N$. We already know from the beginning of the proof that $\|L^N\| \leq \sqrt{c^N}$ so we only need to show that $\|L^Nx\|$ can be made arbitrarily close to $\sqrt{c^N}$ by choosing suitable $x$ from the closed unit ball in $\ell^2(\textsf{A})$.

Let $(e_{\lambda})_{\lambda}$ be an approximate unit for $\textsf{A}$. Observe that $e_{\lambda} \xi \rightarrow \xi$ for each $\xi$ in $H$ since $\textsf{A}$ acts non-degenerately on $H$.

Let us now take, for each $\lambda$, the sequence $x_\lambda=\frac{1}{\sqrt{c^N}}(t_1e_\lambda,\ldots,t_Ne_\lambda,0,0,\ldots) \in \ell^2(\textsf{A})$ (recall that $t_n$'s are left centralizers which ensures that each $t_ne_\lambda$ belongs to $\textsf{A}$).

We now recall a well known inequality that holds true in every $C^*$-algebra: $baa^*b^*\leq \|a\|^2bb^*$. Since $\|e_{\lambda}\|\leq 1$, this gives us
$$
\sum_{n=1}^Nt_ne_\lambda e_\lambda t_n^* \leq \sum_{n=1}^N t_nt_n^*
$$
and hence
$$
\|x_{\lambda}\|=\frac{1}{\sqrt{c^N}}\left\| \sum_{n=1}^Nt_ne_\lambda e_\lambda t_n^*\right\|^{\frac{1}{2}}\leq \frac{1}{\sqrt{c^N}}\left\| \sum_{n=1}^N t_nt_n^*\right\|^{\frac{1}{2}}=1.
$$
Thus, all $x_{\lambda}$'s are in the closed unit ball. Now observe that
$$
\|L^N(x_{\lambda})\|=\frac{1}{\sqrt{c^N}}\left\| \sum_{n=1}^Nt_ne_\lambda t_n^*\right\|=\frac{1}{\sqrt{c^N}}\sup\left\{ \left\| \sum_{n=1}^Nt_ne_\lambda t_n^*\xi\right\|: \xi \in H,\,\,\|\xi\|\leq 1\right\}.
$$
In fact, this is enough to conclude the desired equality $\|L^N\| = \sqrt{c^N}$. One can argue as follows. Fix $\epsilon>0$. There exists $\xi_0 \in H$ such that $\|\xi_0\|\leq 1$ and
$$\left|c^N-\left\| \sum_{n=1}^Nt_n t_n^*\xi_0\right\|\right|<\epsilon.$$
For this $\xi_0$ and $t_1,\ldots,t_n$ there exists $\lambda_0$ such that
$$
\left|  \left\| \sum_{n=1}^Nt_ne_{\lambda_0} t_n^*\xi_0\right\|- \left\| \sum_{n=1}^Nt_n t_n^*\xi_0\right\|\right|<\epsilon.$$
This last approximation was possible since $e_{\lambda} t_n^*\xi_0\rightarrow t_n^*\xi_0$ for all $n=1,2,\ldots,N$.

After all, we conclude that $\|L^N(x_{\lambda_0})\|$ is sufficiently close to $\sqrt{c^N}$. In this way we have proved the implication (a) $\Rightarrow$ (b) and the second assertion of the lemma. It remains to prove that (b) implies (a).

Assume (b) and choose a positive number $c$ such that $\left\|\sum_{n=1}^Nt_nt_n^*\right\|\leq c$ for every $N$ in $\Bbb N$. Take any sequence $(a_n)_n$ from $\ell^2(\textsf{A})$.
If $\epsilon >0$ is given, we can find $N_0$ with the property
$$
N_2>N_1\geq N_0 \Rightarrow \left\| \sum_{n=N_1+1}^{N_2}a_na_n^* \right\| < \epsilon^2.
$$
From this we conclude
$$
\left\| \sum_{n=1}^{N_2}a_nt_n^*-\sum_{n=1}^{N_1}a_nt_n^*\right \|=\left\| \sum_{n=N_1+1}^{N_2}a_nt_n^*\right \|\leq \left\| \sum_{n=N_1+1}^{N_2}a_na_n^*\right \|^{\frac{1}{2}} \, \left\| \sum_{n=N_1+1}^{N_2}t_nt_n^*\right \|^{\frac{1}{2}}< \sqrt{c} \,\epsilon,
$$
where the first inequality is obtained using the Cauchy-Schwarz inequality in the $\Bbb B(H)$-Hilbert module $\Bbb B(H)^{N_2-N_1}$.

Thus, $\left( \sum_{n=1}^{N}a_nt_n^*\right)_N$ is a Cauchy sequence.
\qed


\vspace{.2in}

\begin{remark}
It is much easier to prove the inequality $\|L^N\|\geq \left\| \sum_{n=1}^Nt_n^*t_n \right\|^{\frac{1}{2}}={\sqrt{c^N}}$ when all $t_n$ belong not merely to $LM(\textsf{A})$ but to $\textsf{A}$. If so, one just observes that
$$
\left\| L^N(\frac{1}{\sqrt{c^N}}(t_1,t_2,\ldots,t_N,0,0,\ldots)) \right\|=\frac{1}{\sqrt{c^N}}\left\| \sum_{n=1}^Nt_nt_n^* \right\|=\left\| \sum_{n=1}^Nt_nt_n^* \right\|^{\frac{1}{2}}.
$$
Even if we take $t_n$ from $M(\textsf{A})$ the proof is easy. Namely, in that case we observe that all $L^N$ are adjointable operators; one easily checks that $(L^N)^*:\textsf{A}\rightarrow \ell^2(\textsf{A})$ is then given by $(L^N)^*a=(at_1,at_2,\ldots,at_N,0,0,\ldots)$. Thus, we have $L^N(L^N)^*a=\left( \sum_{n=1}^Nt_nt_n^*\right)a^*$ for all $a$ in $\textsf{A}$ and, since $\textsf{A}$ is an essential ideal in $M(\textsf{A})$, this implies $\|L^N(L^N)^*\|=\left\| \sum_{n=1}^Nt_nt_n^* \right\|$.

However, in the sequel we shall need the full force of the preceding lemma with the sequence $(t_n)_n$ of elements of $LM(\textsf{A})$.
\end{remark}

\vspace{.2in}

Suppose now that either of two equivalent conditions from Lemma \ref{Heuser} is satisfied. Clearly, the operator $L$ is a module map, so we have $L\in \ell^2(\textsf{A})^{\prime}$.
It is now natural to ask whether $L$ is adjointable. We provide the answer (or rather a reformulation of this question) in our next lemma.

Before stating the lemma it is convenient to recall a few facts concerning the multiplier module  $M(\ell^2(\textsf{A}))$ of $\ell^2(\textsf{A})$.
 It is known (see \cite{BG}) that $M(\ell^2(\textsf{A}))$ is the completion of $\ell^2(\textsf{A})$ with respect to the strict topology induced by $\textsf{A}$. This is the topology on $M(\ell^2(\textsf{A}))$ induced by the family of seminorms $x\mapsto \|ax\|$, $a\in \textsf{A}$, and $x\mapsto \|\langle x,y\rangle\|$, $y\in \ell^2(\textsf{A})$.

In general, each Hilbert $\textsf{A}$-module $X$ possesses the strict completion $M(X)$ which is a Hilbert $C^*$-module over $M(\textsf{A})$. As one might expect, when a $C^*$-algebra $\textsf{A}$ is regarded as a Hilbert $C^*$-module over itself, its strict completion coincides with $M(\textsf{A})$. Conveniently enough, each operator $T\in \Bbb B(X,Y)$ extends by strict continuity to a unique operator $T_M \in \Bbb B(M(X),M(Y))$. Moreover, we know that $(T_M)^*=(T^*)_M$.

\vspace{.1in}

Before we state our next result we need to establish one more notational convention. Given an element $a$ from a $C^*$-algebra $\textsf{A}$ and a natural number $n$, we  denote by $a^{(n)}$ the sequence $(0,\ldots,0,a,0,0,\ldots )$ whose only possibly non-trivial entry $a$ is on the $n$-th place. Clearly, $a\mapsto a^{(n)}$ is an embedding of $\textsf{A}$ into $\ell^2(\textsf{A})$.

\vspace{.1in}

\begin{lemma}\label{adjointable dual}
Let $\textsf{A}\subseteq \Bbb B(H)$ be a $C^*$-algebra that acts non-degenerately on $H$. Let $(t_n)_n$ be a sequence of operators in $LM(\textsf{A})$. Then the following two conditions are equivalent:
\begin{itemize}
\item[(a)] The series $\sum_{n=1}^{\infty}a_nt_n^*$ is norm-convergent for all $(a_n)_n$ from $\ell^2(\textsf{A})$ and the operator $L:\ell^2(\textsf{A})\rightarrow \textsf{A}$ defined by $L((a_n)_n)=\sum_{n=1}^{\infty}a_nt_n^*$ is adjointable.
\item[(b)] $(t_n)_n\in M(\ell^2(\textsf{A}))$.
\end{itemize}
\end{lemma}
\proof

Assume (a). So, there exists $L^*\in \Bbb B(\textsf{A},\ell^2(\textsf{A}))$. Consider the extended operators $L_M\in \Bbb B(M(\ell^2(\textsf{A})), M(\textsf{A}))$ and
$(L^*)_M\in \Bbb B(M(\textsf{A}), M(\ell^2(\textsf{A})))$.

Denote by $e$ the unit element in $M(\textsf{A})$ (which is in fact the identity operator on $H$). We first claim that
\begin{equation}\label{e11a}
L_M(e^{(n)})=t_n^*,\,\,\,\forall n \in \Bbb N.
\end{equation}

Indeed, if $(e_{\lambda})_{\lambda}$ is an approximate unit for $\textsf{A}$, we have
$$
L_M(e^{(n)})=\mbox{(strict)}\lim_{\lambda}L_M(e_{\lambda}^{(n)})=\mbox{(strict)}\lim_{\lambda}L(e_{\lambda}^{(n)})=\mbox{(strict)}\lim_{\lambda}e_{\lambda}t_n^*=\mbox{(left strict)}\lim_{\lambda}e_{\lambda}t_n^*=t_n^*.
$$
Note that the first equality above comes from the strict continuity of $L_M$ since $e_{\lambda}^{(n)}$ converges strictly in $M(\ell^2(\textsf{A}))$ to $e^{(n)}$.



\vspace{.1in}

Since $L_M$ takes values in $M(\textsf{A})$, \eqref{e11a} in particular shows us that $t_n^*\in M(\textsf{A})$, so our assumption (a) forces that the original sequence $(t_n)_n$ of left multipliers consists in fact of two-sided multipliers. What is more, we can now obtain a simple formula for $L^*a$. Namely, we know from the discussion preceding this lemma that $L^*a= (L^*)_Ma=(L_M)^*a$ and now we see that the $n$-th component of $L^*a$ is
$$
\langle  L^*a,e^{(n)}\rangle=\langle  (L_M)^*a,e^{(n)}\rangle =\langle a,L_Me^{(n)}\rangle \stackrel{\eqref{e11a}}{=}\langle a,t_n^*\rangle =at_n.
$$

Hence, the adjoint operator $L^*\in \Bbb B(\textsf{A}, \ell^2(\textsf{A}))$ is given by
\begin{equation}\label{e4}
L^*a=(at_1,at_2,at_3,\ldots ), \,\,\,\,a\in \textsf{A}.
\end{equation}
In particular, since $L^*a\in \ell^2(\textsf{A})$, \eqref{e4} implies that
\begin{equation}\label{e5}
\sum_{n=1}^{\infty}at_nt_n^*a^*\,\mbox{ converges in norm in } \,\textsf{A},\,\,\forall a \in \textsf{A}.
\end{equation}

Let $x_N=(t_1,t_2,\ldots,t_N,0,0,\ldots) \in M(\ell^2(\textsf{A}))$, $N \in \Bbb N$. We now claim that the sequence $(x_N)_N$ is strictly Cauchy in $M(\ell^2(\textsf{A}))$.

This is seen as follows. First, for each
$(a_n)_n$ from $\ell^2(\textsf{A})$ we have $\langle  (a_n)_n,x_N\rangle = \sum_{n=1}^Na_nt_n^*$ which is by our assumption (a) norm-convergent and hence Cauchy. Secondly, for each $a$ in $\textsf{A}$ the sequence $(ax_N)_N$ is also Cauchy since (assuming $N_2>N_1$) $\|ax_{N_2}-ax_{N_1}\|^2=\left\|\sum_{n=N_1+1}^{N_2}at_nt_n^*a^* \right\|$ which is small enough by \eqref{e5}.

Since $M(\ell^2(\textsf{A}))$ is strictly complete, there exists $s=(s_1,s_2,s_3,\ldots) \in M(\ell^2(\textsf{A}))$ for which we have $s=\text{strict}\,\lim_{N\rightarrow \infty}x_N$. In particular, this implies $\lim_{N\rightarrow \infty}\langle a^{(n)},x_N\rangle=\langle  a^{(n)},s\rangle$ for all $a \in \textsf{A}$ and $n$ in $\Bbb N$.
In other words, we have $at_n^*=as_n^*$ for all $a$ and $n$. This is enough to conclude $t_n^*=s_n^*$ for all $n$ and this gives us $s=(t_1,t_2,t_3,\ldots )\in M(\ell^2(\textsf{A}))$. After all, we have obtained the implication (a) $\Rightarrow$ (b).

Let us now assume (b). By the definition of $M(\ell^2(\textsf{A}))$, this means that the series $\sum_{n=1}^{\infty}t_nt_n^*$ $\textsf{A}$-strictly converges. In particular, this implies that $\sum_{n=1}^{\infty}t_nt_n^*a$ converges in norm in $\Bbb B(H)$ for all $a$ in $\textsf{A}$. As $A$ acts non-degenerately on $H$, the Cohen-Hewitt factorization theorem tells us that each $\xi \in H$ is of the form $\xi=a\eta$ for some $a\in \textsf{A}$ and $\eta \in H$. Hence $\sum_{n=1}^{\infty}t_nt_n^*$ converges strongly in $\Bbb B(H)$. This in turn implies, via the uniform boundedness principle, that the sequence $\left( \sum_{n=1}^Nt_nt_n^*\right)_N$ is bounded. By Lemma \ref{Heuser} the series $\sum_{n=1}^{\infty}a_nt_n^*$ converges for every $(a_n)_n$ in $\ell^2(\textsf{A})$. Hence, we have the operator $L:\ell^2(\textsf{A}) \rightarrow \textsf{A}$ defined by $L((a_n)_n)=\sum_{n=1}^{\infty}a_nt_n^*$.

Let us now define $R:\textsf{A} \rightarrow \ell^2(\textsf{A})$ by $Ra=(at_1,at_2,at_3,\ldots)$. Observe that $R$ is well defined because by our assumption (b)  $(t_n)_n\in M(\ell^2(\textsf{A}))$ and $\ell^2(\textsf{A})$ is the $\textsf{A}$-ideal submodule of $M(\ell^2(\textsf{A}))$.

To end the proof it remains to show that $L^*=R$. But this is evidently true since
$$
\langle L((a_n)_n), b\rangle=\sum_{n=1}^{\infty}a_nt_n^*b^*=\langle (a_n)_n,Rb\rangle
$$
for all $(a_n)_n$ from $\ell^2(\textsf{A})$ and $b$ in $\textsf{A}$.
\qed

\vspace{.2in}

We are now in position to prove Theorem \ref{dual}.

\vspace{.1in}

{\it Proof of Theorem \ref{dual}.}
We already know from Lemma \ref{Heuser} that
$$
\varphi : \ell^2_{\text{strong}}(LM(\textsf{A})) \rightarrow \ell^2(\textsf{A})^{\prime},\,\,\,\varphi((t_n)_n)=L_{(t_n)_n},\,\,\,L_{(t_n)_n}((a_n)_n)=\sum_{n=1}^{\infty}a_nt_n^*,
$$
is a well defined isometric anti-linear module map. Let us show that $\varphi$ is a surjection.

Take any $L\in \ell^2(\textsf{A})^{\prime}$. For a fixed $n\in \Bbb N$ we have a bounded module map $\textsf{A}\rightarrow \textsf{A}$ defined by $a\mapsto L(a^{(n)})$. By Theorem 1.5 from \cite{Lin} there exists a right multiplier $s_n\in RM(\textsf{A})$ such that $L(a^{(n)})=as_n$. (In fact, Theorem 1.5 in \cite{Lin} establishes a bijective correspondence between $X^{\prime}$ and $LM(\Bbb K(X,\textsf{A}))$. Here we have $X=\textsf{A}$ and $\Bbb K(\textsf{A}) \simeq \textsf{A}$ where "compact" operators on \textsf{A} are described as maps $\textsf{A} \rightarrow \textsf{A}$ of the form $a \mapsto at$ for some fixed $t \in \textsf{A}$. In addition, one should take into account that $RM(\textsf{A})=LM(\textsf{A})^*$.)

Put $s_n=t_n^*$ where $t_n$ is the corresponding left multiplier. In this way we have obtained a sequence $(t_n)_n$ in $LM(\textsf{A})$ such that $L(a^{(n)})=at_n^*$ for each $n$ in $\Bbb N$. In particular, we have for every $N$ in $\Bbb N$
\begin{equation} \label{zadnja}
L((a_1,a_2,\ldots,a_N,0,0,\ldots))=\sum_{n=1}^Na_nt_n^*, \,\,\,\forall a_1,a_2,\ldots,a_N\in\textsf{A}.
\end{equation}

The proof will be finished when we show that $(t_n)_n\in \ell^2_{\text{strong}}(LM(\textsf{A}))$. But this is easy. We know from \eqref{zadnja} by applying Lemma \ref{Heuser} that
$\|\sum_{n=1}^Nt_nt_n^*\|\leq L$, for each $N$ in $\Bbb N$. This shows that the sequence $(t_n)_n$ does belong to $\ell^2_{\text{strong}}(LM(\textsf{A}))$ and now \eqref{zadnja} implies that $L=L_{(t_n)_n}=\varphi((t_n)_n)$.

It remains to show that the restriction
$$
\varphi|_{M(\ell^2(\textsf{A}))} : M(\ell^2(\textsf{A}))\rightarrow \ell^2(\textsf{A})^*
,\,\,\,\varphi((t_n)_n)=L_{(t_n)_n},\,\,\,L_{(t_n)_n}((a_n)_n)=\sum_{n=1}^{\infty}a_nt_n^*,
$$
is an anti-linear isometric isomorphism of Hilbert $M(\textsf{A})$-modules.
By Lemma \ref{adjointable dual}, $\varphi$ maps $M(\ell^2(\textsf{A}))$ into $\ell^2(\textsf{A})^*$.
By the first part of the proof and by the last assertion of Proposition \ref{Paschke} $\varphi|_{M(\ell^2(\textsf{A}))} : M(\ell^2(\textsf{A}))\rightarrow \ell^2(\textsf{A})^*$ is an isometry.
Basically, this is what we had to show because we know from \cite{BG} that $M(\ell^2(\textsf{A}))=\Bbb B(\textsf{A}, \ell^2(\textsf{A}))$ and $\ell^2(\textsf{A})^*=\Bbb B(\ell^2(\textsf{A}),\textsf{A})$.
For reader's convenience we include the argument which shows that $\varphi|_{M(\ell^2(\textsf{A}))}$  is surjective onto $\ell^2(\textsf{A})^*$.

Take arbitrary adjointable $L:\ell^2(\textsf{A}) \rightarrow \textsf{A}$. Since $L^*\in \Bbb B(\textsf{A}, \ell^2(\textsf{A}))$ exists, we can extend it to
$(L^*)_M\in \Bbb B(M(\textsf{A}), M(\ell^2(\textsf{A})))$. Observe that $M(\textsf{A})$ contains the identity operator on $H$ which we again denote by $e$. Since $(L^*)_M$ is a module map, we have
$$
(L^*)_M(m)=(L^*)_M(me)=m\cdot (L^*)_M(e), \,\,\,\forall m\in M(\textsf{A}).
$$
Denote $(L^*)_M(e)=(t_n)_n\in M(\ell^2(\textsf{A}))$.
We now compute:
$$
L_M((y_n)_n)=L_M((y_n)_n)e=\langle  L_M((y_n)_n),e\rangle=\langle  (y_n)_n,(L^*)_M(e)\rangle=\langle  (y_n)_n,(t_n)_n,\rangle,\,\,\,\forall (y_n)_n \in M(\ell^2(\textsf{A})).
$$
In particular, this gives us
$$
L((a_n)_n)=L_M((a_n)_n)=\langle  (a_n)_n, (t_n)_n\rangle,\,\,\,\forall (a_n)_n \in\ell^2(\textsf{A}).
$$
\qed

\vspace{.2in}

We end this section with some examples and comments.

\vspace{.1in}

\begin{example}\label{prvi}
Let $H$ be a separable Hilbert space with an orthonormal basis $(\varepsilon_n)_n$. For each $n\in \Bbb N$ denote by $e_n$ the orthogonal projection to $\text{span}\,\{\varepsilon_n\}$.
Let $\textsf{A}=\Bbb K(H)$ - the $C^*$-algebra of all compact operators on $H$.

Clearly, $e_n\in \Bbb K(H) \subset M(\Bbb K(H))=\Bbb B(H)$. As the sequence $(\sum_{n=1}^Ne_ne_n^*)_N=(\sum_{n=1}^Ne_n)_N$ is bounded, Lemma \ref{Heuser} implies that the operator $L:\ell^2(\Bbb K(H))\rightarrow \Bbb K(H)$ given by $L((a_n)_n)=\sum_{n=1}^{\infty}a_ne_n$ is a well defined bounded module map. Since all bounded module maps of Hilbert $C^*$-modules over $\Bbb K(H)$ are adjointable (\cite{Mag}, Theorem 1), implication (b) $\Rightarrow$ (a) from Lemma \ref{adjointable dual} shows that $(e_n)_n \in M(\ell^2(\Bbb K(H)))$.

Indeed, this last conclusion can also be obtained directly. Anyhow, $L$ belongs to the adjointable dual $\ell^2(\Bbb K(H))^*$. It is useful to observe that $(e_n)_n$ does not belong to $\ell^2(\Bbb K(H))$.
\end{example}

\vspace{.1in}

\begin{remark}\label{objasnjenje}
Corollary \ref{dual e} tells us that, when $\textsf{A}$ contains the unit, $\ell^2(\textsf{A})^*$ is in a bijective correspondence with
$\ell^2(\textsf{A})$. When $\textsf{A}$ is non-unital, this is not the case since then $M(\ell^2(\textsf{A}))$ contains $\ell^2(\textsf{A})$ as a proper subset. (Usually is $M(\ell^2(\textsf{A}))$ way bigger than $\ell^2(\textsf{A}))$. This is demonstrated by the preceding example where we have concluded that $(e_n)_n \in M(\ell^2(\Bbb K(H)))$, but, clearly, $(e_n)_n \not \in \ell^2(M(\Bbb K(H)))=\ell^2(\Bbb B(H))\supseteq \ell^2(\Bbb K(H))$.

It is tempting to try to describe $\ell^2(\textsf{A})^*$ in the non-unital case by passing to $\ell^2(\tilde{\textsf{A}})$, where $\tilde{\textsf{A}}$ is the minimal unitization of $\textsf{A}$. Having in mind the unital case, one could try to establish a correspondence of $\ell^2(\textsf{A})^*$ with $\ell^2(\tilde{\textsf{A}})$. Let us take a closer look to that idea.

Take any
$L\in \ell^2(\textsf{A})^*$. Recall that this means that $L$ is an adjointable operator from $\ell^2(\textsf{A})$ to $\textsf{A}$. But now $\ell^2(\textsf{A})$ can be regarded as a Hilbert $C^*$-module over the unital $C^*$-algebra $\tilde{\textsf{A}}$, so we can regard $L$ as a map, call it $\tilde{L}$, from $\ell^2(\textsf{A})$ into $\tilde{\textsf{A}}$. Is $\tilde{L}$ adjointable? If we assume that the answer is  yes, it would follow that $Lx=\tilde{L}x=\tilde{L}x\cdot e=\langle \tilde{L}x,e\rangle=\langle x,(\tilde{L})^*e\rangle$ for all $x$ in $\ell^2(\textsf{A})$ and this leads to the conclusion that $L$ is represented by the sequence $(\tilde{L})^*e \in \ell^2(\textsf{A})$.

However, this cannot be true in general. Example \ref{prvi} provides $L$ that is represented by the sequence $(e_n)_n$ which does not belong to $\ell^2(\textsf{A})$.
So our innocently looking assumption that the map $\tilde{L}$ is adjointable was wrong. Altough $\tilde{L}$ acts precisely as $L$, it needs not be adjointable.

Alternatively, one can try to extend $L$ to an adjointable map $L^{(e)}:\ell^2(\tilde{\textsf{A}}) \rightarrow \tilde{\textsf{A}}$. If this is possible, $L$ will be represented by the sequence $(L^{(e)})^*e\in\ell^2(\tilde{\textsf{A}})$. But again, this fails in general. Observe that the representing sequence $(e_n)_n$ in the preceding example does not belong to $\ell^2(\tilde{\textsf{A}})$. This is simply because the series $\sum_{n=1}^{\infty}e_ne_n^*=\sum_{n=1}^{\infty}e_n$ does not converge in norm.

The point is that there is no easy way to extend $L$ to an adjointable map of modules over unital $C^*$-algebras. In order to do so, one has to go all the way up to the maximal essential extension $M(\ell^2(\textsf{A}))$ of $\ell^2(\textsf{A})$. Only there, extensions of adjointable operators are available (as we mentioned in the discussion preceding Lemma \ref{adjointable dual}) thanks to the strict continuity.
\end{remark}

\vspace{.1in}

\begin{example}\label{drugi primjer}
Take a separable Hilbert space $H$. Let $H_1$ be a closed subspace od $H$ such that both $H_1$ and $H_1^{\perp}$ are infinite-dimensional. Denote by $p\in \Bbb B(H)$ the orthogonal projection to $H_1$. Consider the $C^*$-algebra $\textsf{A}=C^*(\Bbb K(H),p)$ generated by $p$ and all compact operators. Clearly, $\textsf{A}$ acts non-degenerately on $H$.

Let us now take an orthonormal basis $(\varepsilon_n)_n$ for $H_1$. Again for each $n\in \Bbb N$ we denote  by $e_n$ the orthogonal projection to $\text{span}\,\{\varepsilon_n\}$.

Consider the sequence $(e_1,e_2,e_3,\ldots)$. We first observe that all $e_n$ belong to $LM(\textsf{A})$ (in fact, they belong to $\textsf{A}$) and that the sequence $(\sum_{n=1}^{N}e_ne_n^*)_N=(\sum_{n=1}^{N}e_n)_N$ is bounded. Applying Lemma \ref{Heuser} we conclude that $L:\ell^2(\textsf{A})\rightarrow \textsf{A}$, defined by $L((a_n)_n)=\sum_{n=1}^{\infty}a_ne_n$ is bounded. Thus, $L\in \ell^2(\textsf{A})^{\prime}$.

We claim that $L$ is not adjointable, i.e. $L\in \ell^2(\textsf{A})^{\prime}\setminus \ell^2(\textsf{A})^*$. To see this, it suffices by Lemma \ref{adjointable dual} to prove that
$(e_n)_n \not \in M(\ell^2(\textsf{A}))$. Thus, we must show that the series $\sum_{n=1}^{\infty}e_ne_n^*=\sum_{n=1}^{\infty}e_n$ does not converge in the strict topology on $B(H)$ induced by our $C^*$-algebra $\textsf{A}$. But this is clear. The strict convergence fails because of the presence of $p$ in $\textsf{A}$.
Indeed,
$(\sum_{n=1}^Ne_np)_N$ cannot converge in norm since $\sum_{n=1}^Ne_np=\sum_{n=1}^Ne_n$ for all $N$ and we know that $\sum_{n=1}^{\infty}e_n$ is not norm-convergent.
\end{example}

\vspace{.1in}

\begin{remark}\label{konacni nizovi fail}
Let us also note the following observation. If bounded module maps $L_1$ and $L_2$ from $\ell^2(\textsf{A})^{\prime}$ coincide on the set $c_{00}(\textsf{A})$ of all finite sequences (all but finitely many components are equal to zero), then we must have $L_1=L_2$ simply because the set of finite sequences is norm-dense in $\ell^2(\textsf{A})$.

On the other hand, the set of all finite sequences is not norm-dense in $\ell^2_{\text{strong}}(LM(\textsf{A}))$ (and, in particular, it cannot be dense in $\ell^2_{\text{strong}}(\overline{\textsf{A}}^s)$). Recall that $M(\ell^2(\textsf{A}))$ is closed in $\ell^2_{\text{strong}}(\overline{\textsf{A}}^s)$ and observe that all finite sequences obviously belong to $M(\ell^2(\textsf{A}))$; in fact they are already in $\ell^2(\textsf{A})$.

The sequence $x=(e_1,e_2,e_3,\ldots)$ of one-dimensional orthogonal projections from Example \ref{drugi primjer} may serve as a simple illustration. Recall that this was the representing sequence for a module map $L\in \ell^2(\textsf{A})^{\prime}$. So, we know that
$x=(e_1,e_2,e_3,\ldots)$ must be in $\ell^2_{\text{strong}}(LM(\textsf{A}))$ (which is also easily seen by a direct verification). For $N$ in $\Bbb N$ let $x_N=(e_1,e_2,\ldots,e_N,0,0,\ldots)$. Clearly, the sequence $(x_N)_N$ cannot converge to $x$ in the Hilbert $C^*$-norm on $\ell^2_{\text{strong}}(LM(\textsf{A}))$. In fact, if $N< M$ we have $\|x_M-x_N\|^2=\left\|\langle x_M-x_N, x_M-x_N\rangle\right\|=\left\|\sum_{n=N+1}^Me_n \right\|= 1$.
\end{remark}

\vspace{.3in}

\section{Weak frames}

\vspace{.1in}

Throughout this section $X$ will denote a left Hilbert $C^*$-module over a von Neumann algebra $\textsf{A}$ acting on a Hilbert space $H$.  We point out that the hypothesis that $\textsf{A}$ is a von Neumann algebra is essential.

Recall from the preceding section that, when $\textsf{A}$ is a von Neumann algebra, $\ell^2(\textsf{A})$ - the standard Hilbert $C^*$-module over $\textsf{A}$ - is contained as a closed submodule in its dual
$$
\ell^2_{\text{strong}}(\textsf{A})=\left\{ (a_n)_n: a_n\in \textsf{A},\,\,\sup_{N}\left\|\sum_{n=1}^Na_na_n^*\right\|<\infty\right\},
$$
which is a self-dual Hilbert $\textsf{A}$-module with the inner product
$$\langle (a_n)_n,(b_n)_n\rangle=\text{(strong)}\,\sum_{n=1}^{\infty}a_nb_n^*.$$
that extends the inner product defined on $\ell^2(\textsf{A})$.


\vspace{.1in}

In the first part of this section we additionally assume that $X$ is self-dual.

\vspace{.1in}

\begin{definition}\label{w-s}
Let $X$ be a self-dual Hilbert $C^*$-module over a von Neumann algebra \textsf{A}. The weak-strong topology on $X$ is defined as the weak topology induced by the maps $L_y :X \rightarrow \textsf{A}$, $L_y(x) = \langle x,y\rangle$, $y\in X$, where $\textsf{A}\subseteq \Bbb B(H)$ is regarded with the respect to the strong operator topology.
\end{definition}

A net $(x_{\lambda})_{\lambda}$ in $X$ converges $\text{weak-strong}$ to $x\in X$, which we denote as $x=(\text{weak-strong})\,\lim_{\lambda}x_{\lambda}$, if and only if
$$
\langle x,y\rangle = (\text{strong})\,\lim_{\lambda}\langle x_{\lambda},y\rangle,\,\,\,\,\forall y\in X.
$$

Here again we use our assumption that \textsf{A} is a von Neumann algebra: this guarantees that $(\text{strong})\,\lim_{\lambda}\langle x_{\lambda},y\rangle$, if it exists, belongs to \textsf{A}, so it makes sense to require that this limit is equal to $\langle x,y\rangle$ which is, by the definition of a Hilbert \textsf{A}-module, an element of \textsf{A}.

It should also be observed that the weak-strong topology is Hausdorff. Namely, the strong operator topology on \textsf{A} is Hausdorff and the family $L_y$, $y\in X$, obviously separates points of $X$.

\vspace{.2in}



\begin{definition}\label{very first fd}
Let $X$ be a self-dual Hilbert $C^*$-module over a von Neumann algebra \textsf{A}. A sequence $(x_n)_n$ in $X$ is called a \emph{weak frame} for $X$ if there exist positive constants $A$ and $B$ such that
\begin{equation}\label{weak frame definition}
A \langle x,x\rangle \leq (\text{strong})\,\sum_{n=1}^{\infty}\langle x,x_n\rangle \langle x,x_n\rangle^* \leq B \langle x,x\rangle,\quad\forall x\in X.
\end{equation}
A sequence for which the second inequality in \eqref{weak frame definition} is satisfied for some constant $B$ is said to be a \emph{weak Bessel sequence}.
The constants $A$ and $B$ are called \emph{frame bounds}. If $A=B=1,$ i.e.,  if
\begin{equation}\label{Parseval}
(\text{strong})\,\sum_{n=1}^{\infty}\langle x,x_n\rangle \langle x,x_n^*\rangle = \langle x,x\rangle,\quad\forall x\in X,
\end{equation}
the sequence $(x_n)_n$ is called a \emph{weak Parseval frame}.
\end{definition}

\vspace{.1in}

Note that it is implicitly required in the above definition that the series in \eqref{weak frame definition} converges strongly in \textsf{A}. We know that this is equivalent to the condition $\sup_{N}\left\|\sum_{n=1}^N\langle x,x_n\rangle \langle x,x_n\rangle^*\right\|<\infty$ (this is already noted in establishing equivalence of \eqref{glavni modul} and \eqref{glavni modul 1}). Recall from \cite{FL2} that in the definition of a standard frame one requires norm-convergence of the series $\sum_{n=1}^{\infty}\langle x,x_n\rangle\langle x,x_n^* \rangle$. Hence, each standard frame in $X$ is a weak frame.

However, there are weak frames that are not standard.

\vspace{.1in}

\begin{example}\label{jeftin}
Let $H$ be an infinite-dimensional separable Hilbert space with an orthonormal basis $(\varepsilon_n)_{n\in \Bbb N}$. For each $n$ in $\Bbb N$ denote  by $e_n$ the one-dimensional projection to $\text{span}\{\varepsilon_n\}$. Clearly, we have $(\text{strong}) \sum_{n=1}^{\infty}e_n=e$, where $e$ denotes the identity operator on $H$.

Consider $\Bbb B(H)$ - the algebra of all bounded linear operators on $H$ - as a Hilbert $C^*$-module over itself. We know that $\Bbb B(H)$ is self-dual (cf.~Theorem 1.5 from \cite{Lin}).
We claim that $(e_n)_n$ is a weak non-standard Parseval frame for $\Bbb B(H)$. Indeed, for each $a\in \Bbb(H)$ the series $\sum_{n=1}^{\infty}\langle a,e_n\rangle \langle a,e_n\rangle^*=
\sum_{n=1}^{\infty}ae_na^*=a\left(\sum_{n=1}^{\infty}e_n\right)a^*$ converges strongly to $aa^*=\langle a,a\rangle$.

On the other hand, this is not a standard frame; that is, the series $\sum_{n=1}^{\infty}\langle a,e_n\rangle \langle a,e_n\rangle^*=
\sum_{n=1}^{\infty}ae_na^*$ cannot converge in norm to $aa^*$ for all $a\in \Bbb B(H)$. Indeed, $aa^*=(\text{norm})\sum_{n=1}^{\infty}ae_na^*$ forces $aa^*$ (and hence $a$) to be a compact operator.
\end{example}

\vspace{.1in}

Our first goal is to introduce analysis and synthesis operators for weak frames. In order to ensure adjointability of these operators, the self-duality assumption is needed. On the other hand, all we need can be done even for weak Bessel sequences. Before stating the theorem we recall our notational convention: for $a \in \textsf{A}$ and $n\in \Bbb N$ we denote by $a^{(n)}\in \ell^2(\textsf{A})$ the sequence with $a$ on the $n$-th place and zeros elsewhere. The unit element in \textsf{A} (i.e.~the identity operator on the underlying Hilbert space $H$) is denoted by $e$.

\vspace{.2in}

\begin{theorem}\label{Bessel relaxed}
Let $X$ be a self-dual Hilbert $C^*$-module $X$ over a von Neumann algebra \textsf{A} and let $(x_n)_n$ be a sequence in $X$. Then the following two conditions are equivalent:
\begin{itemize}
\item[(a)] $(x_n)_n$ is a weak Bessel sequence.
\item[(b)] $\sup_{N}\left\|\sum_{n=1}^N\langle x,x_n\rangle \langle x,x_n\rangle^*\right\|<\infty$ for all $x$ in $X$.
\end{itemize}
If $(x_n)_n$ is a weak Bessel sequence the map
$$U:X \rightarrow \ell^2_{\text{strong}}(\textsf{A}),\,\,\,Ux=(\langle x,x_n\rangle)_n$$
is adjointable.
The adjoint operator $U^*:\ell^2_{\text{strong}}(\textsf{A})\rightarrow X$ is given by
\begin{equation}\label{Bessel-adjoint}
U^*((a_n)_n)=(\text{weak-strong})\,\sum_{n=1}^{\infty}a_nx_n,\,\,\, \forall (a_n)_n\in\ell^2_{\text{strong}}(\textsf{A}).
\end{equation}
In particular,  $U^*e^{(n)}=x_n$ for all $n \in \Bbb N.$
\end{theorem}
\proof
Clearly, (a) implies (b). To prove the reverse implication, suppose that (b) is satisfied.
Then $U: X \to \ell^2_{\text{strong}}(\textsf{A})$, $Ux=(\langle x,x_n\rangle)_n$, is a well defined module map.
One now shows, precisely as it is done in the proof of Theorem 2.1 in \cite{AB}, that $U$ is bounded using the closed graph theorem. For reader's convenience we include the proof.

Let $(y,(a_n)_n)=\lim_{k\rightarrow \infty}(y_k,Uy_k),$ where $y_k,y\in X,\,(a_n)_n\in \ell^2_{\text{strong}}(\textsf{A}).$
For each $m\in \Bbb N$ and all $k\in \Bbb N$ we have
\begin{eqnarray*}
\langle (a_n)_n-Uy_k,(a_n)_n-Uy_k\rangle&=&(\text{strong})\sum_{n=1}^{\infty}(a_n-\langle y_k,x_n\rangle)(a_n-\langle y_k,x_n\rangle)^*\\
 &\geq&(a_m-\langle y_k,x_m\rangle)(a_m-\langle y_k,x_m\rangle)^*.
\end{eqnarray*}
It follows
$$ \|a_m-\langle y_k,x_m\rangle\|\leq \|(a_n)_n-Uy_k\|.$$
Since by the assumption we have $(a_n)_n=\lim_{k\rightarrow\infty}Uy_k$ and $y=\lim_{k\rightarrow\infty}y_k,$ this implies
$$
a_m=\lim_{k\rightarrow\infty}\langle y_k,x_m\rangle=\langle y,x_m\rangle.
$$
As $m$ was arbitrary, this shows that $(a_n)_n=Uy.$ So, the graph of $U$ is closed.

Knowing that $U$ is a bounded module map we can apply Theorem 2.8 in \cite{P} to conclude that $\langle Ux,Ux\rangle\leq \|U\|^2\langle x,x\rangle$ for all $x\in X$; thus, $(x_n)_n$ is a weak Bessel sequence. In this way we have proved that (b) implies (a).

Moreover, since $X$ is self-dual, by Proposition 3.4 from \cite{P} we conclude that $U$ is in fact an adjointable operator.

It remains to obtain \eqref{Bessel-adjoint}. We need to prove, for each $(a_n)_n\in\ell^2_{\text{strong}}(\textsf{A})$, that the sequence $\left(\sum_{n=1}^Na_nx_n\right)_N$ converges in the weak-strong topology on $X$  to $U^*((a_n)_n)$. Thus, we must show that
$$
(\text{strong})\,\lim_{N\rightarrow \infty} \left\langle \sum_{n=1}^Na_nx_n,y\right\rangle=\left\langle U^*((a_n)_n),y\right\rangle,\,\,\,\forall y \in X.
$$
Take any $y\in X$. Then we have
\begin{eqnarray*}
(\text{strong})\,\lim_{N\rightarrow \infty} \left\langle \sum_{n=1}^Na_nx_n,y\right\rangle&=&(\text{strong})\,\lim_{N\rightarrow \infty} \sum_{n=1}^Na_n\langle x_n,y\rangle \\
 &=&(\text{strong})\sum_{n=1}^{\infty}a_n\langle y,x_n\rangle^*\\
 &=&\left\langle (a_n)_n,Uy\right\rangle\\
  &=&\left\langle U^*((a_n)_n),y\right\rangle.
\end{eqnarray*}

Finally, the equality $U^*e^{(n)}=x_n$ for all $n$ in $\Bbb N$ now follows from $$\langle x,U^*e^{(n)}\rangle=\langle Ux,e^{(n)}\rangle =\langle x,x_n\rangle,\,\,\,\forall x\in X.$$
\qed

\vspace{.1in}

\begin{definition}
The operators $U$ and $U^{*}$ are called the analysis and the synthesis operator, respectively.
\end{definition}

\vspace{.1in}

\begin{remark}\label{r1}
After concluding in the preceding proof that $U$ is bounded one might try to prove the existence of the adjoint $U^*$ directly, as it is known from literature for standard frames. In the first step one could put $U^*e^{(n)}=x_n$ for every $n$ in $\Bbb N$. And this works fine: for each $x$ in $X$ we would than have
$\langle x,U^*e^{(n)}\rangle = \langle x,x_n\rangle = \langle Ux,e^{(n)}\rangle$. But the next step is precisely where the proof of Theorem 4.1 from \cite{FL2} breaks down in this situation.
The obstruction lies in the fact already observed in Remark \ref{konacni nizovi fail} that the set of all finite sequences is not dense in $\ell^2_{\text{strong}}(\textsf{A})$. And this is the reason why we needed self-duality of $X$ to ensure adjointability of $U$.
\end{remark}

\vspace{.1in}

\begin{remark}\label{r2}
Suppose that $(x_n)_n$ is a weak Bessel sequence in a self-dual Hilbert $C^*$-module over a von Neumann algebra \textsf{A}. Denote by $U$ the analysis operator. We claim that
\begin{equation}\label{Bessel-adjoint ell2}
U^*((a_n)_n)=\text{(norm)}\sum_{n=1}^{\infty}a_nx_n,\,\,\, \forall (a_n)_n\in\ell^2(\textsf{A}).
\end{equation}
To see this we first recall that for each $(a_n)_n\in\ell^2(\textsf{A})$ we have
\begin{equation}\label{norm cvg}
(a_n)_n=\text{(norm)}\lim_{N\rightarrow \infty}\sum_{n=1}^Na_ne^{(n)}=\text{(norm)}\lim_{N \rightarrow \infty}(a_1,a_2,\ldots , a_N,0,0\ldots).
\end{equation}
Since $U^*$ is norm-continuous this gives us  $$U^*(a_n)_n=\text{(norm)}\lim_{N\rightarrow \infty}\sum_{n=1}^Na_nU^*e^{(n)}=\text{(norm)}\sum_{n=1}^{\infty}a_nx_n.$$
In general, we cannot conclude norm-convergence of the series $\sum_{n=1}^{\infty}a_nx_n$ for sequences from $\ell^2_{\text{strong}}(\textsf{A})\setminus \ell^2(\textsf{A})$. This is clearly seen from the following example.

Take $X=\ell^2_{\text{strong}}(\textsf{A})$ and consider the sequence $(e^{(n)})_n$. This sequence is  weak Bessel (in fact, it is a weak Parseval frame) in $X$. Here  norm-convergence of the series
$\sum_{n=1}^{\infty}a_nx_n$, that is, $\sum_{n=1}^{\infty}a_ne^{(n)}$ is in fact equation \eqref{norm cvg} for which we know that it is satisfied  only for those sequences $(a_n)_n\in \ell^2_{\text{strong}}(\textsf{A})$ that belong to $\ell^2(\textsf{A})$.

It is also useful to note that although the inner product in $\ell^2_{\text{strong}}(\textsf{A})$ is defined in terms of strong convergence, we have $\langle (a_n)_n,Ux\rangle =  \text{(norm)}\sum_{n=1}^{\infty}a_n\langle x,x_n\rangle^*$ for all $(a_n)_n\in \ell^2(\textsf{A})$ and $x\in X$. This is because $Ux$, being an element of $\ell^2_{\text{strong}}(\textsf{A})$, induces a bounded module map $L_{Ux}:\ell^2(\textsf{A}) \rightarrow \textsf{A}$ and by Lemma \ref{Heuser} the series
$\sum_{n=1}^{\infty}a_n\langle x,x_n\rangle^*$ is norm-convergent whenever is $(a_n)_n$ an element of $\ell^2(\textsf{A})$.
\end{remark}

\vspace{.1in}

\begin{remark}\label{prva baza}
In fact, the sequence $(e^{(n)})_n$ from  Remark \ref{r2} is more than a weak Parseval frame for $\ell^2_{\text{strong}}(\textsf{A})$. Note that we have $\langle e^{(n)},e^{(m)}\rangle =\delta_{n,m} e$ for all $n,m\in \Bbb N$ which shows that this is an orthogonal sequence. Since it is a weak Parseval frame, we also have $U^*U=I$, i.e.
$$
(c_n)_n=(\text{weak-strong})\,\sum_{n=1}^{\infty}c_ne^{(n)},\,\,\, \forall (c_n)_n\in\ell^2_{\text{strong}}(\textsf{A}).
$$
Thus, $(e^{(n)})_n$ is what can be called a weak orthonormal basis for $\ell^2_{\text{strong}}(\textsf{A})$. Since its analysis operator coincides with the identity operator,  we shall refer in the sequel to $(e^{(n)})_n$ as to the canonical weak basis for $\ell^2_{\text{strong}}(\textsf{A})$. (Observe that $(e^{(n)})_n$ is in the same time a standard orthonormal basis for $\ell^2(\textsf{A})$.) In general, one can define weak Riesz bases as those weak frames whose analysis operators are invertible. In the present paper we omit further discussion on weak Riesz bases.
\end{remark}

\vspace{.1in}

\begin{prop}\label{bessel corr adj}
Let $X$ be a self-dual Hilbert $C^*$-module  over a von Neumann algebra \textsf{A}.
Let $T:\ell^2_{\text{strong}}(\textsf{A}) \rightarrow X$ be a bounded module map. Then the sequence $(x_n)_n$, where $x_n=Te^{(n)}$, $n\in \Bbb N$, is a weak Bessel sequence in $X$ whose analysis operator coincides with $T^*$.
\end{prop}
\proof
First observe that since $\ell^2_{\text{strong}}(\textsf{A})$ is self-dual, $T$ is adjointable by Proposition 3.4 from \cite{P}. We now have
for all $x$ in $X$ and $n$ in $\Bbb N$
$$
\langle x,x_n\rangle = \langle x,Te^{(n)}\rangle = \langle T^*x,e^{(n)}\rangle.
$$
The last term is precisely the $n$-th component of the sequence $T^*x$; hence, $T^*x=(\langle x,x_n\rangle)_n$. Moreover, as $T^*x$ belongs to $\ell^2_{\text{strong}}(\textsf{A})$, we do have $\sup_{N}\left\|\sum_{n=1}^N\langle x,x_n\rangle\langle x,x_n\rangle^* \right\|<\infty$. By Theorem \ref{Bessel relaxed} this implies that $(x_n)_n$ is a weak Bessel sequence.
\qed

\vspace{.1in}

\begin{remark}\label{bessel vs adjointability}
If $X$ is a self-dual Hilbert $C^*$-module  over a von Neumann algebra \textsf{A}, Theorem \ref{Bessel relaxed} and Proposition \ref{bessel corr adj} show us that weak Bessel sequences in $X$ are in a bijective correspondence with bounded module maps $T:\ell^2_{\text{strong}}(\textsf{A}) \rightarrow X$.
\end{remark}

\vspace{.1in}

The following theorem provides us with a sufficient frame condition which is considerably easier to check than the original condition from Definition \ref{very first fd}.

\vspace{.1in}

\begin{theorem}\label{simple minds}
Let $X$ be a self-dual Hilbert $C^*$-module  over a von Neumann algebra \textsf{A}. Then a sequence $(x_n)_n$ in $X$ is a weak frame for $X$ if and only if there exists a constant $A>0$ such that
\begin{equation}\label{sufficient frame condition}
A\|x\|^2\leq \sup_{N}\left\|\sum_{n=1}^N\langle x,x_n\rangle\langle x,x_n\rangle^* \right\|<\infty,\,\,\,\forall x \in X.
\end{equation}
\end{theorem}
\proof
We have already observed that  the series $\sum_{n=1}^{\infty}\langle x,x_n\rangle\langle x,x_n\rangle^*$ converges strongly if and only if $\sup_{N}\left\|\sum_{n=1}^N\langle x,x_n\rangle\langle x,x_n\rangle^*\right\|<\infty$
in which case we have
$$\left\|(\text{strong})\sum_{n=1}^{\infty}\langle x,x_n\rangle\langle x,x_n\rangle^* \right\| =\sup_{N}\left\|\sum_{n=1}^N\langle x,x_n\rangle\langle x,x_n\rangle^*\right\|.$$
The desired conclusion now follows from our Theorem \ref{Bessel relaxed} and Corollary 2.2 from \cite{A}.
\qed

\vspace{.1in}

\begin{cor}\label{onto}
Let $X$ be a self-dual Hilbert $C^*$-module  over a von Neumann algebra \textsf{A}. A sequence $(x_n)_n$ in $X$ is a weak frame for $X$ if and only if there is a surjective bounded module map
$T:\ell^2_{\text{strong}}(\textsf{A}) \rightarrow X$ such that $Te^{(n)}=x_n$ for each $n$ in $\Bbb N$.

If $(x_n)_n$ is a weak frame in $X$ and if $R:X\rightarrow Y$ is a surjective bounded module map, where $Y$ is a self-dual Hilbert \textsf{A}-module, then the sequence
$(Rx_n)_n$ is a weak frame in $Y$.
\end{cor}
\proof
Suppose that $(x_n)_n$ is a weak frame in $X$. Then by \eqref{sufficient frame condition} the analysis operator $U$ is bounded from below and hence has a closed range. This implies that $U^*$ is a surjection.

To prove the converse, suppose that we are given a surjective bounded module map $T:\ell^2_{\text{strong}}(\textsf{A}) \rightarrow X$. First, $T$ is adjointable by Proposition 3.4 from \cite{P}. By Theorem 3.2 from \cite{L} $T^*$ is injective and has a closed range which implies that $T^*$ is bounded from below. As in the proof of Proposition \ref{bessel corr adj} we conclude that $T^*$ is in fact the analysis operator of the sequence $(x_n)_n$. Now Theorem \ref{simple minds} applies.

The second assertion of the corollary is an immediate consequence of the first one.
\qed

\vspace{.1in}

\begin{cor}\label{bezuvjetno}
Let $X$ be a self-dual Hilbert $C^*$-module  over a von Neumann algebra \textsf{A}. If $(x_n)_n$ is a weak frame in $X$, then the sequence $(x_{\sigma(n)})_n$ is also a weak frame  for every permutation $\sigma$ of the set $\Bbb N$.
\end{cor}
\proof
Take any permutation $\sigma$ and fix $N_1\in \Bbb N$.  Then we have
$$
\left\|\sum_{n=1}^{N_1}\langle x,x_{\sigma(n)}\rangle\langle x,x_{\sigma(n)}\rangle^* \right\| \leq \sup_{N}\left\|\sum_{n=1}^N\langle x,x_n\rangle \langle x,x_n\rangle^*\right\|
$$
and
$$
\left\|\sum_{n=1}^{N_1}\langle x,x_n\rangle \langle x,x_n\rangle^*\right\| \leq \sup_{N}\left\|\sum_{n=1}^N\langle x,x_{\sigma(n)}\rangle\langle x,x_{\sigma(n)}\rangle^*\right\|.
$$
This is enough to conclude that
$$
\sup_{N}\left\|\sum_{n=1}^N\langle x,x_{\sigma(n)}\rangle\langle x,x_{\sigma(n)}\rangle^*\right\|=\sup_{N}\left\|\sum_{n=1}^N\langle x,x_n\rangle \langle x,x_n\rangle^*\right\|.
$$
The conclusion now follows from Theorem \ref{simple minds}.
\qed

\vspace{.1in}

\begin{remark}\label{bezuvjetno a}
The conclusion of the preceding corollary applies also to Bessel sequences. This in turn implies that when working with the synthesis operator of a weak Bessel sequence (or a weak frame) $(x_n)_n$ the relevant series $\sum_{n=1}^{\infty}c_nx_n$ converges unconditionally in the weak-strong topology for every sequence $(c_n)_n\in \ell^2_{\text{strong}}(\textsf{A})$.

This is essential in applications since we will often encounter countable systems naturally indexed by sets different from the set of natural numbers. In particular, when dealing with systems indexed by the set of the integers $\Bbb Z$, we will be allowed to restrict our considerations to symmetric partial sums  $\sum_{n=-N}^Nc_nx_n$, $N\in \Bbb N$.
\end{remark}

\vspace{.1in}

\begin{theorem}\label{rf}
Let $X$ be a self-dual Hilbert $C^*$-module  over a von Neumann algebra \textsf{A}. Suppose that $(x_n)_n$ is a weak frame in $X$ with the analysis operator $U$. Then the frame operator $U^*U$ is invertible, the sequence $((U^*U)^{-1}x_n)_n$ is a weak frame in $X$, and the following reconstruction formula holds:
\begin{equation}\label{reconstruction formula}
x=
(\text{weak-strong})\,\sum_{n=1}^{\infty}\langle x,x_n\rangle (U^*U)^{-1}x_n,\,\,\, \forall x\in X.
\end{equation}
\end{theorem}
\proof
As $U$ is bounded from below and has a closed range, one shows that $U^*U$ is invertible by a standard argument. Corollary \ref{onto} tells us that $((U^*U)^{-1}x_n)_n$ is a weak frame in $X$. To prove the reconstruction formula first observe that each adjointable operator $R:X \rightarrow X$ is continuous in the weak-strong topology.
Indeed, suppose that $x=(\text{weak-strong})\lim_{\lambda}x_{\lambda}$ and fix any $y\in X$. Then we have
$$
\langle Rx,y\rangle=\langle x,R^*y\rangle = (\text{strong})\lim_{\lambda}\langle x_{\lambda},R^*y\rangle=(\text{strong})\lim_{\lambda}\langle Rx_{\lambda},y\rangle
$$
which shows that $Rx=(\text{weak-strong})\lim_{\lambda}Rx_{\lambda}$.
Now we have
$$
U^*Ux\stackrel{\eqref{Bessel-adjoint}}{=}(\text{weak-strong})\,\sum_{n=1}^{\infty}\langle x,x_n\rangle x_n,\,\,\, \forall x \in X.
$$
Applying $(U^*U)^{-1}$ to both sides and using its weak-strong continuity we get
$$
x=
(\text{weak-strong})\,\sum_{n=1}^{\infty}\langle x,x_n\rangle (U^*U)^{-1}x_n,\,\,\, \forall x\in X.
$$
\qed

\vspace{.2in}

\begin{definition}
The weak frame $((U^*U)^{-1}x_n)_n$ from the preceding theorem is called the canonical dual of $(x_n)_n$.
\end{definition}

\vspace{.2in}

\begin{remark}\label{mutual duals}
Denote by $V$ the analysis operator of the canonical dual $((U^*U)^{-1}x_n)_n$. Then the reconstruction formula \eqref{reconstruction formula} simply reads $V^*U=I$ where $I$ is the identity operator on $X$. But this is obviously equivalent to $U^*V=I$ which means that we also have
\begin{equation}\label{reconstruction formula bis}
x=
(\text{weak-strong})\,\sum_{n=1}^{\infty}\langle x,(U^*U)^{-1}x_n\rangle x_n,\,\,\, \forall x\in X.
\end{equation}
This tells us that $((U^*U)^{-1}x_n)_n$ and $(x_n)_n$ are dual to each other in a symmetric way.

By a standard argument one also shows that $((U^*U)^{-\frac{1}{2}}x_n)_n$ is a weak Parseval frame. This follows from the the equality
\begin{equation}\label{reconstruction formula parseval}
x=
(\text{weak-strong})\,\sum_{n=1}^{\infty}\langle x,(U^*U)^{-\frac{1}{2}}x_n\rangle (U^*U)^{-\frac{1}{2}}x_n,\,\,\, \forall x\in X
\end{equation}
and the fact that a sequence $(y_n)_n$ in $X$ is a weak Parseval frame if and only if it has the property
\begin{equation}\label{reconstruction formula parseval a}
x=
(\text{weak-strong})\,\sum_{n=1}^{\infty}\langle x,y_n\rangle y_n,\,\,\, \forall x\in X.
\end{equation}
\end{remark}

\vspace{.2in}

Our next proposition basically says that unitary operators of Hilbert $C^*$-modules map weak frames into weak frames. This is something that is certainly expected, but one should observe that in the proof that follows we again encounter a step for which we need  the assumption that the underlying $C^*$-algebra is a von Neumann algebra (i.e. any isomorphism of von Neumann algebras is normal; \cite{Ped}, Proposition 2.5.2).

Suppose that $X$ and $Y$ are left Hilbert $C^*$-modules over $C^*$-algebras $\textsf{A}$ and $\textsf{B}$, respectively, and that $\varphi : \textsf{A} \rightarrow \textsf{B}$ is a morphism of $C^*$-algebras. Recall from \cite{BG} that a map $\Phi :X \rightarrow Y$ is said to be a $\varphi$-morphism of $X$ and $Y$ if it satisfies $\langle \Phi(x_1),\Phi(x_2)\rangle =\varphi(\langle x_1,x_2\rangle)$ for all $x_1,x_2$ from $X$. It turns out that such a map necessarily satisfies $\Phi(ax)=\varphi(a)\phi(x)$ for all $x\in X$ and $a\in \textsf{A}$. We say that $\Phi$ is a unitary operator of Hilbert $C^*$-modules if both $\Phi$ and $\varphi$ are bijections (in fact, it suffices to require that $\Phi$ is a surjection and that $\varphi$ is injective). When this is the case, both $\Phi$ and $\varphi$ are ismotries.

\vspace{.1in}

\begin{prop}\label{slika slabog framea}
Suppose that $X$ and $Y$ are self-dual Hilbert $C^*$-modules over von Neumann algebras $\textsf{A}$ and $\textsf{B}$, respectively, that $\varphi : \textsf{A} \rightarrow \textsf{B}$ is an isomorphism and that $\Phi :X \rightarrow Y$ is a $\varphi$-unitary operator. Then a sequence $(x_n)_n$ is a weak frame (resp.~Bessel sequence) in $X$ with bounds $A$ and $B$ if and only if the sequence $(\Phi(x_n))_n$ is a weak frame (resp.~Bessel sequence) in $Y$ with the same frame (resp.~Bessel) bounds.
\end{prop}
\proof
Essentially this follows from two important properties of isomorphisms of von Neumann algebras. First, we now that $a\leq b$ implies $\varphi(a)\leq \varphi(a)$ (this is already true when  $\varphi$ is a morphism of arbitrary $C^*$-algebras) and secondly, that $\varphi$ is normal. This later property means that when $(a_n)_n$ is an increasing strongly convergent net in $\textsf{A}$ such that $a=\text{(strong)}\lim_{n} a_n$, then $\varphi(a)=\text{(strong)}\lim_{n} \varphi(a_n)$.

Suppose now we are given a weak frame $(x_n)_n$ in $X$ with frame bounds $A$ and $B$. Then we have
$$
A\langle x,x\rangle \leq \text{(strong)} \sum_{n=1}^{\infty}\langle x,x_n\rangle\langle x,x_n\rangle^* \leq A\langle x,x\rangle,\,\,\,\forall x \in X.
$$
We now use all the properties of $\varphi$ to obtain
$$
A\varphi(\langle x,x\rangle) \leq \text{(strong)} \sum_{n=1}^{\infty}\varphi(\langle x,x_n\rangle)\varphi(\langle x,x_n\rangle)^*\leq B\varphi(\langle x,x\rangle),\,\,\,\forall x \in X.
$$
Since $\Phi$ is a $\varphi$-morphism this gives us
$$
A\langle \Phi(x),\Phi(x)\rangle \leq \text{(strong)} \sum_{n=1}^{\infty}\langle \Phi(x),\Phi(x_n)\rangle \langle \Phi(x),\Phi(x_n)\rangle^*\leq B\langle \Phi(x),\Phi(x)\rangle,\,\,\,\forall x \in X.
$$
Finally, since $\Phi$ is a surjection, this tells us that $(\Phi(x_n))_n$ is a weak frame in $Y$.
The same reasoning applied to the maps $\Phi^{-1}$ and $\varphi^{-1}$ proves the converse.
\qed

\vspace{.2in}

We end the section with a brief discussion on a more general situation in which we do not assume self-duality of $X$.
Suppose that $X$ is an arbitrary Hilbert $C^*$-module over a von Neumann algebra $A$. As it is well known (see \cite{P}), each Hilbert $\textsf{A}$-module $X$ can be embedded into a self-dual $\textsf{A}$-module $X^{\prime}$ in such a way that the inner product on $X^{\prime}$ extends the original inner product given on $X$. In fact, $X^{\prime}$ is with a slight abuse of notation a twisted copy of the dual of $X$. As $X^{\prime}$ is self-dual, all what is said in this section concerning with weak Bessel sequences and weak frames applies to $X^{\prime}$.

However, sometimes is the object of our real interest a Hilbert $C^*$-module $X$ that is not self-dual. Then the question arises: are there weak frames for $X$ in the context of a broader ambient module $X^{\prime}$?
Certainly, weak frames for $X$ can be obtained by working in the dual module $X^{\prime}$. In particular, for all weak frames in $X^{\prime}$ the reconstruction formula \eqref{reconstruction formula} is still valid and applies, in particular, to all elements from $X$. The difference is now that such frames may be outer from the $X$-perspective in the sense that some of the frame elements $x_n$ (or even all of them) may belong to $X^{\prime}\setminus X$.

\vspace{.1in}

On the other hand, it is natural to ask what can be said about standard frames for $X$ in this more general context. To provide the answer we first need a couple of auxiliary results which are obtained in (or can be easily concluded from) \cite{P}.

Our first lemma is proved in \cite{P}.  Since it is not explicitly stated there, we include it here for future reference.

\vspace{.1in}

\begin{lemma}\label{trivial complement}\emph{(}\cite{P}\emph{)}
Let $X$ be a Hilbert $C^*$-module over a von Neumann algebra \textsf{A}. Suppose that $f\in X^{\prime}$ is orthogonal to all $x$ in $X$. Then $f=0$. (In other words, the orthogonal complement of $X$ in $X^{\prime}$ is trivial.)
\end{lemma}
\proof
We have $\langle x,f\rangle=0$ for all $x$ in $X$. Then  by Theorem 3.2 from \cite{P} (when $f$ is regarded as a bounded module map from $X$ in \textsf{A}) we have $f(x)=0,\,\forall x\in X$. By the last paragraph of the proof of Theorem 3.2 in \cite{P} this forces $f=0$.
\qed

\vspace{.1in}

We now recall again Propositions 3.4 and 3.6 from \cite{P}. By Proposition 3.6 in \cite{P} if $X$ and $Y$ are Hilbert $C^*$-modules over a von Neumann algebra $\textsf{A}$, each bounded module map $T:X\rightarrow Y$ extends uniquely to a bounded module map $T^{\prime}:X^{\prime} \rightarrow Y^{\prime}$.  (Note that this provides us with an alternative argument for the proof of Lemma \ref {trivial complement}.) Let us call this extended map $T^{\prime}$ the standard extension of $T$.

Another fact that we need is  Proposition 3.4 from  \cite{P}: if $X$ and $Y$ are self-dual Hilbert $C^*$-modules over a von Neumann algebra $\textsf{A}$, each bounded module map $T:X\rightarrow Y$ is adjointable. (In fact, for this conclusion we do not need self-duality of $Y$.)

\vspace{.1in}

\begin{lemma}\label{zic}
Suppose that $X$ and $Y$ are Hilbert $C^*$-modules over a von Neumann algebra $\textsf{A}$ and that $T:X\rightarrow Y$ is an adjointable operator. Then the standard extensions $T^{\prime}:X^{\prime} \rightarrow Y^{\prime}$ and
$(T^*)^{\prime}:Y^{\prime} \rightarrow X^{\prime}$ are adjoint to each other. In other words, $T^{\prime}$ is an adjointable operator and $(T^{\prime})^*=(T^*)^{\prime}$.
\end{lemma}
\proof
To show $(T^{\prime})^*=(T^*)^{\prime}$ it suffices to prove that these operators coincide on $Y$ (since the standard extension is the unique extension). Take any $y\in Y$. We want to prove that $(T^{\prime})^*y=(T^*)^{\prime}y$. As both sides of this equality belong to $X^{\prime}$, we can apply Lemma \ref{trivial complement}. Hence, it is enough to show that $\left\langle (T^{\prime})^*y-(T^*)^{\prime}y,x\right\rangle=0$ for all $x$ in $X$. Fix $x\in X$. Then
$$
\left\langle (T^{\prime})^*y-(T^*)^{\prime}y,x\right\rangle=\left\langle y,T^{\prime}x\right\rangle-\left\langle T^*y,x \right\rangle
=\left\langle y,Tx\right\rangle - \left\langle y,Tx\right\rangle =0
$$

\qed

\vspace{.1in}

\begin{lemma}\label{inv}
Suppose that $X$ is a Hilbert $C^*$-module over a von Neumann algebra $\textsf{A}$ and that $T:X\rightarrow X$ is an invertible bounded module map. Then its standard extension
$T^{\prime}:X^{\prime} \rightarrow X^{\prime}$ is also invertible.
\end{lemma}
\proof
Observe that $T^{-1}$ is also a bounded module map. Consider the standard extensions $T^{\prime}$ and $(T^{-1})^{\prime}$. For each $x$ in $X$ we now have
$(T^{-1})^{\prime}T^{\prime}x=(T^{-1})^{\prime}Tx=T^{-1}Tx=x$. Thus, $(T^{-1})^{\prime}T^{\prime}$ extends the identity operator $I_X$ on $X$. Since the standard extension is unique, and since the identity operator $I_{X^{\prime}}$ on $X^{\prime}$ obviously extends $I_X$, this forces $(T^{-1})^{\prime}T^{\prime}=I_{X^{\prime}}$.

Precisely in the same way we conclude $T^{\prime}(T^{-1})^{\prime}=I_{X^{\prime}}$. Thus, $T^{\prime}$ is invertible.
\qed

\vspace{.1in}

\begin{prop}\label{standard extension}
Let $(x_n)_n$ be a standard frame in a Hilbert $C^*$-module $X$ over a von Neumann algebra $\textsf{A}$. Then $(x_n)_n$ is a weak frame for $X^{\prime}$.
\end{prop}
\proof
Denote by $U:X\rightarrow \ell^2(\textsf{A})$ the analysis operator of $(x_n)_n$. We know from \cite{FL2} that $U$ is adjointable, that we have $U^*e^{(n)}=x_n$ for all $n$, and that $U^*U : X \rightarrow X$ is an invertible operator.
Recall that $\ell^2(\textsf{A})^{\prime}=\ell^2_{\text{strong}}(\textsf{A})$. Consider the standard extensions $U^{\prime}:X^{\prime}\rightarrow \ell^2_{\text{strong}}(\textsf{A})$ and $(U^*)^{\prime}:\ell^2_{\text{strong}}(\textsf{A}) \rightarrow X^{\prime}$.

By Lemma \ref{zic} we have $(U^*)^{\prime}=(U^{\prime})^*$. So $(U^{\prime})^*U^{\prime}=(U^*)^{\prime}U^{\prime}$ is the standard extension of $U^*U$. Since $U^*U$ is an invertible operator, we conclude from Lemma \ref{inv} that $(U^{\prime})^*U^{\prime}$ is invertible. In particular, $(U^{\prime})^*$ must be surjective. So, $(U^{\prime})^*$ is an adjointable surjection. Recall from Remark \ref{r2} that $(e^{(n)})_n$ is a weak frame for $\ell^2_{\text{strong}}(\textsf{A})$. By applying the second assertion of Corollary \ref{onto} we conclude that $((U^{\prime})^*e^{(n)})_n=((U^*)^{\prime}e^{(n)})_n=
(U^*e^{(n)})_n=(x_n)_n$ is a weak frame for $X^{\prime}$.
\qed

\vspace{.1in}

\begin{remark}\label{standard and weak bessel}
Note that $(x_n)_n$ will be a weak Parseval frame in $X^{\prime}$ if $(x_n)_n$ is a standard Parseval frame in $X$. The proof also shows: if $(x_n)_n$ is a standard Bessel sequence in $X$ then $(x_n)_n$ is a weak Bessel sequence in $X^{\prime}$.
\end{remark}

\vspace{.1in}

\begin{prop}\label{kad je weak dosao od standardnog}
Let $X$ be a Hilbert $C^*$-module over a von Neumann algebra $\textsf{A}$. Let $(x_n)_n$ be a weak frame (resp.~Bessel sequence) in $X^{\prime}$ such that $x_n\in X$ for all $n\in \Bbb N$. Let $Y=\{x\in X:(\langle x,x_n\rangle)_n\in \ell^2(\textsf{A})\}$. If $Y$ is dense in $X$ then $(x_n)_n$ is a standard frame (resp.~Bessel sequence) in $X$.
\end{prop}
\proof
Suppose that $(x_n)_n$ is a weak Bessel sequence and denote by $U:X^{\prime}\rightarrow \ell^2_{\text{strong}}(\textsf{A})$ its analysis operator.
We know from Remark \ref{r2} that
$$
U^*((a_n)_n)=\text{(norm)}\sum_{n=1}^{\infty}a_nx_n,\,\,\, \forall (a_n)_n\in\ell^2(\textsf{A}).
$$
Note that
$$
Y=\{x\in X^{\prime}:Ux\in \ell^2(\textsf{A})\}\cap X.
$$
Clearly, $Y$ is a submodule of $X$. We also conclude that $Y$ is closed because $U$ is adjointable and hence a bounded operator. By  the assumption   $Y$ is dense in $X$ and therefore $Y=X$. Thus, we have $Ux\in \ell^2(\textsf{A})$ for all $x$ in $X$ and, since each $x_n$ is in $X$, this is enough to conclude that $(x_n)_n$ is a standard Bessel sequence in $X$ (see Theorem 2.1 in \cite{AB}).

If, in addition, we assume that $(x_n)_n$ is a weak frame (not merely a weak Bessel sequence) then it follows that $(x_n)_n$ is a standard frame in $X$. In fact, its analysis operator is $U_0=U|_X$ and it is bounded from below since $U$ is bounded from below.
\qed

\vspace{.1in}

In the light of  Proposition \ref{standard extension}, Remark \ref{standard and weak bessel} and Proposition \ref{kad je weak dosao od standardnog} one may ask whether there are weak frames (resp.~weak Bessel sequences) for $X^{\prime}$ that are not standard Bessel sequences or frames for $X$. Indeed, there are.

\vspace{.1in}

\begin{example}\label{beskonacni projektori}
Consider an infinite-dimensional separable Hilbert space $H$ decomposed as $H=\oplus_{n=1}^{\infty}H_n$, where $\text{dim}\,H_n=\infty$ for all $n$. Observe that the elements of $H$ can be identified as sequences $(\xi_n)_n$, $\xi_n \in H_n$, satisfying $\sum_{n=1}^{\infty}\|\xi_n\|^2<\infty$.

Let $X=\Bbb K(H)$ be the Hilbert $\Bbb B(H)$-module  consisting of all compact operators acting on $H$. It is known that $X$ is not self-dual; by Theorem 1.5 from \cite{Lin} we know that $X^{\prime}=\Bbb B(H)$ where $\Bbb B(H)$ is regarded as a Hilbert $C^*$-module over itself.

Denote by $q_n$, $n\in \Bbb N$, the orthogonal projections to $H_n$'s. As in Example \ref{jeftin} one easily shows that the series $\sum_{n=1}^{\infty}q_n$ converges in the strong operator topology to the identity operator.

From this
we conclude for each $x\in \Bbb B(H)$ that $x=(\text{strong})\sum_{n=1}^{\infty}xq_n=(\text{strong})\sum_{n=1}^{\infty}\langle x,q_n\rangle q_n$. In particular, this can be rewritten as in \eqref{reconstruction formula parseval a}:
$x=
(\text{weak-strong})\,\sum_{n=1}^{\infty}\langle x,q_n\rangle q_n$ which means that
the sequence $(q_n)_n$
is a weak Parseval frame for $\Bbb B(H)$. On the other hand, since the range of each $q_n$ is infinite-dimensional, $q_n\not \in \Bbb K(H)$ for all $n$ in $\Bbb N$. Therefore, our weak frame $(q_n)_n$ for the dual $X^{\prime}$ does not arise as a standard frame for the Hilbert $C^*$-module $X$ we started from.
\end{example}

\vspace{.4in}

\section{Bounded operators on $\ell^2(\textsf{A})$ and infinite matrices}

\vspace{.1in}

In this section we discuss conditions on a sequence $(x_n)_n$ sufficient to ensure the weak Bessel property.
Recall from Lemma 3.5.1 in \cite{C} that a sequence $(x_n)_n$ is Bessel in a Hilbert space $H$ if and only if the Gram matrix associated with $(x_n)_n$ defines a bounded operator on $\ell^2$. In the setting of Hilbert $C^*$-modules there is no such result for standard frames since it is not enough to have a bounded operator on $\ell^2(\textsf{A})$; one has additionally to know that this operator is adjointable and this does not follow automatically from boundedness. Here we prove a similar result for weak Bessel sequences in self-dual Hilbert $C^*$-modules over von Neumann algebras. For some technical reasons which will become clear from the context we shall additionally assume that the underlying von Neumann algebra is commutative.

Suppose we have a weak Bessel sequence $(x_n)_n$ with a Bessel bound $B$ in a self-dual Hilbert $C^*$-module $X$ over a von Neumann algebra \textsf{A}. Consider the canonical weak basis $(e^{(n)})_n$ for $\ell^2_{\text{strong}}(\textsf{A})$ from Remark \ref{prva baza}  (which is also a standard basis for $\ell^2(\textsf{A})$). We know that the corresponding analysis operator $U:X \rightarrow \ell^2_{\text{strong}}(\textsf{A})$, $Ux=(\langle x,x_n\rangle)_n$, is an adjointable operator. Hence $UU^*\in \Bbb B(\ell^2_{\text{strong}}(\textsf{A}))$.

Recall from Remark \ref{r2} that the $n$-th component $ UU^*((a_k)_k)_n$ of the sequence $UU^*((a_k)_k)$ is
\begin{equation}\label{gram1}
UU^*((a_k)_k)_n=(\text{norm})\,\sum_{k=1}^{\infty}a_k\langle x_k,x_n\rangle \in \textsf{A},\,\,\,\forall (a_k)_k\in \ell^2(\textsf{A})
\end{equation}
where the norm convergence is ensured by the fact that the sequence $(\langle x_n,x_k\rangle)_k$ belongs to $\ell^2_{\text{strong}}(\textsf{A})$ and hence defines an element from $\ell^2(\textsf{A})^{\prime}$, so Lemma \ref{Heuser} applies. Moreover, since $UU^*((a_k)_k)\in \ell^2_{\text{strong}}(\textsf{A})$, we have from \eqref{gram1}
\begin{equation}\label{gram2}
\sup_N\left\| \sum_{n=1}^N\left( \sum_{k=1}^{\infty}a_k\langle x_k,x_n\rangle \right) \left( \sum_{k=1}^{\infty}a_k\langle x_k,x_n\rangle \right)^*\right\|\leq B^2\|(a_k)_k\|^2.
\end{equation}
Note that \eqref{gram1} implies
\begin{equation}\label{gram3}
\langle UU^*e^{(k)},  e^{(n)}\rangle = (UU^*e^{(k)})_n=\langle x_k,x_n\rangle,\,\,\,\forall n,k.
\end{equation}
Let us now introduce  an infinite matrix $\Gamma$ with entries $\Gamma_{nk}$ in  \textsf{A} defined by
$$
\Gamma_{nk}=\langle x_k,x_n\rangle,\,\,\,n,k\in \Bbb N.
$$
The matrix $\Gamma$ is called the Gramian associated with the sequence $(x_n)_n$. We see from \eqref{gram3} that $\Gamma$ is in fact the matrix representation of $UU^*$ with respect to the canonical basis $(e^{(n)})_n$.

We now suppose that the underlying von Neumann algebra \textsf{A} is commutative. Then
\eqref{gram1} can be rewritten as
\begin{equation}\label{gram1 bis}
UU^*((a_k)_k)_n=(\text{norm})\,\sum_{k=1}^{\infty}\langle x_k,x_n\rangle a_k= (\text{norm})\,\sum_{k=1}^{\infty}\Gamma_{nk} a_k ,\,\,\,\forall (a_k)_k\in \ell^2(\textsf{A})
\end{equation}
which
shows us that $UU^*$ acts on elements $(a_k)_k\in \ell^2(\textsf{A})$ simply as matrix multiplication:
\begin{equation}\label{gram4}
UU^*(a_k)_k=\Gamma\cdot (a_k)_k=\left[\begin{array}{cccc}\langle x_1,x_1\rangle &\langle x_2,x_1\rangle&\langle x_3,x_1\rangle&\ldots\\
\langle x_1,x_2\rangle &\langle x_2,x_2\rangle&\langle x_3,x_2\rangle&\ldots\\
\langle x_1,x_3\rangle &\langle x_2,x_3\rangle&\langle x_3,x_3\rangle&\ldots\\
\vdots&\vdots&\vdots&
\end{array}\right]\left[\begin{array}{c}a_1\\a_2\\a_3\\\vdots\end{array}\right]\in \ell^2_{\text{strong}}(\textsf{A}).
\end{equation}
Now we can reverse the process and ask what can be said about the sequence $(x_n)_n$ if its Gramian defines by \eqref{gram4} a bounded operator on $\ell^2(\textsf{A})$ with values in $\ell^2_{\text{strong}}(\textsf{A})$.

\vspace{.1in}

\begin{prop}\label{Bessel by Gram}
Let $(x_n)_n$ be a sequence in a self-dual Hilbert $C^*$-module over a commutative von Neumann algebra \textsf{A}. Suppose that its Gramian $\Gamma$ defines by $(a_k)_k \mapsto \Gamma\cdot (a_k)_k$ as in \eqref{gram4} a bounded operator  $\ell^2(\textsf{A})\rightarrow \ell^2_{\text{strong}}(\textsf{A})$ with bound $B$. Then $(x_n)_n$ is a weak Bessel  sequence in $X$ with a Bessel bound $B$.
\end{prop}
\proof
By the assumption we have  $\|\Gamma\cdot (a_k)_k\|\leq B \|(a_k)_k\|$ for all $(a_k)_k\in \ell^2(\textsf{A})$, that is,
\begin{equation}\label{gram5}
\sup_N\left\| \sum_{n=1}^N\left( \sum_{k=1}^{\infty}\langle x_k,x_n\rangle a_k \right) \left( \sum_{k=1}^{\infty}\langle x_k,x_n\rangle a_k \right)^*\right\|\leq B^2\|(a_k)_k\|^2,\,\,\,\forall (a_k)_k\in \ell^2(\textsf{A}).
\end{equation}
We follow the proof of Lemma 3.5.1 in \cite{C}. Take any $(a_k)_k\in \ell^2(\textsf{A})$ and natural numbers $N,M$ such that $N<M$. Then we have
\begin{eqnarray*}
\left\|\sum_{n=1}^Ma_nx_n - \sum_{n=1}^Na_nx_n \right\|^4&=&\left\|\sum_{n=N+1}^Ma_nx_n \right\|^4\\
 &=&\left\|\left\langle \sum_{n=N+1}^Ma_nx_n,\sum_{m=N+1}^Ma_mx_m
\right\rangle \right\|^2\\
 &=&\left\| \sum_{m=N+1}^M\left\langle  \sum_{n=N+1}^Ma_nx_n,x_m\right\rangle\ a_m^*\right\|^2\\
  &\leq&\mbox{(applying the Cauchy-Schwarz inequality in } \textsf{A}^{M-N})\\
 &\leq&  \left\| \sum_{m=N+1}^M \left( \sum_{n=N+1}^Ma_n\langle x_n,x_m\rangle \right)\left( \sum_{n=N+1}^Ma_n\langle x_n,x_m\rangle \right)^*\right\|\cdot \left\|\sum_{m=N+1}^M a_ma_m^* \right\|\\
 &\leq&\mbox{(passing in the first term to the sequence } (0,\ldots,0,a_{N+1},\ldots,a_M,0,\ldots))\\
 &=&  \left\| \sum_{m=N+1}^M \left( \sum_{n=1}^{\infty}a_n\langle x_n,x_m\rangle \right)\left( \sum_{n=1}^{\infty}a_n\langle x_n,x_m\rangle \right)^*\right\|\cdot \left\|\sum_{m=N+1}^M a_ma_m^* \right\|\\
 &\stackrel{\eqref{gram5}}{\leq}&B^2 \left\|\sum_{m=N+1}^M a_ma_m^* \right\|\cdot \left\|\sum_{m=N+1}^M a_ma_m^* \right\|.
\end{eqnarray*}
So, we have obtained
$$
\left\|\sum_{n=1}^Ma_nx_n - \sum_{n=1}^Na_nx_n \right\|^2\leq B \left\|\sum_{n=N+1}^M a_na_n^* \right\|.
$$
This shows us that $\sum_{n=1}^{\infty}a_nx_n$ converges in norm in $X$ for each sequence $(a_k)_k$ from
$\ell^2(\textsf{A})$.

In other words, we have a well-defined map $T:\ell^2(\textsf{A})\rightarrow X$ by the formula $T((a_k)_k)=\sum_{n=1}^{\infty}a_nx_n$. Clearly, $T$ is a module map.

Repeating the above computation we conclude that $T$ is bounded by $B$.
By Proposition 3.6 from \cite{P} $T$ extends to a bounded module map $T: \ell^2_{\text{strong}}(\textsf{A}) \rightarrow X$. Since we obviously have $Te^{(n)}=x_n$ for every $n\in \Bbb N$, Proposition \ref{bessel corr adj} supra implies that $(x_n)_n$ is a weak Bessel sequence.
\qed

\vspace{.1in}

In the light of the preceding proposition it is of interest
to find practical sufficient conditions on an infinite matrix $M=(m_{ij})$ with entries in \textsf{A} which will ensure that the map defined on $\ell^2(\textsf{A})$ by
\begin{equation}\label{schur2}
L_M((a_k)_k)= M\cdot (a_k)_k=\left[\begin{array}{cccc}m_{11} &m_{12}&m_{13}&\ldots\\
m_{21} &m_{22}&m_{23}&\ldots\\
m_{31} &m_{32}&m_{33}&\ldots\\
\vdots&\vdots&\vdots&
\end{array}\right]\left[\begin{array}{c}a_1\\a_2\\a_3\\\vdots\end{array}\right]
\end{equation}
defines a bounded module map with values in $\ell^2_{\text{strong}}(\textsf{A})$ (or, possibly, even in $\ell^2(\textsf{A})$).

Note that here again we need commutativity of the underlying algebra to ensure that the resulting map is a module map.


Here we provide a generalization of the Schur test for infinite matrices.

\vspace{.1in}

Suppose that $M=(m_{ij})$ is an infinite matrix with entries in a commutative von Neumann algebra \textsf{A}. Let
$$
p_{jk}=(m_{jk}^*m_{jk})^{\frac{1}{2}}\geq 0,\,\,\,\forall j,k\in \Bbb N.
$$
Observe that, since \textsf{A} is a von Neumann algebra, $p_{jk}$ and $p_{jk}^{\frac{1}{2}}$ belong to \textsf{A} for all $j$ and $k$. Consider the following condition: there exist positive real constants $B_c$ and $B_r$ such that
\begin{equation}\label{schur3}
\sum_{j=1}^{\infty}\|m_{jk}\|\leq B_c,\,\,\forall k\in \Bbb N \mbox{ and } \sum_{k=1}^{\infty}\|m_{jk}\|\leq B_r,\,\,\forall j\in \Bbb N.
\end{equation}
There are two more conditions related to the preceding one:
\begin{equation}\label{schur4}
\left\|\text{(norm)}\sum_{j=1}^{\infty}p_{jk}\right\|\leq B_c,\,\,\forall k\in \Bbb N \mbox{ and }\left\|\text{(norm)}\sum_{k=1}^{\infty}p_{jk}\right\|\leq B_r,\,\,\forall j\in \Bbb N,
\end{equation}
and
\begin{equation}\label{schur5}
\left\|\text{(strong)}\sum_{j=1}^{\infty}p_{jk}\right\|\leq B_c,\,\,\forall k\in \Bbb N \mbox{ and }\left\|\text{(strong)}\sum_{k=1}^{\infty}p_{jk}\right\|\leq B_r,\,\,\forall j\in \Bbb N.
\end{equation}
Note that we implicitly assume that the series in \eqref{schur4} and \eqref{schur5} converge in indicated topologies in \textsf{A}. Observe now that
$$
\|p_{jk}\|^2=\|p_{jk}^2\|=\|m_{jk}^*m_{jk}\|=\|m_{jk}\|^2,\,\,\,\forall j,k\in \Bbb N.
$$
Thus, we require in \eqref{schur4} that the series $\sum_{j=1}^{\infty}p_{jk}$ and $\sum_{k=1}^{\infty}p_{jk}$ converge in norm in \textsf{A}, while \eqref{schur3} means that the same series converge absolutely in \textsf{A}.
Hence \eqref{schur3} $\Rightarrow$  \eqref{schur4} $\Rightarrow$  \eqref{schur5}.

\vspace{.1in}

In the proofs of the following two propositions we shall repeatedly use the following facts from general theory of $C^*$-algebras. First, if $(a_k)_k$ and $(b_k)_k$ are sequences in a $C^*$-algebra \textsf{A} such that $0\leq a_k\leq b_k$ for all $k$'s and such that $\sum_{k=1}^{\infty}b_k$ converges in \textsf{A}, then $\sum_{k=1}^{\infty}a_k$ also converges in \textsf{A}. (Reason: $\sum_{k=K_1}^{K_2}a_k\leq \sum_{k=K_1}^{K_2}b_k$ and hence $\left\|\sum_{k=K_1}^{K_2}a_k\right\|\leq \left\|\sum_{k=K_1}^{K_2}b_k\right\|$.)

Secondly, if $a,b\in \textsf{A}$ are such that $0\leq a\leq b$ then we have $c^*ac\leq c^*bc$ for all $c\in \textsf{A}$. Finally, for all $a,c\in \textsf{A}$ we have $a^*c^*ca\leq\|c^*c\|\,a^*a$.

\begin{prop}\label{schur prvi}
Let $M=(m_{ij})$ be an infinite matrix with entries in a commutative von Neumann algebra \textsf{A} that satisfies \eqref{schur4} (with $p_{jk}=(m_{jk}^*m_{jk})^{\frac{1}{2}}$ for all $j,k$). Then the operator $L_M$ defined by \eqref{schur2}, i.e.
\begin{equation}\label{schur6}
(L_M(a_k)_k))_j=\sum_{k=1}^{\infty}m_{jk}a_k,\,\,\,j\in \Bbb N,\,\,(a_k)_k\in \ell^2(\textsf{A}),
\end{equation}
defines a bounded module map $L_M:\ell^2(\textsf{A}) \rightarrow \ell^2(\textsf{A})$ that extends to an adjointable map  $L_M:\ell^2_{\text{strong}}(\textsf{A})\rightarrow \ell^2_{\text{strong}}(\textsf{A})$ such that $\|L_M\|\leq(B_rB_c)^{\frac{1}{2}}$.
\end{prop}
\proof
Take any $(a_k)_k\in \ell^2(\textsf{A})$ and write $m_{jk}=u_{jk}p_{jk}$, $a_k=v_kq_k$, $j,k\in \Bbb N$, where all $u_{jk}$ and $v_k$ are partial isometries and $q_k=(a_ka_k^*)^{\frac{1}{2}}\geq 0$. Observe that all these operators belong to \textsf{A} since \textsf{A} is a von Neumann algebra. In addition, let $z_{jk}=u_{jk}v_k$ for all $j,k\in \Bbb N$. Since \textsf{A} is commutative, we can write
\begin{equation}\label{schur7}
m_{jk}a_k=u_{jk}p_{jk}v_kq_k=(z_{jk}p_{jk}^{\frac{1}{2}})(p_{jk}^{\frac{1}{2}}q_k),\,\,\,\forall j,k\in \Bbb N.
\end{equation}
We first claim that
\begin{equation}\label{schur8}
x=(p_{jk}^{\frac{1}{2}}z_{jk}^*)_k\in \ell^2(\textsf{A}).
\end{equation}
To see this, first observe that for each $k$ we have
$$
(p_{jk}^{\frac{1}{2}}z_{jk}^*)(p_{jk}^{\frac{1}{2}}z_{jk}^*)^*=
z_{jk}^*p_{jk}^{\frac{1}{2}}p_{jk}^{\frac{1}{2}}z_{jk}=
p_{jk}^{\frac{1}{2}}z_{jk}^*z_{jk}p_{jk}^{\frac{1}{2}}\leq \|z_{jk}^*z_{jk}\|\,p_{jk}^{\frac{1}{2}}p_{jk}^{\frac{1}{2}}\leq p_{jk}
$$
which implies that the series
$$
\sum_{k=1}^{\infty}(p_{jk}^{\frac{1}{2}}z_{jk}^*)(p_{jk}^{\frac{1}{2}}z_{jk}^*)^*
$$
converges since by the assumption the series $\sum_{k=1}^{\infty}p_{jk}$ converges in norm in \textsf{A} and we have
\begin{equation}\label{schur8a}
\sum_{k=1}^{\infty}(p_{jk}^{\frac{1}{2}}z_{jk}^*)(p_{jk}^{\frac{1}{2}}z_{jk}^*)^*\leq \sum_{k=1}^{\infty}p_{jk}.
\end{equation}
Next we claim that
\begin{equation}\label{schur9}
y=(p_{jk}^{\frac{1}{2}}q_k)_k\in \ell^2(\textsf{A}).
\end{equation}
Indeed, we have
$$
(p_{jk}^{\frac{1}{2}}q_k)(p_{jk}^{\frac{1}{2}}q_k)^*=q_kp_{jk}^{\frac{1}{2}}p_{jk}^{\frac{1}{2}}q_k=p_{jk}^{\frac{1}{2}}a_ka_k^*p_{jk}^{\frac{1}{2}}\leq p_{jk}^{\frac{1}{2}}\langle (a_k)_k,(a_k)_k\rangle p_{jk}^{\frac{1}{2}}\leq \|\langle (a_k)_k,(a_k)_k\rangle \|\,p_{jk}
$$
which shows us that the series
$$
\sum_{k=1}^{\infty}(p_{jk}^{\frac{1}{2}}q_k)(p_{jk}^{\frac{1}{2}}q_k)^*
$$
converges since $\sum_{k=1}^{\infty}p_{jk}$ is norm-convergent
and
\begin{equation}\label{schur9a}
\sum_{k=1}^{\infty}(p_{jk}^{\frac{1}{2}}q_k)(p_{jk}^{\frac{1}{2}}q_k)^*\leq\sum_{k=1}^{\infty}p_{jk}.
\end{equation}
We now apply the Cauchy-Scwarz inequality $\langle x,y\rangle^*\langle x,y\rangle\leq \| \langle x,x\rangle \|\cdot \langle y,y\rangle$ in the Hilbert $C^*$-module $\ell^2(\textsf{A})$ to the sequences $x$ from \eqref{schur8} and
$y$ from \eqref{schur9} to obtain
\begin{eqnarray}\label{schur10}
\left(\sum_{k=1}^{\infty}m_{jk}a_k \right)\left(\sum_{k=1}^{\infty}m_{jk}a_k \right)^*&\stackrel{\eqref{schur7}}{=}&\left(\sum_{k=1}^{\infty}z_{jk}p_{jk}^{\frac{1}{2}}p_{jk}^{\frac{1}{2}}q_k\right)
\left(\sum_{k=1}^{\infty}z_{jk}p_{jk}^{\frac{1}{2}}p_{jk}^{\frac{1}{2}}q_k\right)^*\nonumber\\
 &=&\langle x,y\rangle^*\langle x,y\rangle\nonumber\\
 &\leq&\| \langle x,x\rangle \| \cdot \langle y,y\rangle\nonumber\\
 &\stackrel{\eqref{schur8a}}{\leq}&\left\| \sum_{k=1}^{\infty}p_{jk}\right\|\left( \sum_{k=1}^{\infty}p_{jk}^{\frac{1}{2}}q_k^2 p_{jk}^{\frac{1}{2}}\right)\nonumber\\
 &\stackrel{\eqref{schur4}}{\leq}&B_r\sum_{k=1}^{\infty}p_{jk}^{\frac{1}{2}}q_k^2 p_{jk}^{\frac{1}{2}}.
\end{eqnarray}
Finally, we now have
\begin{eqnarray*}
\sum_{j=1}^{\infty}(L_M((a_k)_k))_j(L_M((a_k)_k))_j^*&\stackrel{\eqref{schur6}}{=}&\sum_{j=1}^{\infty}\left(\sum_{k=1}^{\infty}m_{jk}a_k\right)\left(\sum_{k=1}^{\infty}m_{jk}a_k\right)^*\\
 &\stackrel{\eqref{schur10}}{\leq}&\sum_{j=1}^{\infty}B_r\sum_{k=1}^{\infty}p_{jk}^{\frac{1}{2}}q_k^2 p_{jk}^{\frac{1}{2}}\\
 &=&\sum_{k=1}^{\infty}B_r\left(\sum_{j=1}^{\infty}p_{jk}\right) q_k^2\\
 &\leq&B_r \left\| \sum_{j=1}^{\infty}p_{jk}\right\|\,\sum_{k=1}^{\infty}q_k^2\\
 &\stackrel{\eqref{schur4}}{\leq}&B_rB_c \sum_{k=1}^{\infty}a_ka_k^*\\
 &=&B_rB_c\,\langle (a_k)_k,(a_k)_k\rangle.
\end{eqnarray*}
Obviously, this implies that the series $\sum_{j=1}^{\infty}(L_M((a_k)_k))_j(L_M((a_k)_k))_j^*$ converges which means that $L_M((a_k)_k)$ is well defined and also $\|L_M((a_k)_k)\|^2\leq B_rB_c\|(a_k)_k\|^2$.
By Proposition 3.6 from \cite{P} $L_M$ extends to a bounded module map $L_M: \ell^2_{\text{strong}}(\textsf{A}) \rightarrow \ell^2_{\text{strong}}(\textsf{A})$. Since $\ell^2_{\text{strong}}(\textsf{A})$ is self-dual, $L_M$ is by Proposition 3.4 from \cite{P} adjointable.
\qed

\vspace{.1in}

\begin{remark}\label{primjena Schura u finijoj varijanti}
Proposition \ref{schur prvi} guarantees that the extended   operator $L_M: \ell^2_{\text{strong}}(\textsf{A}) \rightarrow \ell^2_{\text{strong}}(\textsf{A})$ is adjointable. However, one can ask whether the originally defined operator $L_M: \ell^2(\textsf{A}) \rightarrow \ell^2(\textsf{A})$ is adjointable.

This will be the case if our matrix $M$ satisfies \eqref{schur3}. We have already noted that \eqref{schur3} $\Rightarrow$  \eqref{schur4}; thus, Proposition \ref{schur prvi} applies to infinite matrices $M=(m_{ij})$ that satisfy \eqref{schur3}. Moreover, it is easy to check that in this situation the matrix $M^*=(m_{ji}^*)$ also satisfies \eqref{schur3} (with the roles of the contants $B_r$ and $B_c$ interchanged) and hence defines a bounded operator on $\ell^2(\textsf{A})$. By an easy verification one proves that this defines the adjoint operator to the operator induced by the original matrix $M$.

It is a little bit subtler with those matrices $M$ that satisfy \eqref{schur4}. We now must additionally assume that $M=M^*$ in order to conclude that condition \eqref{schur4} is satisfied also for $M^*$. Then, we again have adjointability of the operator under consideration.

Conveniently enough we will most often use Proposition \ref{schur prvi} applied to the Gramian $G$ associated to some sequence $(x_n)_n$ and then, since for all $n,k$ we have $\langle x_n,x_k\rangle^*=\langle x_k,x_n\rangle^*$, our matrix $G$ has the property  $G=G^*$.
\end{remark}

\vspace{.1in}

\begin{prop}\label{schur drugi}
Let $M=(m_{ij})$ be an infinite matrix with entries from a commutative von Neumann algebra \textsf{A} that satisfies \eqref{schur5} (with $p_{jk}=(m_{jk}^*m_{jk})^{\frac{1}{2}}$ for all $j,k$). Then the operator $L_M$ defined by \eqref{schur2}, i.e.
\begin{equation}\label{schur6w}
(L_M(a_k)_k))_j=\sum_{k=1}^{\infty}m_{jk}a_k,\,\,\,j\in \Bbb N,\,\,(a_k)_k\in \ell^2(\textsf{A}),
\end{equation}
defines a bounded module map $L_M:\ell^2(\textsf{A}) \rightarrow \ell^2_{\text{strong}}(\textsf{A})$ that extends to an adjointable map  $L_M:\ell^2_{\text{strong}}(\textsf{A})\rightarrow \ell^2_{\text{strong}}(\textsf{A})$ such that $\|L_M\|\leq(B_rB_c)^{\frac{1}{2}}$.
\end{prop}
\proof
We must show, for each $x=(a_k)_k\in \ell^2(\textsf{A})$, that
\begin{equation}\label{schur7w}
L_x((d_j)_j)=\left\langle (d_j)_j,\left( \sum_{k=1}^{\infty}m_{jk}a_k \right)_j \right\rangle,\,\,\,(d_j)_j\in \ell^2(\textsf{A}),
\end{equation}
defines a bounded module map on $\ell^2(\textsf{A})$.
Let us first take an arbitrary finite sequence $x=(a_1,\ldots,a_K,0,0,\ldots)$. Consider also finite $d=(d_1,\ldots,c_N,0,0,\ldots)$. Then we have
$$
\|L_x(d)\|^2=\left\| \sum_{j=1}^Nd_j\left( \sum_{k=1}^Km_{jk}a_k\right)^*\right\|^2=\left\| \sum_{j=1}^N\sum_{k=1}^Kd_ja_k^*m_{jk}^*\right\|^2.
$$
Writing again $m_{jk}=u_{jk}p_{jk}$ we get
$$
\|L_x(d)\|^2=\left\| \sum_{j=1}^N\sum_{k=1}^Kd_ja_k^*p_{jk}u_{jk}^*\right\|^2=
\left\| \sum_{j=1}^N\sum_{k=1}^K\left( p_{jk}^{\frac{1}{2}}d_j\right)\left(a_k^* p_{jk}^{\frac{1}{2}}u_{jk}^*\right)\right\|^2.
$$
We now regard this last double sum as the inner product in the Hilbert $C^*$-module $(\textsf{A}^K)^N=(\textsf{A}\times \textsf{A} \times \ldots \times \textsf{A})^N$.
(In order to do so, we understand $a_{jk}=a_k$ for all $j$ and also $d_{jk}=d_j$ for all $k$). Then we use the Cauchy-Schwarz inequality in $(\textsf{A}\times \textsf{A} \times \ldots \times \textsf{A})^N$. In this way, continuing the above computation we obtain

\begin{eqnarray*}
\|L_x(d)\|^2&=&\left\| \sum_{j=1}^N\sum_{k=1}^K\left( p_{jk}^{\frac{1}{2}}d_j\right)\left(a_k^* p_{jk}^{\frac{1}{2}}u_{jk}^*\right)\right\|^2\nonumber\\
&\stackrel{\text{C-S}}{\leq}&\left\| \sum_{j=1}^N\sum_{k=1}^K p_{jk}^{\frac{1}{2}}d_jd_j^*p_{jk}^{\frac{1}{2}}\right\|\cdot \left\| \sum_{j=1}^N\sum_{k=1}^Ku_{jk}p_{jk}^{\frac{1}{2}}a_ka_k^* p_{jk}^{\frac{1}{2}}u_{jk}^*\right\|\nonumber\\
&\leq&\left\| \sum_{j=1}^N\sum_{k=1}^K p_{jk}^{\frac{1}{2}}d_j^*d_jp_{jk}^{\frac{1}{2}}\right\|\cdot \left\| \sum_{j=1}^N\sum_{k=1}^K p_{jk}^{\frac{1}{2}}a_k^*a_kp_{jk}^{\frac{1}{2}}\right\|\nonumber\\
&=&\left\|\sum_{j=1}^Nd_j^*\left(\sum_{k=1}^K p_{jk} \right)d_j \right\|\cdot \left\| \sum_{k=1}^Ka_k^*\left(\sum_{j=1}^N p_{jk}\right) a_k \right\|\nonumber\\
&\leq& \left\|\sum_{j=1}^N\left\|\sum_{k=1}^K p_{jk} \right\|d_j^*d_j \right\|\cdot \left\| \sum_{k=1}^K\left\|\sum_{j=1}^N p_{jk}\right\| a_k^*a_k \right\|\nonumber\\
&\stackrel{\eqref{schur5}}{\leq}& \left\|\sum_{j=1}^NB_rd_j^*d_j \right\|\cdot \left\| \sum_{k=1}^KB_ca_k^* a_k \right\|\nonumber\\
&=&B_cB_r\|x\|^2 \|d\|^2.
\end{eqnarray*}
This proves that $L_x$ is bounded on finite sequences; thus, we can extend it to a bounded module map $L_x$ defined on $\ell^2(\textsf{A})$ such that $\|L_x\|\leq (B_cB_r)^{\frac{1}{2}}\|x\|$. Having in mind that $\ell^2_{\text{strong}}(\textsf{A})$ is the dual of $\ell^2(\textsf{A})$, we conclude that
there exists a unique $x^{\prime}=(a_k^{\prime})_k\in \ell^2_{\text{strong}}(\textsf{A})$ such that
\begin{equation}\label{schur8w}
L_x((d_j)_j)=\langle (d_j)_j,x^{\prime}\rangle,\,\,\,\forall (d_j)_j\in \ell^2(\textsf{A}).
\end{equation}

Comparing this to \eqref{schur7w}, we conclude that $\left(\sum_{k=1}^{\infty}m_{jk}a_k\right)_j=a_j^{\prime}$ for every $j\in \Bbb N$ and therefore $\left( \sum_{k=1}^{\infty}m_{jk}a_k\right)_j \in \ell^2_{\text{strong}}(\textsf{A})$. Moreover, we know that $\left\| \left( \sum_{k=1}^{\infty}m_{jk}a_k\right)_j\right\| =\|(a_j^{\prime})_j\|=\|L_x\|\leq (B_cB_r)^{\frac{1}{2}}\|x\|$. This proves that
\eqref{schur6w} defines a bounded operator $L_M$ on the set of finite sequences which extends to  $L_M:\ell^2(\textsf{A})\rightarrow \ell^2_{\text{strong}}(\textsf{A})$ such that $\|L_M\|\leq (B_cB_r)^{\frac{1}{2}}$.

The proof is finished again by a combined application of Propositions 3.4 and 3.6 from \cite{P}.
\qed

\vspace{.3in}

\section{Gabor systems as modular sequences of translates}

\vspace{.1in}

In this section we describe two families of Hilbert $C^*$-modules that are naturally connected with Gabor analysis. Our goal in this section is to establish a bijective correspondence between Gabor frames (resp.~Bessel sequences) and modular weak frames (resp.~Bessel sequences) of translates. As in the preceding section, an important role will be played by the standard Hilbert  \textsf{A}-module $\ell^2(\textsf{A})$ and its dual $\ell^2_{\text{strong}}(\textsf{A})$, where \textsf{A} is a von Neumann algebra.
All what is said holds, {\it mutatis mutandis}, for $\ell^2$-modules indexed by the set of all integers.
To indicate indexation over $\Bbb Z$ we shall write
$\ell^2_{\Bbb Z}(\textsf{A})$ and $\ell^2_{\Bbb Z, \text{strong}}(\textsf{A})$. Since all relevant sums converge unconditionally, we shall operate with symmetric partial sums of the form $\sum_{n=-N}^Na_nb_n^*$,  $\text{(strong)}\sum_{n=-N}^Na_nb_n^*$, etc.

\vspace{.1in}

Also, here and in the rest of the paper, in each norm or inner product under consideration the ambient space will be indicated in the subscript, e.g.~$\| \cdot \|_{L^2(\Bbb R)}$, $\langle \cdot,\cdot\rangle_{\ell^2_{\Bbb Z, \text{strong}}(\textsf{A})}$, etc.

\vspace{.1in}

Let $a>0$. Recall that $L^{\infty}[0,a]$ is a von Neumann algebra with pointwise operations, complex conjugation playing the role of the involution and the norm $\|f\|_{L^{\infty}[0,a]}=\esssup_{x\in[0,a]}|f(x)|$.

Observe that each function $f\in L^{\infty}[0,a]$ naturally extends to a function $f^a\in L^{\infty}(\Bbb R)$, where $f^a:\Bbb R \rightarrow \Bbb C$ denotes the function which is defined by $f^a(x+ka)=f(x)$ for all $x \in [0,a]$ and $k\in \Bbb Z$.

We also note that $L^{\infty}[0,a]$ is naturally represented on the Hilbert space $L^{2}[0,a]$ via $R:L^{\infty}[0,a] \rightarrow \Bbb B(L^{2}[0,a])$, $R(f)=R_f$, $R_f(g)=fg$, $g\in L^{2}[0,a]$.
It is convenient  to note the following observation for future reference.

\vspace{.1in}

\begin{remark}\label{strong vs measure}
Suppose that $(h_N)_N$ is a sequence in $L^{\infty}[0,a]$ such that $0\leq h_1(x)\leq h_2(x)\leq \ldots$ a.e. and $\sup_N\|h_N\|_{L^{\infty}[0,a]}<\infty$. Then, since $L^{\infty}[0,a]$ is a von Neumann algebra, there exists $h=\text{(strong)}\lim_{N\rightarrow \infty}h_N\in L^{\infty}[0,a]$. Moreover, we also have $h(x)=\lim_{N\rightarrow \infty}h_N(x)$ pointwise a.e.

To see this, first observe that, since $[0,a]$ has finite measure, strong convergence coincides with convergence in measure on bounded sets. Thus, here we conclude that $h$ is the limit in measure of the sequence $(h_N)_N$. In general, convergence in measure does not imply convergence pointwise a.e.  But here the assumption guarantees that we have  $g(x):=\lim_{N\rightarrow \infty}h_N(x)$ for a.e. $x$ and, moreover, $\|g\|_{L^{\infty}[0,a]}\leq \sup_N\|h_N\|_{L^{\infty}[0,a]}<\infty$. Hence, $g\in L^{\infty}[0,a]$ and, since convergence pointwise a.e. implies convergence in measure, we conclude that $h=g$.

Conversely, if a bounded sequence $(h_N)_N$ in $L^{\infty}[0,a]$ (that is, $\|h_n\|_{L^{\infty}[0,a]}\leq C$ for every $N$ and some constant $C$) converges pointwise a.e.~to a function $h\in L^{\infty}[0,a]$, then it also converges in measure and, since it is bounded, converges strongly to $h$.

Notice that in general (i.e. without the assumption that the sequence under consideration is essentially monotone) we cannot conclude that strong convergence in $L^{\infty}[0,a]$ implies convergence pointwise a.e. In this light it is useful to note the following observation.

Suppose we have $f,g\in \ell^2_{\Bbb Z, \text{strong}}(L^{\infty}[0,a])$, $f=(f_n)_n$, $g=(g_n)_n$. Then by the definition of the inner product in $\ell^2_{\Bbb Z, \text{strong}}(L^{\infty}[0,a])$ we have
$\langle f,g\rangle_{\ell^2_{\Bbb Z, \text{strong}}(L^{\infty}[0,a])}=\text{(strong)}\sum_{n\in \Bbb Z}f_n\overline{g_n}$. Moreover, we claim that also
$\langle f,g\rangle_{\ell^2_{\Bbb Z, \text{strong}}(L^{\infty}[0,a])}(x)=\sum_{n\in \Bbb Z}f_n\overline{g_n}(x)$ for a.e. $x\in [0,a]$.

This follows from the polarization formula (that holds in every inner product module)
$$
\langle f,g\rangle_{\ell^2_{\Bbb Z, \text{strong}}(L^{\infty}[0,a])}=\sum_{k=0}^3\frac{i^k}{4}\langle f+i^kg,f+i^kg\rangle_{\ell^2_{\Bbb Z, \text{strong}}(L^{\infty}[0,a])}
$$
and the first part of this remark.
We also observe that convergence of the series $\langle f,g\rangle_{\ell^2_{\Bbb Z, \text{strong}}(L^{\infty}[0,a])}$ in the strong operator topology (and hence also pointwise) is  unconditional.
\end{remark}

\vspace{.1in}

For each $a\in \Bbb R$ let $T_a$ denote translation by $a$; that is, the operator given by $T_af(x)=f(x-a)$, where $f$ is any function on $\Bbb R$. Modulation by $b$ is the operator $M_b$ defined by $M_bf(x)=\text{e}^{2\pi ibx}f(x)$. Note that $T_a$ and $M_b$ are unitary operators on $L^2(\Bbb R)$ for all $a,b\in \Bbb R$.

Let us now fix $g\in L^2(\Bbb R)$ and $a,b\in \Bbb R$. We denote by $G(g,a,b)$ the Gabor system generated by $g$ with the lattice parameters $a$ and $b$, i.e.~the sequence $\left(M_{mb}T_{na}g \right)_{m,n\in \Bbb Z}$. For general facts about Gabor systems we refer the reader to \cite{C}, \cite{G} and \cite{H}; see also \cite{G1}.

\vspace{.1in}

Consider now
an arbitrary function $f\in L^2(\Bbb R)$ and $a>0$. For each integer $n$ and for every $x$ we denote by  $(f_n(x))_n$ the sequence defined by $f_n(x)=f(x-na)$. Using the standard periodization trick (see e.g.~\cite{G}, Lemma 1.4.1)
$$
\int_0^a\sum_{n\in \Bbb Z}|f(x-na)|^2dx=\sum_{n\in \Bbb Z}\int_0^a|f(x-na)|^2dx=\int_{-\infty}^{\infty}|f(x)|^2dx=\|f\|_{L^2(\Bbb R)}^2<\infty
$$
we conclude that $\sum_{n\in \Bbb Z}|f(x-na)|^2<\infty$ a.e. which means that
$(f_n(x))_n\in \ell^2_{\Bbb Z}$ for a.e. $x$.
Moreover, by the general $\ell^2$-theory it now follows that the series
\begin{equation}\label{first inner}
\langle f,g\rangle_a(x)=\sum_{n\in \Bbb Z}f(x-na)\overline{g(x-na)}
\end{equation}
converges absolutely for a.e. $x$ and for all $f,g$ from $L^2(\Bbb R)$. The resulting function $\langle f,g\rangle_a$ is $a$-periodic and belongs to $L^{1}[0,a]$. Thus, we have
a map $\langle \cdot,\cdot\rangle_a:L^2(\Bbb R)\times L^2(\Bbb R)\rightarrow L^{1}[0,a]$ and it is known that this map possesses all properties of a vector-valued inner product. For the details we refer the reader to  \cite{CasL}. We also point out that this map under the name the bracket product has been successfully used in the study of shift invariant systems; see \cite{DDR1},\cite{DDR2}, \cite{RS1},\cite{RS2}.

However, $L^{1}[0,a]$ is not a $C^*$-algebra and if want to end up with a Hilbert $C^*$-module we must restrict ourselves to a suitable class of functions (as it has been done in
\cite{CocoL}).

\vspace{.1in}

\begin{definition}\label{bl}
For any $a>0$ let
\begin{equation}\label{bl space}
L_a^{\infty}(\ell^2)=\{f:\Bbb R \rightarrow \Bbb C: f \mbox{ measurable and }\|f\|_{L_a^{\infty}(\ell^2)}^2=\esssup_{x\in [0,a]}\sum_{n\in \Bbb Z}|f(x-na)|^2<\infty\}.
\end{equation}
\end{definition}

First, we note that $C_c(\Bbb R)\subseteq L_a^{\infty}(\ell^2) \subseteq L^2(\Bbb R) \cap L^{\infty}(\Bbb R)$.

Secondly, one easily verifies that  $L_a^{\infty}(\ell^2)$ is
a vector space with pointwise operations. (Some work is needed to check that $L_a^{\infty}(\ell^2)$ is closed under addition, but we omit a verification.) It is also a left $L^{\infty}[0,a]$-module with the action $(h,f) \mapsto h^af$, $h\in L^{\infty}[0,a]$, $f\in L_a^{\infty}(\ell^2)$.

Having in mind \eqref{bl space} it is now natural to define a map
$$
\langle f,g\rangle_{L_a^{\infty}(\ell^2)}:L_a^{\infty}(\ell^2) \times L_a^{\infty}(\ell^2) \rightarrow L^{\infty}[0,a]
$$
by
\begin{equation}\label{bl inner}
\langle f,g\rangle_{L_a^{\infty}(\ell^2)}(x)=\sum_{n\in \Bbb Z}f(x-na)\overline{g(x-na)}=\sum_{n\in \Bbb Z}T_{na}f(x)\overline{T_{na}g(x)}.
\end{equation}
Note that for $g=f$ this reduces to
\begin{equation}\label{bl inner square}
\langle f,f\rangle_{L_a^{\infty}(\ell^2)}(x)=\sum_{n\in \Bbb Z}|f(x-na)|^2
\end{equation}
which is, by \eqref{bl space}, a function from $L^{\infty}[0,a]$.
Using polarization we now conclude that $\langle f,g\rangle_{L_a^{\infty}(\ell^2)}$ introduced in \eqref{bl inner} does belong $L^{\infty}[0,a]$, for all $f,g\in L_a^{\infty}(\ell^2)$.

Note that the defining formulae for $\langle \cdot,\cdot \rangle_{L_a^{\infty}(\ell^2)}$ and $\langle \cdot,\cdot \rangle_{a}$ coincide; by writing $\langle \cdot,\cdot \rangle_{L_a^{\infty}(\ell^2)}$ we just emphasize the fact that both factors belong to  $L_a^{\infty}(\ell^2)$ and hence the result is a function from $L^{\infty}[0,a]$.

As we already observed,  $\langle \cdot,\cdot\rangle_{L_a^{\infty}(\ell^2)}$ has all necessary properties of an $L^{\infty}[0,a]$-valued inner product.
So, in order to conclude that $L_a^{\infty}(\ell^2)$ is a Hilbert $L^{\infty}[0,a]$-module, it only remains to see that $L_a^{\infty}(\ell^2)$ is complete with respect to $\|\cdot\|_{L_a^{\infty}(\ell^2)}$ introduced in \eqref{bl space}. This is also already known; see \cite{CocoL}. However, in Theorem \ref{dual discovery} below we will prove more. But first, three remarks are in order.

\vspace{.1in}

\begin{remark}\label{usporedba normi}
$L_a^{\infty}(\ell^2)$ embeds continuously in $L^2(\Bbb R)$:
\begin{equation}\label{usporedba normi eq}
\|f\|_{L^2(\Bbb R)}\leq \sqrt{a}\,\|f\|_{L_a^{\infty}(\ell^2)},\,\,\,\forall f \in L_a^{\infty}(\ell^2).
\end{equation}
Indeed,
$$
\|f\|_{L^(\Bbb R)}^2=\int_{-\infty}^{\infty}|f(x)|^2dx=\sum_{n\in \Bbb Z}\int_0^a|f(x-na)|^2dx=\int_0^a\sum_{n\in \Bbb Z}|f(x-na)|^2dx\leq a\|f\|_{L_a^{\infty}(\ell^2)}^2<\infty.
$$
\end{remark}

\vspace{.1in}

\begin{remark}\label{taman}
Suppose that $G(g,a,b)$ is a Bessel sequence with a Bessel bound $B$. Then it is well known (see Theorem 11.6 in \cite{H} or Proposition 8.3.2 in \cite{C}) that
$\sum_{n\in \Bbb Z}|g(x-na)|^2<Bb$ a.e. Thus, $g\in L_a^{\infty}(\ell^2)$. Moreover, by Lemma 9.2.2 in \cite{C} $G(g,a,b)$ is a Bessel sequence with a Bessel bound $B$ if and only if $G(g, \frac{1}{b},\frac{1}{a})$ is a Bessel sequence with a Bessel bound $Bab$. Hence, we also have $g\in L_{\frac{1}{b}}^{\infty}(\ell^2)$.
This shows that all generators $g$ of Gabor Bessel sequences $G(g,a,b)$ are contained in $L_a^{\infty}(\ell^2) \cap L_{\frac{1}{b}}^{\infty}(\ell^2)$.
\end{remark}

\vspace{.1in}

\begin{remark}\label{4 prostora}
In fact, $L_a^{\infty}(\ell^2)$ is the largest of four spaces that naturally appear in this context. For $a>0$ consider the following spaces of measurable functions $f:\Bbb R\rightarrow \Bbb C$:
\begin{equation}\label{wa1 space}
W_a(L^{\infty},\ell^1)=\left\{f: \sum_{n\in \Bbb Z}\left(\esssup_{x\in [0,a]}|f(x-na)|\right)<\infty\right\},
\end{equation}
\begin{equation}\label{wa2 space}
W_a(L^{\infty},\ell^2)=\left\{f: \sum_{n\in \Bbb Z}\left(\esssup_{x\in [0,a]}|f(x-na)|\right)^2<\infty\right\},
\end{equation}
\begin{equation}\label{Wienerlike am space}
X_a(L^{\infty},\ell^2)=\left\{f: \sum_{n\in \Bbb Z}\left|T_{-na}f\chi_{[0,a]}\right|^2 \mbox{ converges in norm in } L^{\infty}[0,a]\right\}.
\end{equation}
The first two are well known Wiener amalgam spaces. We claim that $W_a(L^{\infty},\ell^1)\subseteq W_a(L^{\infty},\ell^2)\subseteq X_a(L^{\infty},\ell^2)\subseteq L_a^{\infty}(\ell^2)$.

It is clear that $W_a(L^{\infty},\ell^1)\subseteq W_a(L^{\infty},\ell^2)$.
To prove that $W_a(L^{\infty},\ell^2)\subseteq X_a(L^{\infty},\ell^2)$ take any $f\in W_a(L^{\infty},\ell^2)$ and observe that
\begin{equation}\label{wa var space}
\left(\esssup_{x\in [0,a]}|f(x-na)|\right)^2=\esssup_{x\in [0,a]}|f(x-na)|^2,\,\,\,\forall n\in \Bbb Z.
\end{equation}
Now we have
$$
\sum_{n\in \Bbb Z}\left\| \left|T_{-n}f\chi_{[0,a]}\right|^2 \right\|_{L^{\infty}[0,a]}=\sum_{n\in \Bbb Z}\esssup_{x\in [0,a]}\left| T_{-n}f(x)\right|^2
 =\sum_{n\in \Bbb Z}\esssup_{x\in [0,a]}\left|f(x+na) \right|^2\stackrel{\eqref{wa2 space}, \eqref{wa var space}}{<}\infty.
$$
This tells us that  $\sum_{n\in \Bbb Z}\left|T_{-na}f\chi_{[0,a]}\right|^2$ is an absolutely convergent series in $L^{\infty}[0,a]$. Since $L^{\infty}[0,a]$ is a Banach space we conclude that $\sum_{n\in \Bbb Z}\left|T_{-na}f\chi_{[0,a]}\right|^2$ converges in norm in in $L^{\infty}[0,a]$; hence, $f \in X_a(L^{\infty},\ell^2)$.

Finally, for each $f \in X_a(L^{\infty},\ell^2)$ the series $\sum_{n\in \Bbb Z}\left|T_{-na}f\chi_{[0,a]}\right|^2$ converges in norm in $L^{\infty}[0,a]$ to a function, call it $h$, from $L^{\infty}[0,a]$. This implies pointwise a.e.~convergence to $h$. But when the series $\sum_{n\in \Bbb Z}\left|T_{-na}f\chi_{[0,a]}\right|^2$ converges to a function from $L^{\infty}[0,a]$ pointwise a.e.~then, clearly, $f \in L_a^{\infty}(\ell^2)$.

The functions $$f_1=\sum_{n=0}^{\infty}\frac{1}{n+1}\chi_{[n+1-\frac{1}{2^n},n+1-\frac{1}{2^{n+1}}]}\in W_1(L^{\infty},\ell^2)\setminus W_1(L^{\infty},\ell^1),$$ $$f_2=\sum_{n=0}^{\infty}\frac{1}{\sqrt{n+1}}\chi_{[n+1-\frac{1}{2^n},n+1-\frac{1}{2^{n+1}}]}\in X_1(L^{\infty},\ell^2) \setminus W_1(L^{\infty},\ell^2),$$ and
$$f_3=\sum_{n=0}^{\infty}\chi_{[n+1-\frac{1}{2^n},n+1-\frac{1}{2^{n+1}}]}\in L_1^{\infty}(\ell^2) \setminus X_1(L^{\infty},\ell^2)$$ show us that the above inclusions are strict.
\end{remark}

\vspace{.1in}

Observe that, for each $f\in X_a(L^{\infty},\ell^2)$, the sequence $\left( T_{-na}f\chi_{[0,a]}\right)_n$ belongs to $\ell^2_{\Bbb Z}(L^{\infty}[0,a])$. Conversely, given a sequence
$(f_n)_n$ in $\ell^2_{\Bbb Z}(L^{\infty}[0,a])$, we can define a function $f\in L^{\infty}(\Bbb R)$ by $f|_{[na,(n+1)a]}=T_{na}f_n$, for all $n\in \Bbb Z$.
Note that $f$ is in fact defined by $f(x+na)=f_n(x)$, for all $x\in [0,a]$ and $n\in \Bbb Z$.
Then, clearly, $f\in
X_a(L^{\infty},\ell^2)$. This shows that $X_a(L^{\infty},\ell^2)$ is  a copy of $\ell^2_{\Bbb Z}(L^{\infty}[0,a])$ and hence can be endowed with the structure of a left Hilbert $L^{\infty}[0,a]$-module. This brings us to the following theorem which shows the nature of $X_a(L^{\infty},\ell^2)$ and $L_a^{\infty}(\ell^2)$ as Hilbert $C^*$-modules.

\vspace{.1in}

\begin{theorem}\label{dual discovery}
$L_a^{\infty}(\ell^2)$, $a>0$,
is a Hilbert $C^*$-module over $L^{\infty}[0,a]$ with the pointwise operations, the action of $L^{\infty}[0,a]$ defined by
$(h,f) \mapsto h^af$, $f\in L_a^{\infty}(\ell^2)$, $h\in L^{\infty}[0,a]$, and the inner product defined by \eqref{bl inner}:
$$
\langle f,g\rangle_{L_a^{\infty}(\ell^2)}(x)=\sum_{n\in \Bbb Z}f(x-na)\overline{g(x-na)},
$$
that contains $X_a(L^{\infty},\ell^2)$ as a closed Hilbert $L^{\infty}[0,a]$-submodule.
The map
\begin{equation}\label{analysis}
U_a:L_a^{\infty}(\ell^2) \rightarrow \ell^2_{\Bbb Z, \text{strong}}(L^{\infty}[0,a]),\,\,\,U_af=(T_{-na}f\chi_{[0,a]})_n,
\end{equation}
is a unitary operator of Hilbert $L^{\infty}[0,a]$-modules whose adjoint $U^*$ is given by $U^*((f_n)_n)=f$, where $f$ is defined by $f|_{[na,(n+1)a]}=T_{na}f_n$ for all $n$ in $\Bbb Z$. Moreover, $$U_a|_{X_a(L^{\infty},\ell^2)}: X_a(L^{\infty},\ell^2)\rightarrow \ell^2_{\Bbb Z}(L^{\infty}[0,a])$$ is also a unitary operator of Hilbert $C^*$-modules.

In particular, $L_a^{\infty}(\ell^2)$ is self-dual.
Finally, $L_a^{\infty}(\ell^2)$ is naturally (anti-linerly) isomorphic to the dual $X_a(L^{\infty},\ell^2)^{\prime}$ of $X_a(L^{\infty},\ell^2)$.
\end{theorem}
\proof
First observe that
$$
\esssup_{x\in [0,a]}\sum_{n\in \Bbb Z}|f(x-na)|^2<\infty
$$
if and only if
$$
\sup_{N\in \Bbb N}
\esssup_{x\in [0,a]}\sum_{n=-N}^N|f(x-na)|^2<\infty
$$
in which case these two expressions coincide. (We omit a verification.)
Therefore, we can rewrite \eqref{bl space} as
$$
L_a^{\infty}(\ell^2)=\left\{f: \|f\|_{L_a^{\infty}(\ell^2)}^2=\sup_{N\in \Bbb N}\left\|\sum_{n=-N}^N\left|T_{-na}f\chi_{[0,a]}\right|^2\right\|_{L^{\infty}[0,a]}<\infty\right\}.
$$
This means that for each $f\in L_a^{\infty}(\ell^2)$ we have $U_af=(T_{-na}f\chi_{[0,a]} )_n\in \ell^2_{\Bbb Z, \text{strong}}(L^{\infty}[0,a])$.

For the same reasons, the map $V: \ell^2_{\Bbb Z, \text{strong}}(L^{\infty}[0,a]) \rightarrow L_a^{\infty}(\ell^2)$, $V((f_n)_n)=f$, where $f$ is given by $f|_{[na,(n+1)a]}=T_{na}f_n$ for all $n$ in $\Bbb Z$ is also well defined. By a routine verification which we also omit one shows that both $U_a$ and $V$ are module maps and that $V=U_a^{-1}$.

Recall that for all $(f_n)_n\in \ell^2_{\Bbb Z, \text{strong}}(L^{\infty}[0,a])$ we have
\begin{equation}\label{cvg in measure 1}
\langle (f_n)_n,(f_n)_n \rangle_{\ell^2_{\Bbb Z, \text{strong}}(L^{\infty}[0,a])} = \text{(strong)}\sum_{n\in \Bbb Z}|f_n|^2,
\end{equation}
where the series converges in the strong operator topology on $L^{\infty}[0,a]$. Using Remark \ref{strong vs measure} we see that the series that defines the inner product in $\ell^2_{\Bbb Z, \text{strong}}(L^{\infty}[0,a])$ that is, $\langle (f_n)_n,(g_n)_n \rangle_{\ell^2_{\Bbb Z, \text{strong}}(L^{\infty}[0,a])}=\text{(strong)}\sum_{n\in \Bbb Z}f_n\overline{g_n}$ also converges pointwise a.e.

On the other hand, we have
\begin{equation}\label{cvg in measure 2}
\langle U_a^{-1}(f_n)_n,U_a^{-1}(f_n)_n \rangle_{L_a^{\infty}(\ell^2)}(x) = \sum_{n\in \Bbb Z}|f(x-na)|^2=\sum_{n\in \Bbb Z}|f_n(x)|^2,
\end{equation}
where this series converges pointwise a.e. Hence, sums in \eqref{cvg in measure 1} and \eqref{cvg in measure 2} coincide.
Thus, we have obtained
\begin{equation}\label{izometrija}
\langle U_a^{-1}(f_n)_n,U_a^{-1}(f_n)_n \rangle_{L_a^{\infty}(\ell^2)}=\langle (f_n)_n,(f_n)_n \rangle_{\ell^2_{\Bbb Z, \text{strong}}(L^{\infty}[0,a])},\,\,\,\forall (f_n)_n\in \ell^2_{\Bbb Z, \text{strong}}(L^{\infty}[0,a]).
\end{equation}
This shows that $U_a^{-1}$ is an isometry. In particular, since $\ell^2_{\Bbb Z, \text{strong}}(L^{\infty}[0,a])$ is complete, we conclude that
$L_a^{\infty}(\ell^2)$ is also complete.

Furthermore, \eqref{izometrija} by polarization implies
$$
\langle U_a^{-1}(f_n)_n,U_a^{-1}(g_n)_n \rangle_{L_a^{\infty}(\ell^2)}=\langle (f_n)_n,(g_n)_n \rangle_{\ell^2_{\Bbb Z, \text{strong}}(L^{\infty}[0,a])},\,\,\,\forall (f_n)_n, (g_n)_n\in \ell^2_{\Bbb Z, \text{strong}}(L^{\infty}[0,a]).
$$
If we now take $U_a^{-1}(f_n)_n=f$ this gives us
$$
\langle f,U_a^{-1}(g_n)_n \rangle_{L_a^{\infty}(\ell^2)}=\langle U_af,(g_n)_n \rangle_{\ell^2_{\Bbb Z, \text{strong}}(L^{\infty}[0,a])},\,\,\,\forall f\in L_a^{\infty}(\ell^2),\,\forall (g_n)_n\in \ell^2_{\Bbb Z, \text{strong}}(L^{\infty}[0,a]);
$$
thus, $U_a^{-1}=U_a^*$.

We now know that $L_a^{\infty}(\ell^2)$ and $\ell^2_{\Bbb Z, \text{strong}}(L^{\infty}[0,a])$ are unitary equivalent Hilbert $C^*$-modules. Since $\ell^2_{\Bbb Z, \text{strong}}(L^{\infty}[0,a])$ is self-dual, $L_a^{\infty}(\ell^2)$ has the same property.

At the same time we also see that the restriction of $U_a$ to $X_a(L^{\infty},\ell^2)$ is a bijection between the submodules $X_a(L^{\infty},\ell^2)$ and $\ell^2_{\Bbb Z}(L^{\infty}[0,a])$. In particular, $X_a(L^{\infty},\ell^2)$ is complete and hence a Hilbert $C^*$-module over $L^{\infty}[0,a]$. Finally, using unitary equivalence of $X_a(L^{\infty},\ell^2)$ and $L_a^{\infty}(\ell^2)$ with $\ell^2_{\Bbb Z}(L^{\infty}[0,a])$ and $\ell^2_{\Bbb Z, \text{strong}}(L^{\infty}[0,a])$, respectively, we conclude that
the dual of $X_a(L^{\infty},\ell^2)$ is $L_a^{\infty}(\ell^2)$.
\qed

\vspace{.2in}

\begin{remark}\label{norms}
Note that the norm on $X_a(L^{\infty},\ell^2)$ arising from its inner product is given by
\begin{equation}\label{prvi modul norm}
\|f\|_{X_a(L^{\infty},\ell^2)} =\left\|\langle f,f\rangle_{X_a(L^{\infty},\ell^2)}\right\|_{L^{\infty}[0,a]}^{\frac{1}{2}}=\left\|\sum_{n\in \Bbb Z} \left|T_{-na}f\chi_{[0,a]}\right)|^2\right\|_{L^{\infty}[0,a]}^{\frac{1}{2}}.
\end{equation}
which can be written as
$$
\|f\|_{X_a(L^{\infty},\ell^2)}=\left(\esssup_{x\in [0,a]}\sum_{n\in \Bbb Z}|f(x-na)|^2 \right)^{\frac{1}{2}}.
$$
This, of course, agrees with \eqref{bl space}, but the difference is that here, for functions in $X_a(L^{\infty},\ell^2)$, we have norm-convergence of the series, while the series in \eqref{bl space} converges only pointwise a.e. Precisely the same is true for inner products.

Observe that this reflects the analogous situation in $\ell^2$-modules. The inner product on $\ell^2_{\Bbb Z, \text{strong}}(\textsf{A})$ is given by
$\langle (f_n)_n,(g_n)_n\rangle_{\ell^2_{\Bbb Z, \text{strong}}(\textsf{A})}=\text{(strong)}\sum_{n\in \Bbb Z}f_ng_n^*$ but this series converges actually in norm in the underlying algebra \textsf{A} for all $(f_n)_n$ and $(g_n)_n$ from $\ell^2_{\Bbb Z}(\textsf{A})$.
\end{remark}

\vspace{.1in}

\begin{remark}\label{Wienerove norme}
Recall that $W_a(L^{\infty},\ell^1)\subseteq W_a(L^{\infty},\ell^2) \subseteq X_a(L^{\infty},\ell^2)\subseteq L_a^{\infty}(\ell^2)$ and that
\begin{equation}\label{was}
\|f\|_{W_a(L^{\infty},\ell^2)}=\left(\sum_{k\in \Bbb Z}\esssup_{x\in [0,a]}|f(x-ka)|^2 \right)^{\frac{1}{2}}
\end{equation}
is a norm on $W_a(L^{\infty},\ell^2)$. It is easy to see that $\|f\|_{L_a^{\infty}(\ell^2)}\leq \|f\|_{W_a(L^{\infty},\ell^2)}$. One can also check that
$\|f\|_{L_a(L^{\infty},\ell^1)}\leq \|f\|_{W_a(L^{\infty},\ell^2)}$ for all $f$ from $W_a(L^{\infty},\ell^1)$.

Note also that $W_a(L^{\infty},\ell^1)$ is dense in $X_a(L^{\infty},\ell^2)$ with respect to  $\|\cdot\|_{L_a^{\infty}(\ell^2)}$. In fact, its subspace consisting of all essentially bounded functions with compact support is already dense in $X_a(L^{\infty},\ell^2)$ with respect to  $\|\cdot\|_{L_a^{\infty}(\ell^2)}$. Observe that the set of all essentially bounded functions with compact support is via the unitary operator $U_a$ from Theorem \ref{dual discovery} in a bijective correspondence with the set of all finite sequences in $\ell^2_{\Bbb Z}(L^{\infty}[0,a])\subseteq \ell^2_{\Bbb Z, \text{strong}}(L^{\infty}[0,a])$.
\end{remark}

\vspace{.2in}

We are now ready for the main result in this section in which we establish a correspondence of Gabor Bessel sequences (weak frames) and weak Bessel sequences (frames) of translates in our Hilbert $C^*$-module $L_a^{\infty}(\ell^2)$. In fact, it turns out that the translation paremeter $a$ continues to play the same role, while the modulation parameter $b$ determines the ambient module.

\vspace{.1in}

\begin{theorem}\label{CasCocLem}
Let $g\in L^2(\Bbb R)$ and $a,b>0$. Then $G(g,a,b)$ is a Bessel sequence in $L^2(\Bbb R)$ with a Bessel bound $B$ (a frame with frame bounds $A$ and $B$) if and only if the sequence $(T_{na}(\frac{1}{\sqrt{b}}g))_{n\in \Bbb Z}$ is a weak Bessel sequence with a Bessel bound $B$ (a weak frame with frame bounds $A$ and $B$) in the Hilbert $L^{\infty}[0,\frac{1}{b}]$-module $L_{\frac{1}{b}}^{\infty}(\ell^2)$.
\end{theorem}
\proof
Recall that for all $f,g\in L^2(\Bbb R)$ and $a>0$ we have a function $\langle f,g\rangle_a\in L^{1}[0,a]$ defined by \eqref{first inner}.

Let us take arbitrary $g\in L^2(\Bbb R)$ and $a,b>0$.
By Corollary 4.6.17 from \cite{CasL} $G(g,a,b)$ is a frame with frame bounds $A$ and $B$ if and only if
\begin{equation}\label{casazza lammers}
A\langle f,f\rangle_{\frac{1}{b}}(x)\leq \sum_{n\in \Bbb Z}\left|\left\langle f,T_{na}(\frac{1}{\sqrt{b}}g)\right\rangle_{\frac{1}{b}}(x)\right|^2\leq B\langle f,f\rangle_{\frac{1}{b}}(x),\mbox{ for a.e. }x,\,\forall f\in L^2(\Bbb R).
\end{equation}

By Theorem 5.5 from \cite{CocoL}, \eqref{casazza lammers} is equivalent to
\begin{equation}\label{coco lammers}
A\langle f,f\rangle_{L_{\frac{1}{b}}^{\infty}(\ell^2)}(x)\leq \sum_{n\in\Bbb Z}\left|\left\langle f,T_{na}(\frac{1}{\sqrt{b}}g\right\rangle_{L_{\frac{1}{b}}^{\infty}(\ell^2)}(x) \right|^2\leq B\langle f,f\rangle_{L_{\frac{1}{b}}^{\infty}(\ell^2)}(x),\mbox{for a.e.}x,
\forall f\in L_{\frac{1}{b}}^{\infty}(\ell^2).
\end{equation}
It should be mentioned, however, that the factor $\frac{1}{\sqrt{b}}$ that multiplies $g$ in \eqref{casazza lammers} and \eqref{coco lammers} is missing in both Corollary 4.6.17 from \cite{CasL} and Theorem 5.5 from \cite{CocoL}\footnote{This really simple mistake  - a missing $\sqrt{b}$ in the denominator - origins actually from Theorem 4.6.3 in \cite{CasL}. Unfortunately, the key result (Corollary 4.6.17) from \cite{CasL} is quoted and used in \cite{CocoL} in this wrong form without the factor $\frac{1}{\sqrt{b}}$.}.

Having obtained \eqref{coco lammers} we are just one step from the end of the proof. First, in the sum over $n\in \Bbb Z$ in \eqref{coco lammers} we have convergence pointwise a.e. Observe that that sum is, by the righthand side inequality, an essentially bounded function. Since pointwise a.e. convergence in $L^{\infty}[0,\frac{1}{b}]$ implies convergence in measure which is in turn equivalent to convergence in the strong operator topology in $L^{\infty}[0,\frac{1}{b}]$, we can rewrite \eqref{coco lammers} in the form
\begin{equation}\label{db}
A\langle f,f\rangle_{L_{\frac{1}{b}}^{\infty}(\ell^2)} \leq \text{(strong)}\sum_{n\in\Bbb Z}\left|\left\langle f,T_{na}(\frac{1}{\sqrt{b}}g\right\rangle_{L_{\frac{1}{b}}^{\infty}(\ell^2)} \right|^2\leq B\langle f,f\rangle_{L_{\frac{1}{b}}^{\infty}(\ell^2)},\,\forall f \in L_{\frac{1}{b}}^{\infty}(\ell^2).
\end{equation}
Thus, by Definition \ref{very first fd} the sequence $(T_{na}(\frac{1}{\sqrt{b}}g))_{n\in \Bbb Z}$ is a weak frame in $L_{\frac{1}{b}}^{\infty}(\ell^2)$ with frame bounds $A$ and $B$.
Conversely, if we have \eqref{db},  we conclude using Remark \ref{strong vs measure} that the series over $n\in \Bbb Z$ in \eqref{db} converges also pointwise a.e. and this gives us
\eqref{coco lammers}.
\qed

\vspace{.4in}

\section{Hilbert $C^*$-modules in Gabor analysis}

\vspace{.2in}

We open the section with demonstrating a couple of useful properties of the Hilbert $C^*$-module $L_a^{\infty}(\ell^2)$, $a>0$. For any $r>0$ we denote by $D_r$ the dilation operator defined by $D_rf(x)=f(rx)$.
Observe that here we work with the dilation operator $D_r$ without the normalizing factor that is usually used when the dilation is regarded as an operator on $L^2(\Bbb R)$.

\vspace{.1in}

\begin{prop}\label{ekvivalentni moduli}
Let $a,c>0$. Then $D_{\frac{a}{c}}:L^{\infty}[0,a] \rightarrow L^{\infty}[0,c]$, $D_{\frac{a}{c}}f(x)=f(\frac{a}{c}x)$, is an isomorphism of  von Neumann algebras and also $D_{\frac{a}{c}}:L_a^{\infty}(\ell^2) \rightarrow L_c^{\infty}(\ell^2)$ is a unitary operator of Hilbert $C^*$-modules. Finally, the restriction $D_{\frac{a}{c}}|_{X_a(L^{\infty},\ell^2)}$ is also a unitary operator of Hilbert $C^*$-modules $X_a(L^{\infty},\ell^2)$ and $X_c(L^{\infty},\ell^2)$.
\end{prop}
\proof
One checks that $D_{\frac{a}{c}}:L^{\infty}[0,a] \rightarrow L^{\infty}[0,c]$, $D_{\frac{a}{c}}f(x)=f(\frac{a}{c}x)$ is an isomorphism of  von Neumann algebras by an easy verification.

Now observe that the map $\varphi : \ell^2_{\Bbb Z, \text{strong}}(L^{\infty}[0,a]) \rightarrow \ell^2_{\Bbb Z, \text{strong}}(L^{\infty}[0,c])$ defined by $\varphi((f_n)_n)=(D_{\frac{a}{c}}f_n)_n$ is a unitary operator (actually a $D_{\frac{a}{c}}$-morphism) of Hilbert $C^*$-modules. This is also seen by a routine verification which we omit; however, let us only mention that here we use again (as in the proof of Proposition \ref{slika slabog framea}) the fact that $D_{\frac{a}{c}}$ is, being an isomorphism of von Neumann algebras, a normal map.

Applying Theorem \ref{dual discovery} we now conclude that the map $\Phi=U_c^*\varphi U_a :L_a^{\infty}(\ell^2) \rightarrow L_c^{\infty}(\ell^2)$ is a unitary operator. We claim that $\Phi$ acts as dilation by $\frac{a}{c}$.
First recall that $U_af=(f_n)_n$ where $f_n=T_{-na}f\chi_{[0,a]}$, i.e.~$f_n(x)=f(x+na)$. Then $\varphi U_af=(D_{\frac{a}{c}}f_n)_n$, where $(D_{\frac{a}{c}}f_n)(x)=f_n(\frac{a}{c}x)=f(\frac{a}{c}x+na)$, $x\in [0,c]$. Finally,
$\Phi(f)=U_c^*\varphi U_a(f)$ means that $\Phi(f)(x+nc)=f(\frac{a}{c}x+na)=D_{\frac{a}{c}}f(x+nc)$, for all $x\in [0,c]$ and $n\in \Bbb Z$; thus, $\Phi=D_{\frac{a}{c}}$.

The last assertion of the proposition follows from the fact that the restriction $\varphi|_{\ell^2_{\Bbb Z}(L^{\infty}[0,a])}$ is a unitary operator of Hilbert $C^*$-modules $\ell^2_{\Bbb Z}(L^{\infty}[0,a])$ and
$\ell^2_{\Bbb Z}(L^{\infty}[0,c])$ together with the corresponding statement of Theorem \ref{dual discovery}.
\qed

\vspace{.1in}

\begin{prop}\label{racionalno isti}
Let $a,c>0$. Then $W_a(L^{\infty},\ell^1)=W_c(L^{\infty},\ell^1)$ and $W_a(L^{\infty},\ell^2)=W_c(L^{\infty},\ell^2)$.  If $\frac{a}{c}\in \Bbb Q$ we also have
$X_a(L^{\infty},\ell^2)=X_c(L^{\infty},\ell^2)$ and $L_a^{\infty}(\ell^2)=L_c^{\infty}(\ell^2)$ and the corresponding norms on $L_a^{\infty}(\ell^2)$ and $L_c^{\infty}(\ell^2)$ are equivalent. Conversely, each of last two equalities implies $\frac{a}{c}\in \Bbb Q$.
\end{prop}
\proof
The first two equalities for Wiener amalgam spaces are well known (see e.g.~\cite{H}, Section 11.4). This allows us to write $W(L^{\infty},\ell^1)$ and $W(L^{\infty},\ell^2)$ without specifying any particular parameter $a$ and we will adopt this convention in the rest of the paper.

Suppose now that $\frac{a}{c}=\frac{p}{q}$ with $p,q\in \Bbb N$, $(p,q)=1$.

For $f\in L_a^{\infty}(\ell^2)$ put $\esssup_{x\in [0,a]}\sum_{l\in \Bbb Z}|f(x-la)|^2=C$.
Take any $y\in [0,c]$. Then the numbers $y,c+y, 2c+y,\ldots,(p-1)c+y$ are contained in the interval $[0,pc]=[0,qa]$ and therefore they are of the form $x_j+k_ja$, $j=1,\ldots,p$ for some
$x_1,\ldots,x_p\in [0,a]$ and $k_1,\ldots,k_p\in \{0,1,\ldots, q-1\}$. Hence
$$
\sum_{k=0}^{p-1}|f(kc+y)|^2=\sum_{j=1}^p|f(x_j+k_ja)|^2.
$$
Now we consider next cycle. We have $y+pc=x_1+k_1a+qa$, $y+(p+1)c=x_2+k_2a+qa$, ..., $y+(2p-1)c=x_p+k_pa+qa$ which together with the preceding equality gives us
$$
\sum_{k=0}^{2p-1}|f(kc+y)|^2=\sum_{j=1}^p\left(|f(x_j+k_ja)|^2+|f(x_j+k_ja+qa)|^2\right).
$$
It now follows that for a.e. $y$ we have
$$
\sum_{k\in \Bbb Z}|f(kc+y)|^2\leq \sum_{l\in \Bbb Z}|f(x_1+la)|^2+\sum_{l\in \Bbb Z}|f(x_2+la)|^2+\ldots +\sum_{l\in \Bbb Z}|f(x_p+la)|^2\leq pC;
$$
thus,
$$
\|f\|_{L_c^{\infty}(\ell^2)}\leq \sqrt{p}\|f\|_{L_a^{\infty}(\ell^2)}.
$$
By symmetry we also conclude
$$
\|f\|_{L_a^{\infty}(\ell^2)}\leq \sqrt{q}\|f\|_{L_c^{\infty}(\ell^2)}.
$$

To prove that $\frac{a}{c}\in \Bbb Q$ is also a necessary condition for the equality $L_a^{\infty}(\ell^2)=L_c^{\infty}(\ell^2)$ we first claim that
\begin{equation}\label{samo racionalno}
L_a^{\infty}(\ell^2)=L_c^{\infty}(\ell^2) \Rightarrow L_{ar}^{\infty}(\ell^2)=L_{cr}^{\infty}(\ell^2),\,\,\,\forall r>0.
\end{equation}
To see this, take any $r>0$ and the corresponding dilation  $D_rf(x)=f(rx)$ and observe that
$$
\sum_{n\in \Bbb Z}|f(x+ncr)|^2=\sum_{n\in \Bbb Z}|f(r(\frac{x}{r}+nc))|^2=\sum_{n\in \Bbb Z}|D_rf(\frac{x}{r}+nc))|^2.
$$
This shows that $f\in L_{cr}^{\infty}(\ell^2)$ if and only if $D_rf\in L_{c}^{\infty}(\ell^2)$ or, equivalently, $D_r^{-1}L_{c}^{\infty}(\ell^2)=L_{cr}^{\infty}(\ell^2)$.
For the same reason we also have $D_r^{-1}L_{a}^{\infty}(\ell^2)=L_{ar}^{\infty}(\ell^2)$ and this two equalities prove
\eqref{samo racionalno}.

Suppose now that that we have $L_a^{\infty}(\ell^2)=L_c^{\infty}(\ell^2)$ for $a,c$ such that $\frac{a}{c}\not \in \Bbb Q$. Then  \eqref{samo racionalno} implies $L_1^{\infty}(\ell^2)=L_{\frac{c}{a}}^{\infty}(\ell^2)$. But this is impossible as demonstrated by an example from \cite{FJ} (see the proof of Proposition 9.6.2 in \cite{C}).

We now turn to submodules  $X_a(L^{\infty},\ell^2)$ and $X_c(L^{\infty},\ell^2)$. Suppose first $\frac{a}{c}\in \Bbb Q$ and take any $f\in X_a(L^{\infty},\ell^2)$. Then by Remark \ref{Wienerove norme} there is a sequence $(f_n)_n$ of bounded functions with compact support such that $\|f-f_n\|_{L_a^{\infty}(\ell^2)}\rightarrow 0$. Since $\| \cdot \|_{L_a^{\infty}(\ell^2)}$ and $\| \cdot \|_{L_c^{\infty}(\ell^2)}$ are equivalent norms on $L_a^{\infty}(\ell^2)=L_c^{\infty}(\ell^2)$, we also have $\|f-f_n\|_{L_c^{\infty}(\ell^2)}\rightarrow 0$. By Remark \ref{Wienerove norme} this implies $f\in X_c(L^{\infty},\ell^2)$.

To end the proof we need to show that $X_a(L^{\infty},\ell^2)=X_c(L^{\infty},\ell^2)$ implies $\frac{a}{c}\in \Bbb Q$. Suppose that $X_a(L^{\infty},\ell^2)=X_c(L^{\infty},\ell^2)$. Since $X_a(L^{\infty},\ell^2)$ and $X_c(L^{\infty},\ell^2)$ are by Proposition \ref{ekvivalentni moduli} unitary equivalent modules, we have two equivalent norms - $\|\cdot \|_{L_a^{\infty}(\ell^2)}$ and $\|\cdot \|_{L_c^{\infty}(\ell^2)}$ on the same set
$X_a(L^{\infty},\ell^2)=X_c(L^{\infty},\ell^2)$. Consider, for any $h\in L_a^{\infty}(\ell^2)$, the map $l_h$ on  $X_a(L^{\infty},\ell^2)$ defined by $l_h(f)=\langle f,h\rangle_{L_a^{\infty}(\ell^2)}$. This map is also a bounded module map on $X_c(L^{\infty},\ell^2)$ and therefore there exists $\tilde{h}\in L_c^{\infty}(\ell^2)$ such that $l_h(f)=\langle f,\tilde{h}\rangle_{L_c^{\infty}(\ell^2)}$, for all $f\in X_c(L^{\infty},\ell^2)$.
So, we have $\sum_{n\in \Bbb Z}f(x-na)\overline{h(x-na)}=\sum_{n\in \Bbb Z}f(x-nc)\overline{\tilde{h}(x-nc)}$ for all $f$ in $X_a(L^{\infty},\ell^2)=X_c(L^{\infty},\ell^2)$; in particular, for all $f\in C_c(\Bbb R)$. This is enough to conclude $h=\tilde{h}$; that is, $h\in L_c^{\infty}(\ell^2)$ which proves $L_a^{\infty}(\ell^2)\subseteq L_c^{\infty}(\ell^2)$. The reverse inclusion is proved in the same way; hence, $L_a^{\infty}(\ell^2)= L_c^{\infty}(\ell^2)$. By the preceding part of the proof this implies $\frac{a}{c}\in \Bbb Q$.
\qed

\vspace{.2in}

It is well known that the Wiener spaces $W(L^{\infty},\ell^1)$ and $W(L^{\infty},\ell^2)$ are translation invariant. The same is true for our Hilbert $C^*$-modules   $X_{a}(L^{\infty},\ell^2)$ and $L_{a}^{\infty}(\ell^2)$.

\vspace{.1in}

\begin{prop}\label{translation invariant}
Let $a>0$. Then $X_{a}(L^{\infty},\ell^2)$ and $L_{a}^{\infty}(\ell^2)$ are invariant under all translations $T_c$, $c\in \Bbb R$. Each translation $T_c$ is an isometry on $L_{a}^{\infty}(\ell^2)$.
\end{prop}
\proof
Fix $a>0$ and $c\in \Bbb R$. Let $f \in L_{a}^{\infty}(\ell^2)$ and $x\in [0,a]$. Consider $x-c$ and find an integer $n_0$ with the property $x-c-n_0a=y\in [0,a]$. Then we have
$$
\sum_{n\in \Bbb Z}|T_cf(x-na)|^2=\sum_{n\in \Bbb Z}|f(x-c-n_0a+n_0a-na)|^2=\sum_{n\in \Bbb Z}|f(y-(n-n_0)a)|^2=\sum_{n^{\prime}\in \Bbb Z}|f(y-n^{\prime}a)|^2.
$$
From this we conclude that $T_cf\in L_{a}^{\infty}(\ell^2)$ and $\|T_cf\|_{L_{a}^{\infty}(\ell^2)}\leq \|f\|_{L_{a}^{\infty}(\ell^2)}$. The same conclusion applied for $T_{-c}$ and $T_cf$ gives us the opposite inequality:
$\|f\|_{L_{a}^{\infty}(\ell^2)}=\|T_{-c}T_cf\|_{L_{a}^{\infty}(\ell^2)}\leq \|T_cf\|_{L_{a}^{\infty}(\ell^2)}$.

Let us now take any $f\in X_{a}(L^{\infty},\ell^2)$. This means that the series $\sum_{n\in \Bbb Z}\left| T_{-na}f\chi_{[0,a]}\right|^2$ converges in norm in $L^{\infty}[0,a]$ and we need to show that this implies norm convergence in $L^{\infty}[0,a]$ of the series $\sum_{n\in \Bbb Z}\left| T_{-na}T_cf\chi_{[0,a]}\right|^2$.
Recall that is enough to consider only symmetric partial sums since if a series of positive elements in a $C^*$-algebra converges, it converges unconditionally.

In the $n$-th summand we have the function $|f(x+na-c)|^2$ for $x\in [0,a]$; i.e.~$|f(y)|^2$ for $na-c\leq y\leq (n+1)a-c$. Observe that
\begin{equation}\label{pomaknute sume}
[na-c,(n+1)a-c]\subseteq [(n^{\prime}-1)a,n^{\prime}a]\cup [n^{\prime}a,(n^{\prime}+1)a],
\end{equation}
where $n-n^{\prime}=d_{a,c}$ and the integer $d_{a,c}$ is independent of $n$ and depends only on $a$ and $c$.

Let $\epsilon >0$ be given. Since $\sum_{n\in \Bbb Z}\left| T_{-na}f\chi_{[0,a]}\right|^2$ converges in norm in $L^{\infty}[0,a]$ there exists a natural number $N_0$ with the property
$$
N_2>N_1\geq N_0 \Longrightarrow \left\|\sum_{N_1+1\leq|n|\leq N_2}\left| T_{-na}f\chi_{[0,a]}\right|^2\right\|_{L^{\infty}[0,a]}<\epsilon.
$$
From this we conclude, by letting $N_2\rightarrow \infty$, that
$$
\left\|\sum_{|n|> N_0}\left| T_{-na}f\chi_{[0,a]}\right|^2\right\|_{L^{\infty}[0,a]}\leq\epsilon.
$$
Now for $N_2>N_1\geq N_0+|d_{a,c}|$ we have
\begin{eqnarray*}
\left\|\sum_{N_1+1\leq|n|\leq N_2}\left| T_{-na}T_cf\chi_{[0,a]}\right|^2\right\|_{L^{\infty}[0,a]}&\leq&\left\|2\sum_{N_1+1-d_{a,c}\leq|n^{\prime}|\leq N_2-d_{a,c}}\left| T_{-n^{\prime}a}f\chi_{[0,a]}\right|^2\right\|_{L^{\infty}[0,a]}\\
 &\leq&\left\|2\sum_{|n^{\prime}|> N_0}\left| T_{-n^{\prime}a}f\chi_{[0,a]}\right|^2\right\|_{L^{\infty}[0,a]}\\
 &\leq&2\epsilon
\end{eqnarray*}
\qed

\vspace{.1in}

\begin{remark}\label{not adjointable}
We note that $T_c$ for $c\not = na$, $n\in \Bbb Z$, is not a module map on $L_{a}^{\infty}(\ell^2)$. In particular, this shows us that $T_c$ cannot be an adjointable operator on the Hilbert $C^*$-module $L_{a}^{\infty}(\ell^2)$, except for those $c$ that are integer multiples of $a$.
\end{remark}

\vspace{.2in}

We are now in position to apply our results from Sections 3 and 4 in Gabor analysis.

\vspace{.1in}

Suppose we are given a function $g\in L^2(\Bbb R)$ and $a,b>0$. Let
\begin{equation}\label{Gabor gammas}
\Gamma_{jk}(x)=\langle T_{ka}g,T_{ja}g\rangle_{\frac{1}{b}}(x)=\sum_{l\in \Bbb Z}g(x-ka-\frac{l}{b})\overline{g(x-ja-\frac{l}{b})},\,\,\,j,k\in \Bbb Z,
\end{equation}
 and
\begin{equation}\label{Gabor gammas 1}
\Gamma_j=\Gamma_{j0}(x)=\langle g,T_{ja}g\rangle_{\frac{1}{b}}(x)=\sum_{l\in \Bbb Z}g(x-\frac{l}{b})\overline{g(x-ja-\frac{l}{b})},\,\,\,j\in \Bbb Z.
\end{equation}
In addition, let
\begin{equation}\label{Gabor gs 1}
G_k(x)=\langle g,T_{\frac{k}{b}}g\rangle_a(x)=\sum_{n\in \Bbb Z}g(x-na)\overline{g(x-na-\frac{k}{b})},\,\,\,k\in \Bbb Z.
\end{equation}

As we observed in Section 5 (recall the equation \eqref{first inner}), the functions $\Gamma_{jk}$, $\Gamma_{j}$, and $G_k$ are well defined for a.e.~$x$.

\vspace{.1in}

The following theorem is known. We include the proof to demonstrate how the theory of weak Bessel modular sequences applies.

\begin{theorem}\label{CC revisited}
Suppose that for $g\in L^2(\Bbb R)$ and $a,b>0$ either of the following two conditions is satisfied:
\begin{equation}\label{dual cc}
\sup_{N\in \Bbb N}\left\| \sum_{j=-N}^N |\Gamma_j|\right\|_{L^{\infty}[0,\frac{1}{b}]}\leq B,
\end{equation}
\begin{equation}\label{original cc}
\sup_{N\in \Bbb N}\left\| \sum_{k=-N}^N |G_k|\right\|_{L^{\infty}[0,a]}\leq B.
\end{equation}
Then $G(g,a,b)$ is a Bessel sequence in $L^2(\Bbb R)$ with a Bessel bound $\frac{1}{b}B$.
\end{theorem}
\proof
Suppose that \eqref{dual cc} is satisfied. By Theorem \ref{CasCocLem} we need to show that
$(T_{na}(\frac{1}{\sqrt{b}}g))_{n\in \Bbb Z}$ is a weak Bessel sequence with a Bessel bound $\frac{1}{b}B$ in the Hilbert $L^{\infty}[0,\frac{1}{b}]$-module $L_{\frac{1}{b}}^{\infty}(\ell^2)$. This is equivalent to the property that $(T_{na}g)_{n\in \Bbb Z}$ is a weak Bessel sequence with a Bessel bound $B$.
By Proposition \ref{Bessel by Gram} it suffices to prove that the Gram matrix $\Gamma$ of the sequence $(T_{na}g)_{n\in \Bbb Z}$  defines a bounded module map $\ell^2(L^{\infty}[0,\frac{1}{b}]) \rightarrow \ell^2_{\text{strong}}(L^{\infty}[0,\frac{1}{b}])$ with a bound $B$.

Notice that \eqref{dual cc} implies in particular that $\Gamma_0\in L^{\infty}[0,\frac{1}{b}]$ which means that $g\in L_{\frac{1}{b}}^{\infty}(\ell^2)$. This in turn implies by Proposition \ref{translation invariant} that all $T_{na}g$ belong to $L_{\frac{1}{b}}^{\infty}(\ell^2)$ and hence formulae \eqref{Gabor gammas} and \eqref{Gabor gammas 1} can be rewritten with $\langle \cdot,\cdot\rangle_{L_{\frac{1}{b}}^{\infty}(\ell^2)}$ instead of $\langle \cdot,\cdot\rangle_{\frac{1}{b}}$.

Observe now that the matrix coefficients of the Gramm matrix $\Gamma$ are
$$
\langle T_{ka}g,T_{ja}g\rangle_{L_{\frac{1}{b}}^{\infty}(\ell^2)}(x)=\sum_{l\in \Bbb Z}g(x-ka-\frac{l}{b})\overline{g(x-ja-\frac{l}{b})}=\Gamma_{jk}(x),\,\,\,j,k\in \Bbb Z.
$$
By Proposition \ref{schur drugi} it is enough to see that the matrix coefficients $\Gamma_{jk}$ satisfy condition \eqref{schur5}:
\begin{equation}\label{schur5 gabor}
\left\|\text{(strong)}\sum_{j\in \Bbb Z}|\Gamma_{jk}|\right\|_{L^{\infty}[0,\frac{1}{b}]}\leq B,\,\,\forall k\in \Bbb Z \mbox{ and }\left\|\text{(strong)}\sum_{k\in \Bbb Z}|\Gamma_{jk}|\right\|_{L^{\infty}[0,\frac{1}{b}]}\leq B,\,\,\forall j\in \Bbb Z.
\end{equation}
First observe that $\Gamma_{kj}=\overline{\Gamma_{jk}}$ for all $k$ and $j$. Therefore, it is enough to check the first inequality in \eqref{schur5 gabor}. Secondly, we have
\begin{eqnarray*}
\Gamma_{jk}(x)&=&\langle T_{ka}g,T_{ja}g\rangle_{L_{\frac{1}{b}}^{\infty}(\ell^2)}(x)\\
 &=&\sum_{l\in \Bbb Z}g(x-ka-\frac{l}{b})\overline{g(x-ja-\frac{l}{b})}\\
 &=&\sum_{l\in \Bbb Z}g(x-ka-\frac{l}{b})\overline{g(x-ka-(j-k)a-\frac{l}{b})}\\
 &=&\Gamma_{j-k}(x-ka).
\end{eqnarray*}
Since all $\Gamma_{j-k}$ are periodic functions, our assumption \eqref{dual cc} now implies the desired first inequality in \eqref{schur5 gabor}.

Let us now assume \eqref{original cc}. Observe that this means that the function $g$ satisfies \eqref{dual cc} with parameters $\frac{1}{b}$ and $\frac{1}{a}$ playing the roles of $a$ and $b$, respectively. Hence, by the first part of the proof the sequence $(T_{na}(\sqrt{a}g))_{n\in \Bbb Z}$ is a weak Bessel sequence with a Bessel bound $aB$ in the Hilbert $L^{\infty}[0,a]$-module $L_a^{\infty}(\ell^2)$. By Theorem \ref{CasCocLem} this means that the sequence $G(g,\frac{1}{b},\frac{1}{a})$ is Bessel with a Bessel bound $aB$. Finally, by Lemma 9.2.2 from \cite{C} we conclude that $G(g,a,b)$ is Bessel with a Bessel bound $\frac{1}{ab}aB=\frac{1}{b}B$.
\qed

\vspace{.1in}

\begin{remark}\label{evo cc}
Suppose that $g$ satisfies \eqref{original cc}. This means that the sequence of symmetric partial sums of the series $\sum_{k\in \Bbb Z}|G_k|$ is bounded which is equivalent to the strong convergence of this series, which is in turn, by Remark \ref{strong vs measure}, equivalent to its pointvise a.e. convergence to a function in $L^{\infty}[0,a]$.
Hence, we may rewrite \eqref{original cc} in the form
\begin{equation}\label{really original cc}
\esssup_{x\in [0,a]} \sum_{k\in \Bbb Z} |G_k(x)|\leq B.
\end{equation}
which is what is usually called the CC condition (see \cite{C}, Section 9.1). Therefore if $g$ satisfies \eqref{original cc} with the  parameters $a$ and $b$ the above theorem actually restates the first assertion of Theorem 8.4.4 from \cite{C}.

We also note that a stronger assumption, namely
\begin{equation}\label{really original Daubechies}
\sum_{k\in \Bbb Z} \|G_k\|_{L^{\infty}[0,a]}\leq B
\end{equation}
ensures the same conclusion. To see this, simply observe that \eqref{really original Daubechies} implies \eqref{really original cc} or use Proposition \ref{schur prvi}. The condition \eqref{really original Daubechies} is well known (cf.~equation (8.13) and Theorem 8.4.1 in \cite{C}); it was first used in the late 1980's by I.~Daubechies.
\end{remark}

\vspace{.1in}

We now turn to the frame operator.

Suppose that $G(g,a,b)$ and $G(h,a,b)$ are Bessel sequences for some $g,h\in L^2(\Bbb R)$ and $a,b>0$. Then $(T_{na}(\frac{1}{\sqrt{b}}g))_{n\in \Bbb Z}$ and $(T_{na}(\frac{1}{\sqrt{b}}h))_{n\in \Bbb Z}$ are weak Bessel sequences in the Hilbert $L^{\infty}[0,\frac{1}{b}]$-module $L_{\frac{1}{b}}^{\infty}(\ell^2)$. Denote by $U$ and $V$ the corresponding analysis operators. To avoid confusion we will denote by $U_0$ and $V_0$ the analysis operators of the original sequences $G(g,a,b)$ and $G(h,a,b)$. By Theorem \ref{Bessel relaxed} we have
\begin{equation}\label{Bessel-adjoint Gabor}
V^*Uf=(\text{weak-strong})\,\frac{1}{b}\sum_{n\in \Bbb Z}\langle f,T_{na}g\rangle_{L_{\frac{1}{b}}^{\infty}(\ell^2)}T_{na}h.
\end{equation}

This means that
\begin{equation}\label{BaG advanced}
\langle V^*Uf,p\rangle_{L_{\frac{1}{b}}^{\infty}(\ell^2)}=(\text{strong})\,\frac{1}{b}\sum_{n\in \Bbb Z}\langle f,T_{na}g\rangle_{L_{\frac{1}{b}}^{\infty}(\ell^2)}
\langle T_{na}h,p\rangle_{L_{\frac{1}{b}}^{\infty}(\ell^2)},\,\,\, \forall f,p\in L_{\frac{1}{b}}^{\infty}(\ell^2).
\end{equation}
By Remark \ref{strong vs measure} we obtain the pointwise a.e. convergence for all $f,p\in L_{\frac{1}{b}}^{\infty}(\ell^2)$:
\begin{equation}\label{BaG advanced pointwise}
\langle V^*Uf,p\rangle_{L_{\frac{1}{b}}^{\infty}(\ell^2)}(x)=\frac{1}{b}\sum_{n\in \Bbb Z}\langle f,T_{na}g\rangle_{L_{\frac{1}{b}}^{\infty}(\ell^2)}(x)
\langle T_{na}h,p\rangle_{L_{\frac{1}{b}}^{\infty}(\ell^2)}(x) \mbox{ for a.e. }x\in [0,\frac{1}{b}].
\end{equation}

Denote by $e$ the constant function $1$ on the interval  $[0,\frac{1}{b}]$. Since $e$ is the unit element in our von Neumann algebra $L^{\infty}[0,\frac{1}{b}]$, we know from Remark \ref{prva baza} that the sequence $(e^{(n)})_{n\in \Bbb Z}$ is the canonical weak basis for $\ell^2_{\Bbb Z,\text{strong}}(L^{\infty}[0,\frac{1}{b}])$. Theorem \ref{dual discovery} tells us now that the sequence $(\chi_{[\frac{n}{b},\frac{n+1}{b}]})_{n\in \Bbb Z}$ is the canonical weak basis in the Hilbert $L^{\infty}[0,\frac{1}{b}]$-module $L_{\frac{1}{b}}^{\infty}(\ell^2)$. In our next proposition we compute the matrix of the frame operator in this basis.

Let
\begin{equation}\label{joint Gk}
G_k^{h,g}(x)=\left(\sum_{n\in \Bbb Z}T_{na}hT_{\frac{k}{b}}T_{na}\overline{g}\right)(x)=\sum_{n\in \Bbb Z}h(x-na)\overline{g(x-na-\frac{k}{b})},\,\,\,k\in \Bbb Z.
\end{equation}
Recall from Remark \ref{taman} that $g,h\in L_a^{\infty}(\ell^2) \cap L_{\frac{1}{b}}^{\infty}(\ell^2)$, so we can write
\begin{equation}\label{joint Gk bis}
G_k^{h,g}(x)=\langle h,T_{\frac{k}{b}}g\rangle_{L_a^{\infty}(\ell^2)},\,\,\,k\in \Bbb Z.
\end{equation}
Observe also that for $h=g$ we have $G_k^{h,g}=G_k$ where $G_k$ are introduced in \eqref{Gabor gs 1}.

\vspace{.1in}

\begin{prop}\label{matrica frame operatora prop}
Let $g,h\in L^2(\Bbb R)$ and $a,b>0$ be such that $G(g,a,b)$ and $G(h,a,b)$ are Bessel sequences in $L^2(\Bbb R)$. Denote by $U$ and $V$ the analysis operators of  weak Bessel sequences $(T_{na}(\frac{1}{\sqrt{b}}g))_{n\in \Bbb Z}$ and $(T_{na}(\frac{1}{\sqrt{b}}h))_{n\in \Bbb Z}$  in the Hilbert $L^{\infty}[0,\frac{1}{b}]$-module $L_{\frac{1}{b}}^{\infty}(\ell^2)$. Let $(m_{kj})$ be the matrix of the operator $V^*U$ in the canonical weak basis $(\chi_{[\frac{n}{b},\frac{n+1}{b}]})_{n\in \Bbb Z}$ of
$L_{\frac{1}{b}}^{\infty}(\ell^2)$. Then
\begin{equation}\label{matrica cross operatora}
m_{kj}(x)=\frac{1}{b}G_{k-j}^{h,g}(x+\frac{k}{b}),\,\,\,x\in [0,\frac{1}{b}],\,k,j\in \Bbb Z.
\end{equation}
In particular, if $(s_{kj})$ is the matrix of the frame operator $S=U^*U$ of $(T_{na}(\frac{1}{\sqrt{b}}g))_{n\in \Bbb Z}$ with respect to $(\chi_{[\frac{n}{b},\frac{n+1}{b}]})_{n\in \Bbb Z}$, then
\begin{equation}\label{matrica frame operatora}
s_{kj}(x)=\frac{1}{b}G_{k-j}(x+\frac{k}{b}),\,\,\,x\in [0,\frac{1}{b}],\,k,j\in \Bbb Z.
\end{equation}
\end{prop}
\proof
We first note that, for any $f\in L_{\frac{1}{b}}^{\infty}(\ell^2)$ and $j\in \Bbb Z$,
\begin{equation}\label{matricni koeficijenti}
\langle f,\chi_{[\frac{j}{b},\frac{j+1}{b}]}\rangle_{L_{\frac{1}{b}}^{\infty}(\ell^2)}(x)=f(x+\frac{j}{b}),\,\,\,x\in [0,\frac{1}{b}].
\end{equation}
Indeed, we have
$$
\langle f,\chi_{[\frac{j}{b},\frac{j+1}{b}]}\rangle_{L_{\frac{1}{b}}^{\infty}(\ell^2)}(x)=\sum_{k\in \Bbb Z}f(x-\frac{k}{b})\chi_{[\frac{j}{b},\frac{j+1}{b}]}(x-\frac{k}{b}) ,\,\,\,x\in [0,\frac{1}{b}].
$$
Clearly, for all $k\not = -j$ the corresponding terms vanish, while for $k = -j$ we get $f(x+\frac{j}{b})$.

Now we have for all $k,j\in\Bbb Z$ and $x\in [0,\frac{1}{b}]$
\begin{eqnarray*}
m_{kj}(x)&=&\langle V^*U\chi_{[\frac{j}{b},\frac{j+1}{b}]},\chi_{[\frac{k}{b},\frac{k+1}{b}]}\rangle_{L_{\frac{1}{b}}^{\infty}(\ell^2)}(x)\\
 &\stackrel{\eqref{BaG advanced pointwise}}{=}&\frac{1}{b}\sum_{n\in \Bbb Z}\langle \chi_{[\frac{j}{b},\frac{j+1}{b}]},T_{na}g\rangle_{L_{\frac{1}{b}}^{\infty}(\ell^2)}(x)
\langle T_{na}h,\chi_{[\frac{k}{b},\frac{k+1}{b}]}\rangle_{L_{\frac{1}{b}}^{\infty}(\ell^2)}(x)\\
 &\stackrel{\eqref{matricni koeficijenti}}{=}&\sum_{n\in \Bbb Z}h(x-na+\frac{k}{b})\overline{g(x-na+\frac{j}{b})}\\
 &\stackrel{\eqref{joint Gk}}{=}&\frac{1}{b}G_{k-j}^{h,g}(x+\frac{k}{b}).
\end{eqnarray*}
\qed

\vspace{.1in}

\begin{remark}\label{karakterizacija Gabor Parseval}
Now one can easily reobtain a well known characterization of Parseval Gabor frames. By a result from \cite{RS2} (see also Theorem 3.2 in \cite{CCJ1}), a sequence $G(g,a,b)$ is a Parseval frame in $L^2(\Bbb R)$ if and only if the following two conditions are satisfied:
\begin{itemize}
\item[(i)] $G_0(x)=b$ a.e.,
\item[(ii)] $G_k(x)=0$ a.e.~for all $k\not =0$.
\end{itemize}
To see this, assume first that $G(g,a,b)$ is a Parseval frame. Then by Theorem \ref{CasCocLem}
$(T_{na}(\frac{1}{\sqrt{b}}g))_{n\in \Bbb Z}$ is a weak Parseval frame in $L_{\frac{1}{b}}^{\infty}(\ell^2)$. This implies that $U^*U$ is the identity operator on $L_{\frac{1}{b}}^{\infty}(\ell^2)$ and hence its matrix coefficients $s_{jk}$ satisfy $s_{jk}=\delta_{jk}e$ for all $j,k$ (with the Kronecker $\delta_{jk}$). Now the preceding proposition gives us (i) and (ii).

Conversely, if (i) and (ii) are satisfied we first conclude using Theorem \ref{CC revisited} that $G(g,a,b)$ is a Bessel sequence and hence by Theorem \ref{CasCocLem}
$(T_{na}(\frac{1}{\sqrt{b}}g))_{n\in \Bbb Z}$ is a weak Bessel sequence. In particular, we have a well defined and bounded analysis operator $U$ and now using (i) and (ii) and the preceding proposition we see that the frame operator $S=U^*U$ is the identity operator. Hence $(T_{na}(\frac{1}{\sqrt{b}}g))_{n\in \Bbb Z}$ is a weak Parseval frame and applying Theorem \ref{CasCocLem} again we conclude that $G(g,a,b)$ is a Parseval frame.
\end{remark}

\vspace{.2in}

Another important feature of the correspondence of Bessel Gabor sequences in $L^2(\Bbb R)$ and weak Bessel sequences in $L_{\frac{1}{b}}^{\infty}(\ell^2)$ is that it preserves duality.
We first prove a lemma that is important in its own.

\vspace{.1in}

\begin{lemma}\label{frame operatori isti}
Let $G(g,a,b)$ and $G(h,a,b)$ be Gabor Bessel sequences in $L^2(\Bbb R)$ with the analysis operators $U_0$ and $V_0$. Denote by $U$ and $V$ the analysis operators of the corresponding weak Bessel sequences $(T_{na}(\frac{1}{\sqrt{b}}g))_{n\in \Bbb Z}$ and $(T_{na}(\frac{1}{\sqrt{b}}h))_{n\in \Bbb Z}$.
Then $V_0^*U_0f=V^*Uf$ for every $f\in L_{\frac{1}{b}}^{\infty}(\ell^2)$.
\end{lemma}
\proof
Consider again the equation \eqref{BaG advanced pointwise}:
$$
\langle V^*Uf,p\rangle_{L_{\frac{1}{b}}^{\infty}(\ell^2)}(x)=\frac{1}{b}\sum_{n\in \Bbb Z}\langle f,T_{na}g\rangle_{L_{\frac{1}{b}}^{\infty}(\ell^2)}(x)
\langle T_{na}h,p\rangle_{L_{\frac{1}{b}}^{\infty}(\ell^2)}(x) \mbox{ for a.e. }x\in [0,\frac{1}{b}]
$$
which holds for all $f,p\in L_{\frac{1}{b}}^{\infty}(\ell^2)$.
If we take $p=\chi_{[\frac{j}{b},\frac{j+1}{b}]}$ then, using \eqref{matricni koeficijenti}, we get for all $f\in L_{\frac{1}{b}}^{\infty}(\ell^2)$ and every $j\in \Bbb Z$
$$
V^*Uf(x+\frac{j}{b})=\frac{1}{b}\sum_{n\in \Bbb Z}\langle f,T_{na}g\rangle_{L_{\frac{1}{b}}^{\infty}(\ell^2)}(x)T_{na}h(x+\frac{j}{b}), \mbox{ for a.e. }x\in [0,\frac{1}{b}].
$$
Since $\langle f,T_{na}g\rangle_{L_{\frac{1}{b}}^{\infty}(\ell^2)}$ is a $\frac{1}{b}$-periodic, we can write
$\langle f,T_{na}g\rangle_{L_{\frac{1}{b}}^{\infty}(\ell^2)}(x)=\langle f,T_{na}g\rangle_{L_{\frac{1}{b}}^{\infty}(\ell^2)}(x+\frac{j}{b})$ so that the preceding equality can be written as
\begin{equation}\label{frame operator tockovno}
V^*Uf(x)=\frac{1}{b}\sum_{n\in \Bbb Z}\langle f,T_{na}g\rangle_{L_{\frac{1}{b}}^{\infty}(\ell^2)}(x)T_{na}h(x), \mbox{ for a.e. }x\in \Bbb R,\,\forall f\in L_{\frac{1}{b}}^{\infty}(\ell^2).
\end{equation}
We now recall Theorem 4.6.8 from \cite{CasL} which states that
\begin{equation}\label{cas lammers frame operator}
V_0^*U_0f=\frac{1}{b}\sum_{n\in \Bbb Z}\langle f,T_{na}g\rangle_{\frac{1}{b}}T_{na}h,\,\,\,\forall f \in L^2(\Bbb R)
\end{equation}
unconditionally.
In fact, this formula is stated and proved in \cite{CasL} only for $h=g$, but an inspection of the proof shows that the same result is valid in this more general form for two Gabor windows $g$ and $h$.
Observe a subtle difference between \eqref{Bessel-adjoint Gabor} (which led us to \eqref{frame operator tockovno}) and \eqref{cas lammers frame operator}. In \eqref{cas lammers frame operator} we have convergence in $\|\cdot\|_{L^2(\Bbb R)}$ but we cannot conclude a similar relation in $\|\cdot\|_{L_{\frac{1}{b}}^{\infty}(\ell^2)}$ since $\|f\|_{L^2(\Bbb R)}\leq \frac{1}{b}\|f\|_{L_{\frac{1}{b}}^{\infty}(\ell^2)}$ for all $f$ in $L_{\frac{1}{b}}^{\infty}(\ell^2)$. In addition, note that in \eqref{cas lammers frame operator} we have the $\frac{1}{b}$-product $\langle f,T_{na}g\rangle_{\frac{1}{b}}$ defined in \eqref{first inner}. Although the defining formula for $\langle \cdot,\cdot\rangle_{L_{\frac{1}{b}}^{\infty}(\ell^2)}$ and $\langle \cdot,\cdot\rangle_{{\frac{1}{b}}}$
is the same, $\langle \cdot,\cdot\rangle_{{\frac{1}{b}}}$ makes sense for all functions from $L^2(\Bbb R)$, but the result is a function in $L^1[0,\frac{1}{b}]$ and not necessarily in $L^{\infty}[0,\frac{1}{b}]$.

However, by comparing \eqref{frame operator tockovno} and \eqref{cas lammers frame operator} we conclude that \eqref{cas lammers frame operator} also holds pointwise for all functions from $L_{\frac{1}{b}}^{\infty}(\ell^2)$, i.e.
\begin{equation}\label{cas lammers frame operator tockovno}
V_0^*U_0f(x)=\frac{1}{b}\sum_{n\in \Bbb Z}\langle f,T_{na}g\rangle_{\frac{1}{b}}(x)T_{na}h(x),\mbox{ for a.e. }x\in \Bbb R,\,\forall f \in L_{\frac{1}{b}}^{\infty}(\ell^2)
\end{equation}
A comparison of the last equality and \eqref{frame operator tockovno} gives us the desired conclusion.
\qed

\vspace{.2in}

\begin{remark}\label{klasicno skoro svuda}
Let us retain notations from the preceding proof. We claim that \eqref{cas lammers frame operator tockovno} holds for all $f\in L^2(\Bbb R)$, i.e.~that
\begin{equation}\label{cas lammers frame operator tockovno 1}
V_0^*U_0f(x)=\frac{1}{b}\sum_{n\in \Bbb Z}\langle f,T_{na}g\rangle_{\frac{1}{b}}(x)T_{na}h(x),\mbox{ for a.e. }x\in \Bbb R,\,\forall f \in L^2(\Bbb R).
\end{equation}
To see this take any $f\in L^2(\Bbb R)$ and define
$$
f_0(x)=\left\{
\begin{array}{cl}0,&\mbox{if }\langle f,f\rangle_{\frac{1}{b}}(x)=0\\
\frac{f(x)}{\langle f,f\rangle_{\frac{1}{b}}(x)^{1/2}},&\mbox{if }\langle f,f\rangle_{\frac{1}{b}}(x)\not=0
\end{array}\right..
$$
Obviously, $f_0\in L_{\frac{1}{b}}^{\infty}(\ell^2)$ and \eqref{cas lammers frame operator tockovno} applied to $f_0$ yields \eqref{cas lammers frame operator tockovno 1}.
\end{remark}

\vspace{.1in}

\begin{prop}\label{dualnost ocuvana}
Let $G(g,a,b)$ and $G(h,a,b)$ be Gabor Bessel sequences in $L^2(\Bbb R)$. $G(g,a,b)$ and $G(h,a,b)$ are mutually dual if and only if the weak Bessel sequences $(T_{na}(\frac{1}{\sqrt{b}}g))_{n\in \Bbb Z}$ and $(T_{na}(\frac{1}{\sqrt{b}}h))_{n\in \Bbb Z}$ are dual to each other.
\end{prop}
\proof
Denote again by $U_0$ and $V_0$ the original analysis operators and by $U$ and $V$ the corresponding modular analysis operators. We must show that $V^*U=I_{L_{\frac{1}{b}}^{\infty}(\ell^2)}$ if and only $V_0^*U_0=I_{L^2(\Bbb R)}$, where $I$ denotes the identity operator on the indicated ambient space.

Suppose first that $G(g,a,b)$ and $G(h,a,b)$ are dual to each other. Then we have $V_0^*U_0=I_{L^2(\Bbb R)}$
which by Lemma \ref{frame operatori isti} immediately implies $V^*U=I_{L_{\frac{1}{b}}^{\infty}(\ell^2)}$.

Conversely, if $(T_{na}(\frac{1}{\sqrt{b}}g))_{n\in \Bbb Z}$ and $(T_{na}(\frac{1}{\sqrt{b}}h))_{n\in \Bbb Z}$ are dual to each other we have $V^*U=I_{L_{\frac{1}{b}}^{\infty}(\ell^2)}$. By Lemma \ref{frame operatori isti} we conclude that $V_0^*U_0f=f$ for all $f \in L_{\frac{1}{b}}^{\infty}(\ell^2)$.
This is enough to conclude $V_0^*U_0=I_{L^2(\Bbb R)}$ since these two bounded operators coincide on $L_{\frac{1}{b}}^{\infty}(\ell^2)$ which is a dense subspace of
$(L^2(\Bbb R)$ with respect to $\|\cdot\|_{L^2(\Bbb R)})$.
\qed

\vspace{.2in}

\begin{remark}\label{WR}
We are now again in position to show how an important classical result can easily be reobtained by passing from Gabor Bessel sequences to the corresponding weak Bessel sequences.

The Wexler-Raz theorem (\cite{C}, Theorem 9.3.5) states that Gabor Bessel sequences $G(g,a,b)$ and $G(h,a,b)$ in $L^2(\Bbb R)$ are dual to each other if and only if the following two conditions are satisfied:
\begin{itemize}
\item[(i)] $\langle h,M_{\frac{m}{a}}T_{\frac{k}{b}}g\rangle_{L^2(\Bbb R)}=0$ for all $m,k\in \Bbb Z$ such that $m^2+k^2>0$,
\item[(ii)] $\langle h,g\rangle_{L^2(\Bbb R)}=ab$.
\end{itemize}
To see this, first recall a result from \cite{DDR1}; see also Proposition 4.4.5 in \cite{CasL}:
if $g,h$ are in $L^2(\Bbb R)$ then
\begin{equation}\label{ortogonalost u a produktu}
\text{span}\,\{M_{\frac{m}{a}}h:m\in \Bbb Z\} \perp \text{span}\,\{M_{\frac{m}{a}}g:m\in \Bbb Z\} \Leftrightarrow \langle h,g\rangle_{a}=0 \mbox{ a.e.}
\end{equation}
where $\perp$ indicates orthogonality with respect to the inner product in ${L^2(\Bbb R)}$.
To see this, just observe that using the periodization trick we have
$$
\langle h,M_{\frac{m}{a}}g\rangle_{L^2(\Bbb R)}=\int_0^a\langle h,g\rangle_{a}(x)\text{e}^{-2\pi i \frac{m}{a}x}dx.
$$
Suppose now that $G(g,a,b)$ and $G(h,a,b)$ are Bessel sequences dual to each other. Then by the preceding proposition weak frames
$(T_{na}(\frac{1}{\sqrt{b}}g))_{n\in \Bbb Z}$ and $(T_{na}(\frac{1}{\sqrt{b}}h))_{n\in \Bbb Z}$ in $L_{\frac{1}{b}}^{\infty}(\ell^2)$ are dual to each other. If $U$ and $V$ again denote the corresponding analysis operators, we have $V^*U=I_{L_{\frac{1}{b}}^{\infty}(\ell^2)}$. By Proposition \ref{matrica frame operatora prop} we conclude that
$$
\frac{1}{b}G_{k-j}^{h,g}(x+\frac{k}{b})=\left\{\begin{array}{ll}0&\mbox{if }k\not = j\\1&\mbox{if }k= j\end{array}\right.,\mbox{ for a.e. }x\in [0,\frac{1}{b}]
$$
where $G_k^{h,g}(x)=\sum_{n\in \Bbb Z}h(x-na)\overline{g(x-na-\frac{k}{b})}=\langle h,T_{\frac{k}{b}}g\rangle_{a}$ are the functions introduced in \eqref{joint Gk}.
This implies $G_k^{h,g}(x)=0$ a.e.~for all $k\not=0$ and $G_0^{h,g}(x)=b$ a.e. Using \eqref{ortogonalost u a produktu} we immediately get (i) and (ii).

Conversely, (i) and (ii) imply via \eqref{ortogonalost u a produktu} $G_k^{h,g}(x)=0$ a.e.~for all $k\not=0$. For $k=0$ we have $\int_0^aG_0^{h,g}(x)dx=ab$ and $\int_0^aG_0^{h,g}(x)\text{e}^{-2\pi i \frac{m}{a}x}dx=ab$ for all $m\not = 0$. This is enough to conclude that $G_0^{h,g}(x)=b$ a.e.  Proposition \ref{matrica frame operatora prop} now implies that $V^*U=I_{L_{\frac{1}{b}}^{\infty}(\ell^2)}$.
\end{remark}

\vspace{.2in}

We now turn to the Walnut representation.

Suppose again that $G(g,a,b)$ and $G(h,a,b)$ are Bessel sequences in $L^2(\Bbb R)$ for some $g,h\in L^2(\Bbb R)$ and $a,b>0$. Then $(T_{na}(\frac{1}{\sqrt{b}}g))_{n\in \Bbb Z}$ and $(T_{na}(\frac{1}{\sqrt{b}}h))_{n\in \Bbb Z}$ are weak Bessel sequences in the Hilbert $L^{\infty}[0,\frac{1}{b}]$-module $L_{\frac{1}{b}}^{\infty}(\ell^2)$. Denote by $U$ and $V$ the corresponding analysis operators.
We first recall the equation \eqref{Bessel-adjoint Gabor}:
$$
V^*Uf=(\text{weak-strong})\,\frac{1}{b}\sum_{n\in \Bbb Z}\left(\sum_{k\in \Bbb Z}T_{\frac{k}{b}}fT_{\frac{k}{b}}T_{na}\overline{g}\right)T_{na}h
,\,\,\, \forall f\in L_{\frac{1}{b}}^{\infty}(\ell^2).
$$
Since here
both series converge unconditionally we may write
\begin{equation}\label{Bessel-adjoint Gabor 1}
V^*Uf=(\text{weak-strong})\,\frac{1}{b}\sum_{k\in \Bbb Z}\left(\sum_{n\in \Bbb Z}T_{na}hT_{\frac{k}{b}}T_{na}\overline{g}\right)T_{\frac{k}{b}}f,\,\,\, \forall f\in L_{\frac{1}{b}}^{\infty}(\ell^2).
\end{equation}
Using \eqref{joint Gk} we may rewrite \eqref{Bessel-adjoint Gabor 1} as
\begin{equation}\label{Bessel-adjoint Gabor 2}
V^*Uf=(\text{weak-strong})\,\frac{1}{b}\sum_{k\in \Bbb Z}G_k^{h,g}T_{\frac{k}{b}}f,\,\,\, \forall f\in L_{\frac{1}{b}}^{\infty}(\ell^2).
\end{equation}
We know from \eqref{joint Gk bis} that $G_k^{h,g} \in L^{\infty}[0,a]$. Extended by $a$-periodicity these functions can be viewed as elements in $L^{\infty}(\Bbb R)$. In particular, in \eqref{Bessel-adjoint Gabor 2} the functions $G_k^{h,g}$ are  understood as elements of $L^{\infty}[0,\frac{1}{b}]$ acting on $T_{\frac{k}{b}}f$ in the Hilbert $L^{\infty}[0,\frac{1}{b}]$-module $L_{\frac{1}{b}}^{\infty}(\ell^2)$.

In this way we have obtained a general modular form of the modular Walnut representation.

\vspace{.1in}

\begin{theorem}\label{walnut mod}
Lat $G(g,a,b)$ and $G(h,a,b)$ be Bessel sequences in $L^2(\Bbb R)$. Denote by $U$ and $V$ the analysis operators of the weak Bessel sequences $(T_{na}(\frac{1}{\sqrt{b}}g))_{n\in \Bbb Z}$ and $(T_{na}(\frac{1}{\sqrt{b}}h))_{n\in \Bbb Z}$ in the Hilbert $L^{\infty}[0,\frac{1}{b}]$-module $L_{\frac{1}{b}}^{\infty}(\ell^2)$. Then
\begin{equation}\label{walnut mod 1}
\langle V^*Uf,p\rangle_{L_{\frac{1}{b}}^{\infty}(\ell^2)}=(\text{strong})\,\frac{1}{b}\left\langle\sum_{k\in \Bbb Z}G_k^{h,g}T_{\frac{k}{b}}f,p\right \rangle_{L_{\frac{1}{b}}^{\infty}(\ell^2)},\,\,\, \forall f,p\in L_{\frac{1}{b}}^{\infty}(\ell^2),
\end{equation}
and
\begin{equation}\label{walnut mod 2}
\langle V^*Uf,p\rangle_{L_{\frac{1}{b}}^{\infty}(\ell^2)}(x)=\frac{1}{b}\left\langle\sum_{k\in \Bbb Z}G_k^{h,g}T_{\frac{k}{b}}f,p\right\rangle_{L_{\frac{1}{b}}^{\infty}(\ell^2)}\,\,\hspace{-.1in}(x),\,\mbox{for a.e.}\,x\in [0,\frac{1}{b}], \forall f,p\in L_{\frac{1}{b}}^{\infty}(\ell^2).
\end{equation}
\end{theorem}
\proof
The first formula is obtained directly from \eqref{Bessel-adjoint Gabor 2} and  the second formula follows from Remark \ref{strong vs measure}.
\qed

\vspace{.1in}

An an immediate consequence we get the following result concerning pointwise convergence of the Walnut representation for Gabor Bessel sequences in $L^2(\Bbb R)$.

\vspace{.1in}

\begin{cor}\label{walnut tockovno}
Lat $G(g,a,b)$ and $G(h,a,b)$ be Bessel sequences in $L^2(\Bbb R)$ with the analysis operators $U_0$ and $V_0$.  Then
\begin{equation}\label{walnut mod 3}
V_0^*U_0f(x)=\frac{1}{b}\sum_{k\in \Bbb Z}G_k^{h,g}T_{\frac{k}{b}}f(x),\mbox{ for a.e. }x\in \Bbb R,\,\,\, \forall f\in L^2(\Bbb R).
\end{equation}
\end{cor}
\proof
Denote by $U$ and $V$ the analysis operators of the weak Bessel sequences $(T_{na}(\frac{1}{\sqrt{b}}g))_{n\in \Bbb Z}$ and $(T_{na}(\frac{1}{\sqrt{b}}h))_{n\in \Bbb Z}$ in the Hilbert $L^{\infty}[0,\frac{1}{b}]$-module $L_{\frac{1}{b}}^{\infty}(\ell^2)$.
The equality \eqref{walnut mod 2} with $p=\chi_{[\frac{j}{b},\frac{j+1}{b}]}$ for all $j\in \Bbb Z$ gives us
\begin{equation}\label{walnut petljam}
V^*Uf(x)=\frac{1}{b}\sum_{k\in \Bbb Z}G_k^{h,g}T_{\frac{k}{b}}f(x),\mbox{ for a.e. }x\in \Bbb R,\,\,\, \forall f\in L_{\frac{1}{b}}^{\infty}(\ell^2).
\end{equation}
By Lemma \ref{frame operatori isti} we then also have
\begin{equation}\label{walnut petljam 1}
V_0^*U_0f(x)=\frac{1}{b}\sum_{k\in \Bbb Z}G_k^{h,g}T_{\frac{k}{b}}f(x),\mbox{ for a.e. }x\in \Bbb R,\,\,\, \forall f\in L_{\frac{1}{b}}^{\infty}(\ell^2).
\end{equation}
Take any $f\in L^2(\Bbb R)$. As in Remark \ref{klasicno skoro svuda} let
$$
f_0(x)=\left\{
\begin{array}{cl}0,&\mbox{if }\langle f,f\rangle_{\frac{1}{b}}(x)=0\\
\frac{f(x)}{\langle f,f\rangle_{\frac{1}{b}}(x)^{1/2}},&\mbox{if }\langle f,f\rangle_{\frac{1}{b}}(x)\not=0
\end{array}\right..
$$
Obviously, $f_0\in L_{\frac{1}{b}}^{\infty}(\ell^2)$ and \eqref{walnut petljam 1} applied to $f_0$ gives us \eqref{walnut mod 3}.
\qed

\vspace{.3in}

In particular, when we take $h=g$ in the preceding corollary we get a formula for the frame operator $S$.
Recall the classical Walnut representation for Gabor Bessel sequences (\cite{C}, Theorem 9.2.1):
\begin{equation}\label{klasicni walnut}
Sf=\frac{1}{b}\sum_{k\in \Bbb Z}G_kT_{\frac{k}{b}}f,\,\,\,\forall f\in L^2(\Bbb R).
\end{equation}
This is originally proved by D.~Walnut for Bessel Gabor sequences with the generating function $g$ from the Wiener space $W(L^{\infty},\ell^1)$ and then extended to Bessel Gabor sequences whose generators $g$ satisfy \eqref{Gabor gs 1} (i.e. the CC condition).
The series in \eqref{klasicni walnut} converges absolutely when $g\in W(L^{\infty},\ell^1)$ and unconditionally in case that $g$ satisfies \eqref{Gabor gs 1}.
A natural question arises whether the Walnut representation in the modular setting can be obtained with some stronger form of convergence than in \eqref{walnut mod 1}.

It turns out that the answer is positive for those weak Bessel sequences in $L_{\frac{1}{b}}^{\infty}(\ell^2)$ that are standard in $X_{\frac{1}{b}}(L^{\infty},\ell^2)$.

\vspace{.1in}

\begin{prop}\label{bolji walnut}
Lat $G(g,a,b)$ and $G(h,a,b)$ be Bessel sequences in $L^2(\Bbb R)$. Denote by $U$ and $V$ the analysis operators of the weak Bessel sequences $(T_{na}(\frac{1}{\sqrt{b}}g))_{n\in \Bbb Z}$ and $(T_{na}(\frac{1}{\sqrt{b}}h))_{n\in \Bbb Z}$. Suppose that $(T_{na}(\frac{1}{\sqrt{b}}g))_{n\in \Bbb Z}$ is a standard Bessel sequence in $X_{\frac{1}{b}}(L^{\infty},\ell^2)$. Then
\begin{equation}\label{bolji walnut jedn 1}
V^*Uf=\frac{1}{b}\sum_{n\in \Bbb Z}\langle f,T_{na}g\rangle_{L_{\frac{1}{b}}^{\infty}(\ell^2)}T_{na}h,\,\,\,\forall f\in X_{\frac{1}{b}}(L^{\infty},\ell^2),
\end{equation}
and
\begin{equation}\label{bolji walnut jedn 2}
V^*Uf=\frac{1}{b}\sum_{k\in \Bbb Z}G_k^{h,g}T_{\frac{k}{b}}f,\,\,\,\forall f\in X_{\frac{1}{b}}(L^{\infty},\ell^2),
\end{equation}
where both series converge unconditionally in norm in $L_{\frac{1}{b}}^{\infty}(\ell^2)$.
\end{prop}
\proof
Since $(T_{na}(\frac{1}{\sqrt{b}}g))_{n\in \Bbb Z}$ is a standard Bessel sequence in $X_{\frac{1}{b}}(L^{\infty},\ell^2)$, we know that
$$
Uf=\left(\langle f,\frac{1}{\sqrt{b}}T_{na}g\rangle_{L_{\frac{1}{b}}^{\infty}(\ell^2)}\right)_n\in \ell^2_{\Bbb Z}(L^{\infty}[0,\frac{1}{b}]),\,\,\,\forall f\in X_{\frac{1}{b}}(L^{\infty},\ell^2).
$$
Using Remark \ref{r2} we conclude  \eqref{bolji walnut jedn 1}. Since the series in \eqref{bolji walnut jedn 1} converges unconditionally, we can write
$$
V^*Uf=\frac{1}{b}\sum_{k\in \Bbb Z}\left(\sum_{n\in \Bbb Z}T_{na}hT_{\frac{k}{b}}T_{na}\overline{g}\right)T_{\frac{k}{b}}f,\,\,\,\forall f\in X_{\frac{1}{b}}(L^{\infty},\ell^2)
$$
which is in fact the desired equality \eqref{bolji walnut jedn 2}.
\qed

\vspace{.1in}

As a corollary we now have the Walnut representation of the frame operator of Gabor Bessel sequences in $L^2(\Bbb R)$ generated by windows which need not belong to the Wiener space nor necessarily satisfy the CC condition.

\vspace{.1in}

\begin{theorem}\label{mozda novi walnut}
Let $G(g,a,b)$ be a Bessel sequence in $L^2(\Bbb R)$, $g\in L^2(\Bbb R)$, $a,b>0$, with the frame operator $S$. If $g\in X_{\frac{1}{b}}(L^{\infty},\ell^2)\cap X_{a}(L^{\infty},\ell^2)$ then
\begin{equation}\label{zadnja tvrdnja}
Sf=\frac{1}{b}\sum_{n\in \Bbb Z}\langle f,T_{na}g\rangle_{L_{\frac{1}{b}}^{\infty}(\ell^2)}T_{na}g,\,\,\,\forall f\in X_{\frac{1}{b}}(L^{\infty},\ell^2)
\end{equation}
where this series converges unconditionally in norm in $L^2(\Bbb R)$.
\end{theorem}
\proof
Denote by $U$ the analysis operator of the corresponding weak Bessel sequence $(T_{na}(\frac{1}{\sqrt{b}}g))_{n\in \Bbb Z}$ and by $U_0$ the analysis operator of $G(g,a,b)$. Notice that $S=U_0^*U_0$. Since by Lemma \ref{frame operatori isti} $U^*U$ and $U_0^*U_0$  coincide on $L_{\frac{1}{b}}^{\infty}(\ell^2)$ and since $\|f\|_{L^2(\Bbb R)}\leq\frac{1}{b}\|f\|_{L_{\frac{1}{b}}^{\infty}(\ell^2)}$ for all $f\in L_{\frac{1}{b}}^{\infty}(\ell^2)$, it is enough to see that
$$
U^*Uf=\frac{1}{b}\sum_{n\in \Bbb Z}\langle f,T_{na}g\rangle_{L_{\frac{1}{b}}^{\infty}(\ell^2)}T_{na}g,\,\,\,\forall f\in X_{\frac{1}{b}}(L^{\infty},\ell^2)
$$
in norm in $L_{\frac{1}{b}}^{\infty}(\ell^2)$. This in turn will follow from
\eqref{bolji walnut jedn 1} in Proposition \ref{bolji walnut} if we can prove that $(T_{na}(\frac{1}{\sqrt{b}}g))_{n\in \Bbb Z}$ is a standard Bessel sequence in $X_{\frac{1}{b}}(L^{\infty},\ell^2)$.

To end this, first observe that the assumption $g\in X_{\frac{1}{b}}(L^{\infty},\ell^2)$ implies by Proposition \ref{translation invariant} that the whole sequence $(T_{na}(\frac{1}{\sqrt{b}}g))_{n\in \Bbb Z}$ belongs to $X_{\frac{1}{b}}(L^{\infty},\ell^2)$.
By proposition \ref{kad je weak dosao od standardnog} it suffices to show that $\left( \langle f,T_{na}g\rangle_{L_{\frac{1}{b}}^{\infty}(\ell^2)}\right)_n\in \ell^2_{\Bbb Z}(L^{\infty}[0,\frac{1}{b}])$ for all $f$ from a dense subset of
$X_{\frac{1}{b}}(L^{\infty},\ell^2)$.

We shall show that
$$
\left( \langle \chi_{[\frac{j}{b},\frac{j+1}{b}]},T_{na}g\rangle_{L_{\frac{1}{b}}^{\infty}(\ell^2)}\right)_n\in \ell^2_{\Bbb Z}(L^{\infty}[0,\frac{1}{b}]),\,\,\,\forall j\in \Bbb Z.
$$
Take any $j\in \Bbb Z$. We must show that the series
$$
\sum_{n\in\Bbb Z}|\langle T_{na}g,\chi_{[\frac{j}{b},\frac{j+1}{b}]}\rangle_{L_{\frac{1}{b}}^{\infty}(\ell^2)}|^2
$$
converges in norm in $L^{\infty}[0,\frac{1}{b}]$.
Using \eqref{matricni koeficijenti}, the above series can be written as
\begin{equation}\label{ciljna konvergencija}
\sum_{n\in \Bbb Z}|g(x-na+\frac{j}{b})|^2.
\end{equation}
This is a series of positive elements in a $C^*$-algebra; thus, it is enough to consider symmetric partial sums.

Now we use the other part of our assumption; namely that $g\in X_{a}(L^{\infty},\ell^2)$. Again by Proposition \ref{translation invariant} we have $T_{-\frac{j}{b}}g\in X_{a}(L^{\infty},\ell^2)$.
Thus, the series
$$
\langle T_{-\frac{j}{b}}g,T_{-\frac{j}{b}}g \rangle_{L_{a}^{\infty}(\ell^2)}=\sum_{n\in \Bbb Z}|g(x-na+\frac{j}{b})|^2
$$
converges in norm in $L^{\infty}[0,a]$.
In particular, the sequence of partial sums is Cauchy. For $\epsilon>0$ there exists $N_0$ such that
$$
M>N\geq N_0 \Rightarrow \left\|\sum_{n=-M}^M|g(x-na+\frac{j}{b})|^2-\sum_{n=-N}^N|g(x-na+\frac{j}{b})|^2 \right\|_{L^{\infty}[0,a]}<\epsilon
$$
wherefrom
\begin{equation}\label{essup krajnji}
\esssup_{x\in [0,a]} \sum_{|n|\geq N_0}|g(x-na+\frac{j}{b})|^2 \leq \epsilon.
\end{equation}
If $\frac{1}{b}<a$ this obviously implies
$$
\esssup_{x\in [0,\frac{1}{b}]} \sum_{|n|\geq N_0}|g(x-na+\frac{j}{b})|^2 \leq \epsilon.
$$
If $a<\frac{1}{b}$ let $N_1$ be the greatest natural number for which $N_1a<b$. Then for every $x\in [0,\frac{1}{b}]$ there exists a natural number $N(x)\in \{0,1,\ldots,N_1\}$ such that $x-N(x)a\in [0,a]$. This allows us to conclude from \eqref{essup krajnji} that
$$
\esssup_{x\in [0,\frac{1}{b}]} \sum_{|n|\geq N_0+N_1}|g(x-na+\frac{j}{b})|^2 \leq \epsilon.
$$
Thus, the sequence of partial sums of the series \eqref{ciljna konvergencija} is Cauchy in $L^{\infty}[0,\frac{1}{b}]$.
\qed

\vspace{.1in}

\begin{remark}\label{zavrsni}
Observe that \eqref{zadnja tvrdnja} extends by Corollary \ref{walnut tockovno} to all functions from $L^2(\Bbb R)$, but only in the sense of pointwise a.e. convergence for those $f$ that do not belong to $X_{\frac{1}{b}}(L^{\infty},\ell^2)$.

The condition $g\in X_{\frac{1}{b}}(L^{\infty},\ell^2)\cap X_{a}(L^{\infty},\ell^2)$ from the preceding theorem is satisfied for all $g\in W(L^{\infty},\ell^1)$ - this follows immediately from Remark \ref{4 prostora}. Also, if $g\in X_{\frac{1}{b}}(L^{\infty},\ell^2)$ or $g\in X_{a}(L^{\infty},\ell^2)$ with the additional assumption $ab\in \Bbb Q$, then
Proposition \ref{racionalno isti} implies $g\in X_{\frac{1}{b}}(L^{\infty},\ell^2)\cap X_{a}(L^{\infty},\ell^2)$.
\end{remark}

\vspace{.3in}

{\it Conclusion.} In the first part of the paper the dual of the standard Hilbert $C^*$-module over an arbitrary $C^*$-algebra is described. In particular, a concrete description is obtained when the underlying algebra is a von Neumann algebra. As pointed out by M.~Frank in private communication, some of the results concerning modules over for von Neumann algebras can be extended at least to the class of Hilbert $C^*$-modules over monotone complete $C^*$-algebras. This is left for future investigation.

We have introduced a concept of a weak Bessel sequence and a weak frame in Hilbert $C^*$-modules over von Neumann algebras. Fundamental properties of such systems are obtained. It turned out that such weak modular systems behave similar to standard Bessel sequences and frames with respect to certain weak topology. Moreover, if the underlying von Neumann algebra is commutative, this weak modular systems are naturally described and represented by their Gram matrices.

Weak modular Bessel sequences and frames naturally appear in Gabor analysis. In fact,
standard Gabor Bessel sequences and Gabor frames in $L^2(\Bbb R)$ may be interpreted as weak Bessel sequences resp.~weak frames of translates in certain Hilbert $C^*$-module over the commutative von Neumann algebra $L^{\infty}(I)$ where $I=[0,\frac{1}{b}]\subseteq \Bbb R$ is the interval  determined by the modulation parameter $b$. This correspondence enabled us to reobtain (and to reinterpret) in a natural way some of the classical results from Gabor analysis. Some of the results, e.g.~that which is concerned with the Walnut representation, appear to broaden the scope of the corresponding classical results on Gabor systems. Certainly, this line of investigation owns very much to the approach of P.~Casazza and M.C.~Lammers and uses the ideas, but in a different language, from the work of A.~Ron and Z.~Shen. However, simplicity of this new proofs of some of the classical results (whose original proofs are very involved), suggests that our Hilbert $C^*$-module technique might serve as a promising tool in the study of Gabor systems.

\vspace{.2in}

{\it Acknowledgement.} The author is indebted to Michael Frank for careful reading the whole manuscript, useful discussions and several helpful remarks that improved the original exposition.

\end{document}